\input amstex
\documentstyle{amsppt}
\loadmsbm
\loadbold
\vsize=7.5in
\UseAMSsymbols
\magnification=\magstep1
\font\rom cmr10

\rom
\nopagenumbers

%
%
%
\input epsf  

%
%
%
%
%
\ifx\curlyloaded\relax\let\next 
  \else\let\curlyloaded\relax\let\next\relax
\fi
\next
%
\xdef\recoveratcodezqrz{\catcode`\noexpand\@=\the\catcode`\@}
\catcode`\@=11
\edef\innernewfam{\expandafter\noexpand\csname newfam\endcsname}
%
\ifx\amstexloaded@\relax
  \loadmathfont{rsfs}
  \def\curly{\mathfont@\curly}
  
\else
  \font\tenrsfs=rsfs10
  \font\sevenrsfs=rsfs7
  \font\fiversfs=rsfs5
  \innernewfam\rsfsfam \def\curly{\fam\rsfsfam\tenrsfs}
  \textfont\rsfsfam=\tenrsfs  \scriptfont\rsfsfam=\sevenrsfs 
  \scriptscriptfont\rsfsfam=\fiversfs
  \def\curly{\fam\rsfsfam}
\fi
\let\innernewfam\undefined
\recoveratcodezqrz
 %

\def\({\left(}
\def\){\right)}
\def\pfrac#1#2{\(\frac{#1}{#2}\)}
\def\kk{\frac12 k(k+1)}

\def\xx{\boldkey x}
\def\zz{\boldkey z}
\def\yy{\boldkey y}
\def\hh{\boldkey h}

\def\mm{\boldkey m}
\def\nn{\boldkey n}
\def\uu{\boldkey u}
\def\vv{\boldkey v}
\def\ww{\boldkey w}
\def\cc{\boldkey c}
\def\bc{\boldkey c}
\def\bb{\boldkey b}
\def\dd{\boldkey d}
\def\ff{\boldkey f}

\redefine\AA{\curly A}
\def\BB{\curly B}
\def\DD{\curly D}
\def\CC{\curly C}
\def\FT{\widetilde{F}}

\def\a{\alpha}
\def\eps{\varepsilon}
\def\e{\varepsilon}
\def\g{\gamma}
\redefine\d{\delta}
\def\th{\theta}
\def\om{\omega}

\def\lm{\lambda}
\def\lam{\lambda}
\def\del{\delta}
\def\delbar{\bar{\delta}}

\define\balph{\boldsymbol\alpha}
\def\bbeta{\boldsymbol \beta}
\define\baa{\boldsymbol\alpha}
\define\bpsi{\boldsymbol\Psi}
\define\bphi{\boldsymbol\Phi}

\define\bups{\boldsymbol\Upsilon}
\define\Ho{H\"older's inequality}
\define\CS{the Cauchy-Schwarz inequality }
\def\CPR{\CC(P,R)}

\define\fl#1{\lfloor #1 \rfloor}
\define\bigfl#1{\left\lfloor #1 \right\rfloor}

%
%

\topmatter

\title Vinogradov's Integral and bounds for the Riemann
Zeta Function \endtitle
\author Kevin Ford \endauthor
\rightheadtext{Vinogradov's Integral and Riemann Zeta Function}
\thanks
Research supported in part by National Science Foundation grant 
DMS--0070618 
\endthanks
\email ford\@math.sc.edu \endemail
\date October 25, 2001 \enddate
\address Department of Mathematics, University of Illinois at Urbana-Champaign,
Urbana, IL 61801
\endaddress
\subjclass\nofrills{2000 {\it Mathematics Subject Classification.}}
Primary 11M06, 11N05, 11L15; Secondary 11D72, 11M35
\endsubjclass
\endtopmatter

\document
\head 1. Introduction \endhead


The methods of Korobov [11] and Vinogradov [28]  produce a zero-free region
for the Riemann zeta function $\zeta(s)$ of the following strength: for
some $c>0$, there
are no zeros of $\zeta(s)$ for $s=\sigma+it$ with $|t|\ge 3$ and
$\sigma>1-c(\log |t|)^{-2/3} (\log\log |t|)^{-1/3}$.  The principal tool
is an upper bound for $|\zeta(s)|$ near the line $\sigma=1$.  In 1967,
Richert [22] used this method to give 
the bound 
$$
|\zeta(\sigma+it)| \le A |t|^{B(1-\sigma)^{3/2}} \log^{2/3} |t|
\qquad (|t|\ge 2, \tfrac12 \le \sigma \le 1) \tag{1.1}
$$
with $B=100$ and $A$ and unspecified absolute constant.
Similar results with smaller $B$ values
have been proven subsequently by several authors, the best being $B=18.4974$
 and due to Kulas [13].  Recently,
Y. Cheng [3] has given a completely explicit version of this bound,
with $A=175$ and $B=46$. 

In this paper, we improve substantially the value of $B$, while also keeping
the bound entirely explicit.  More generally, we bound the Hurwitz zeta
function, defined for $\Re s > 1$ and $0<u\le 1$ by
$\zeta(s,u)= \sum_{n=0}^\infty (n+u)^{-s}$. 
The Hurwitz zeta function may be used to bound Dirichlet L-functions via
the identity
$L(s,\chi) = q^{-s} \sum_{m=1}^{q} \chi(m) \zeta(s,m/q)$,
where $\chi$ is a Dirichlet character modulo $q$. 
Notice that $\zeta(s)=\zeta(s,1)$.   Since $\zeta(\bar{s},u)=
\overline{\zeta(s,u)}$, we may restrict our attention to $s$ lying in the
upper half-plane.

\proclaim{Theorem 1} The inequalities
$$
\split
| \zeta(\sigma+it)| &\le A t^{B(1-\sigma)^{3/2}} \log^{2/3} t
 \qquad (t\ge 3, \tfrac12 \le \sigma \le 1), \\
| \zeta(\sigma+it,u)-u^{-s}| &\le A t^{B(1-\sigma)^{3/2}} \log^{2/3} t
 \qquad (0<u\le 1, t\ge 3, \tfrac12 \le \sigma \le 1)
\endsplit
$$
hold with $B=4.45$ and $A=76.2$.
\endproclaim

If the Riemann Hypothesis is true, then the conclusion of Theorem 1 holds
with any positive $B$, with the constant $A$ depending on $B$.
Bounds of the type (1.1) with explicit values of $B$
have numerous applications, including (i) explicit zero-free regions
for $\zeta(s)$; (ii) explicit error bounds for the prime number theorem;
(iii) zero density bounds for $\zeta(s)$; (iv) mean value theorems for
$\zeta(s)$; (v) bounds for error terms in the Dirichlet divisor problem.
We briefly indicate the consequences of Theorem 1 for each of these five
problems.

(i) One can use (1.1) to give explicit values for the constant $c$ in the
zero-free region mentioned in the opening paragraph.  In a separate paper
[6], the author shows that $\zeta(\beta+it)\ne 0$ for $t$ sufficiently
large and
$$
1-\beta \le \frac{0.05507 B^{-2/3}}{(\log t)^{2/3}(\log\log t)^{1/3}}.
$$
Moreover, using the full strength of Theorem 1, in [6] the zero-free
region
$$
t\ge 3, \quad 1-\beta \le \frac{c}{(\log t)^{2/3}(\log\log t)^{1/3}},
\quad c=\frac{1}{57.54}
$$
is proved.  By comparison, Popov [20] showed that the above holds with
holds with $c=0.00006888$, and Cheng [4] proved a zero-free region with
$c=1/990$.

(ii) A corollary of Theorem 1, the work in [6], and Theorem 8 of Pintz [19],
is the following error bound in the prime number theorem:
$$
\pi(x) - \text{li}(x) = O \( x \exp \{ - c (\log x)^{3/5} (\log\log x)^{-1/5}
\} \), \quad c=0.2098.
$$

(iii) Let $N(\sigma,T)$ denote the number of zeros of
$\zeta(s)$ in the rectangle $\sigma \le \Re s \le 1$, $|\Im s|\le T$.
If (1.1) holds, then for $\frac9{10} \le \sigma \le 1$, we have
$$
N(\sigma,T) \ll T^{13.043B (1-\sigma)^{3/2}} \log^{15} T.
$$
This follows from Theorem 12.3 of Montgomery [17], taking
$1-\alpha = 4.93(1-\sigma)$; see also \S 11.4 of [8]. 
Incidentally, there is an error in Corollary 12.5 of [17], where it
is stated that $B=100$ implies
$$
N(\sigma,T)\ll T^{167(1-\sigma)^{3/2}} \log^{17} T.
$$
As a corollary, Theorem 1 gives
$$
N(\sigma,T) \ll T^{58.05 (1-\sigma)^{3/2}} \log^{15} T.
$$

(iv) Let 
$$
M_k(\sigma,T) = \frac{1}{T} \int_0^T |\zeta(\sigma+it)|^{2k}\, dt.
$$
Let $\sigma_k$ be the infimum of the numbers $\sigma$ with $M_k(\sigma,T)
=O(1)$, and let $\mu_k(\sigma)$ be the infimum of the numbers $\xi$ such
that $M_k(\sigma,T)=O(T^\xi)$.  If $\sigma>\sigma_k$, we have an
asymptotic formula for $M_k(\sigma,T)$ ([25], \S 7.8):
$$
M_k(\sigma,T) \sim \sum_{n=1}^\infty \frac{d_k(n)^2}{n^{2\sigma}},
$$
where $d_k(n)$ is the number of $k$-tuples of positive integers
$(b_1,b_2,\cdots,b_k)$ with $b_1\cdots b_k = n$.  In particular, $d_2(n)$
is the number of positive divisors of $n$.  Also, when $\Re s > 1$,
$(\zeta(s))^k = \sum_{n=1}^\infty d_k(n) n^{-s}$. 
Upper bounds on $\sigma_k$ can be deduced from upper bounds on $\zeta(s)$
inside the critical strip by means of a Theorem of Carlson ([25], Theorem
7.9): for any $0<\a<1$, we have
$$
\sigma_k \le \max \( \frac12, \a, 1 - \frac{1-\a}{1+\mu_k(\a)} \). \tag{1.2}
$$
By (1.1), we have trivially $\mu_k(\sigma) \le 2Bk(1-\sigma)^{3/2}$.  Taking
$\a = 1 - (Bk)^{-2/3}$ in (1.2) gives $\sigma_k \le 1 - \frac13 (Bk)^{-2/3}$.
For more on mean value theorems, see Chapter VII of [25] and Chapter 8 of
[8].

(v)  Denote by $\Delta_k(x)$
the usual error term in the Dirichlet divisor problem, i.e.
$$
\Delta_k(x) = \sum_{n\le x} d_k(n) - \operatornamewithlimits{Res}_{s=1}
 x^s (\zeta(s))^k  s^{-1} =  \sum_{n\le x} d_k(n) -xP_k(\log x),
$$
where $P_k$ is a certain polynomial.  Let $\a_k$ be the infimum of
numbers $\a$ with $\Delta_k(x)=O(x^\a)$.  Dirichlet in 1849 proved that
$\a_2 \le \frac12$ and his method can be used to deduce
$\a_k \le 1 - \frac{1}{k}$.  Modern treatments make use of Perron's 
formula in the form
$$
 \sum_{n\le x} d_k(n) = \frac{1}{2\pi i} \int_{c-i\infty}^{c+i\infty} 
\zeta^k(s) \frac{x^s}{s}\, ds, \quad c>1.
$$
Then the contour is moved inside the critical strip, the main term coming
from the pole at $s=1$, and the error term coming from upper bounds for
$\zeta(s)$. In 1960, Richert [21] proved that $\a_k \le 1 - ck^{-2/3}$
for some positive constant $c$.
Subsequently, the value of $c$ was made explicit as a function of the
constant $B$ in (1.1) by Karatsuba [10] ($c=\frac12 (2B)^{-2/3} \approx
0.31498B^{-2/3}$).  Writing $c=dB^{-2/3}$, the value of $d$ was improved by
Ivi\'c and Ouellet [9] to $d=\frac13 2^{2/3} \approx 0.52913$.
There are two claims for larger $d$, but both arguments are flawed.
Fujii [7] claims $d=2^{-1/2} (\sqrt{8}-1)^{-1/3} \approx 0.57826$, but
the details are omitted (the method appears to give $d=\frac12$);
Panteleeva [18] claims $d=2^{-2/3} \approx
0.62996$, but the proof of this result (Theorem 3 of [18]) has a flaw,
namely the differentiation of (14) in invalid.

For the mean square of $\Delta_k(x)$,
Ivi\'c and Ouellet [9] proved that
$$
\int_1^x \Delta^2_k(y)\, dy  \ll_{\varepsilon,k} x^{1+2b_k+\varepsilon},
\qquad b_k =  1 - \frac23 \pfrac{1}{Bk}^{2/3}.
$$
More information may be found in Chapter XII of [25] and Chapter 13 of [8].
\bigskip

Theorem 1 depends primarily
on upper bounds for the following exponential sum:
$$
S(N,t) = \max_{0<u\le 1} \max_{N < R \le 2N} \left| \sum_{N<n\le R} 
(n+u)^{-it} \right|,
$$
where $N$ is a positive integer and $t\ge N$.  We shall prove the following.

\proclaim{Theorem 2} Suppose $N$ is a positive integer,
$N \le t$ and set $\lambda = \frac{\log t}{\log N}$.  Then
$$
S(N,t) \le 9.463 N^{1-1/(133.66 \lambda^2)}.
$$
\endproclaim

By comparison, Kulas [12] proved that $S(N,t) \ll N^{1-1/(2309.525\lambda^2)}$
for $\lambda\ge 1000$.

\proclaim{Corollary 2A} Suppose $\chi$ is a Dirichlet character modulo
$q$, where $q\le N$ and $2\le N \le qt$.  Then
$$
\max_{N<R\le 2N} \left| \sum_{N<n\le R} \chi(n) n^{-it} \right| \le
10.463 \frac{\phi(q)}{q} N e^{-\frac{\log^3 (N/q)}{133.66 \log^2 t}}.
$$
\endproclaim

\demo{Proof}  Suppose the maximum on the left occurs at $R=R_0$.  Then
$$
\sum_{N<n\le R_0} \chi(n) n^{-it} = \sum_{\ell=1 \atop (\ell,q)=1}^{q}
\chi(\ell) \sum_{N<n\le R_0 \atop n\equiv \ell \pmod{q}} n^{-it}.
$$
Writing $n=mq+\ell$  gives
$$
\split
\left| \sum_{N<n\le R_0 \atop n\equiv \ell \pmod{q}} n^{-it} \right|
 &\le 1 + \left| \sum_{\frac{N-\ell+q}{q} < m \le \frac{R_0-\ell}{q}}
(m+\ell/q)^{-it} \right| \\
&\le 1 + S\( \frac{N-\ell+q}{q}, t \).
\endsplit
$$
Theorem 2 then gives
$$
\left| \sum_{N<n\le R_0} \chi(n) n^{-it} \right| \le \phi(q) \( 1 +
9.463 \pfrac{N}{q}^{1 - \frac{\log^2(N/q)}{133.66 \log^2 t}} \).
$$
Lastly, $N/q \ge 1$, and the result follows.
\qed\enddemo

As with prior treatments, 
Theorem 2 in turn depends on explicit bounds for Vinogradov's integral,
defined as
$$
J_{s,k}(P) = \int_{[0,1]^k}
\left| \sum_{1\le x\le P} e(\a_1 x + \cdots + \a_k x^k)\right|^{2s} 
d{\boldsymbol \alpha}, \tag{1.3}
$$
where $\baa=(\a_1,\ldots,\a_k)$ and $e(z)=e^{2\pi i z}$.
Equivalently, $J_{s,k}(P)$ is the number of solutions of the simultaneous
equations
$$
\sum_{i=1}^s (x_i^j - y_i^j) = 0 \qquad (1\le j\le k); \quad 1\le x_i,y_i
\le P. \tag{1.4}
$$
For $\hh=(h_1,\ldots,h_k)$, let $J_{s,k}(P;\hh)$ be the number of
solutions of
$$
\sum_{i=1}^s (x_i^j-y_i^j)=h_j \quad (1\le j\le k); \quad 1\le x_i,y_i\le P.
$$
In particular,
$$
\split
J_{s,k}(P;\hh) &=  \int_{[0,1]^k}
\left| \sum_{1\le x\le P} e(\a_1 x + \cdots + \a_k x^k)\right|^{2s} 
e(-\a_1 h_1 - \cdots - \a_k h_k) d{\boldsymbol \alpha} \\
&\le J_{s,k}(P;(0,\ldots,0)) = J_{s,k}(P).
\endsplit
$$
Hence, writing $Q=\lfloor P \rfloor$, we obtain
$$
Q^{2s}=\sum_{\hh} J_{s,k}(P;\hh)\le \sum_{\Sb \hh \\ |h_j| \le s(Q^j-1) 
\endSb} J_{s,k}(P) \le (2s)^k Q^{k(k+1)/2} J_{s,k}(P).
$$
Also, counting only the solutions of (1.4) with $x_i=y_i$ for each $i$
gives $J_{s,k}(P) \ge Q^s$.  Therefore
$$
J_{s,k}(P) \ge \max\( (2s)^{-k} \lfloor P\rfloor ^{2s-\kk}, \lfloor P \rfloor
^s \). \tag{1.5}
$$
Upper bounds take the form of
$$
J_{s,k}(P) \le D(s,k) P^{2s - \kk + \eta(s,k)}, \tag{1.6}
$$
where $\eta(s,k) \ge 0$ and $D(s,k)$ is independent of $P$.
Stechkin in 1975 [24] proved (1.6) with
$$
\eta(rk,k) = \frac12 k^2 (1-1/k)^r, \qquad
D(rk,k)=\exp\{C\min(r,k)k^2\log k\}
$$
for an absolute constant $C$.  The constant factor was improved by Wooley [31].
Small improvements to the exponents of $P$ were subsequently made by 
Arkhipov and Karatsuba [1] and Tyrina [26] (significant for $s\ll k^2$).
Also significant is Wooley's [32] result when 
$s\ll k^{3/2-\eps}$, which is very close to the ``ideal'' bounds
$C(k,s)P^s$ in that range of $s$.  For our purposes, the most important 
improvement comes from
Wooley [30], who improved the exponents substantially in a wide range of $s$,
showing that (1.6) holds with $\eta(k,s) \approx \frac12 k^2
e^{1/2-2s/k^2}$ valid for $s\ll k^2\log k$ (see [5], Lemma 5.2).
In Theorem 3 below, we combine Wooley's method with the main idea from
[1] to
improve this to $\eta(k,s) \approx \frac38 k^2 e^{1/2-2s/k^2}$.
In the application to
bounding the Riemann zeta function, we will take $s$ to be of order
$k^2$, so this small improvement is significant.

\proclaim{Theorem 3} Let $k$ and $s$ be integers with $k\ge 1000$ and
$2k^2 \le s \le \frac{k^2}{2} ( \frac12 + \log \frac{3k}{8} )$. 
Then
$$
J_{s,k}(P) \le  k^{2.055k^3-5.91k^2+3s} 1.06^{sk+2s^2/k-9.7278k^3}
 P^{2s-\frac12 k(k+1) + \Delta_s} \quad (P\ge 1),
$$
where
$$
\Delta_s = \tfrac38 k^2 e^{1/2 - 2s/k^2+1.7/k}.
$$
Further, if $k\ge 129$, there is an integer $s\le \rho k^2$
such that for $P\ge 1$,
$$
J_{s,k}(P) \le k^{\theta k^3} P^{2s-\kk + 0.001k^2},
$$
with
$$
(\rho,\theta) = \cases (3.21432, 2.3291) & (k\ge 200) \\
(3.21734, 2.3849) & (150\le k\le 199) \\
(3.22313, 2.4183) & (129\le k\le 149) \endcases \tag{1.7}
$$
\endproclaim
By itself, Theorem 3 implies the inequalities in
Theorem 1 with $B$ a bit more than 10.4.

The most significant new idea is to bound $S(N,t)$ in terms of both
$J_{s,k}(P)$ and another quantity which counts the number of solutions of
{\it incomplete} Diophantine systems (where we regard (1.4) to be 
{\it complete} because the powers of the variables range from 1 to $k$).
Define $J_{s,k,h}(\BB)$ to be the number of solutions of the system
$$
\sum_{i=1}^s (x_i^j-y_i^j) = 0 \qquad (h\le j\le k); \qquad
x_i,y_i \in \BB. \tag{1.8}
$$
Incomplete systems were first studied my Mardzhanishvili ([15], [16]),
who gave sufficient conditions for the existence of solutions of the
system 
$$
\sum_{i=1}^s x_i^j = N_j \qquad (j\in \curly J),
$$
where $\curly J$ is an arbitrary finite subset of positive integers.
More general systems of Diophantine equations and associated trigonometric
sums are treated in [2].

The Vinogradov method [28], when applied to bounding a more general sum
$$
\sum_{N<n\le 2N} e(p(n)), \quad p(n)= \a_1 n + \cdots + \a_k n^k,
$$
ultimately depends on having good rational approximations for
a subset of the coefficients of $p(n)$, say for
$\a_i,\a_{i+1},\ldots, \a_{j}$.  By applying
trivial estimates to sums involving the other coefficients, we may 
restrict attention to associated mean-values over 
$\a_i,\a_{i+1},\ldots, \a_{j}$ which are equivalent to $J_{s,j,i}(\BB)$.
The core of the argument is given in Lemma 5.1.

When $\BB \subseteq [1,P]$, we have a trivial bound
$$
J_{s,k,h}(\BB) \le s^{h-1} P^{h(h-1)/2} J_{s,k}(P). \tag{1.9}
$$
In the application to bounding $S(N,t)$, however, (1.9) gives nothing better
than if $J_{s,k,h}(\BB)$ were replaced by $J_{s,k}(P)$ from the outset.
By a more sophisticated method, which is a generalization of the author's
work ([5]) on mean values of complete Weyl sums, one can bound
$J_{s,k,h}([1,P])$ in terms of $J_{s',k}(P)$ (with $s'<s$), and attain
superior bounds for $S(N,t)$.  When $\BB = \AA(P,R)$, the set of numbers
$\le P$ with no prime factors exceeding $R$ ($R$-``smooth'' numbers),
 $R$ is a sufficiently small
power of $P$ (depending on $k,h,s$), and $h$
close to $k$, Wooley's ``efficient
differencing'' method ([29], [30], [34]) produces even better
exponents of $P$.  However,  the implied constants coming from
the bounds in [34] grow too fast as functions
of $k,h,s$, and thus are
inadequate for bounding $S(N,t)$ for the entire range
$1\le \lambda\ll \sqrt{\log N}$.
The principal problem is that elements of $\AA(P,R)$ may contain
a very large number of divisors.  We overcome this by taking $\BB=\CPR$,
the set of integers $\le P$ composed only of prime factors in
$(\sqrt{R},R]$.  We thus retain all of the advantages gained by using
$R$-smooth numbers, but now the number of prime factors of each
such number is bounded above by $2\frac{\log P}{\log R}$.  The next
theorem, which will be used for the proof of Theorem 2, is an example of
what can be proved.

\proclaim{Theorem 4}
Suppose $k\ge 60$, $0.9k \le h\le k-2$,
$2t \le s \le \lfloor h/2 \rfloor t$, and
$P\ge e^{Dk^2}$ where $D \ge 10$.
Further assume that 
$$
\frac{2}{k^3} < \eta \le \frac{1}{2k}, \quad
\frac{18}{k} \le \frac{4\log k}{Dk^2 \eta} \le 0.4. \tag{1.10}
$$
Then
$$
J_{s,k,h}(\CC(P,P^\eta)) \le e^C P^{2s-\frac{t}{2}(h+k)+\frac{t(t-1)}{2}
+\eta s^2/(2t) + ht\exp \{ -s/(ht) \}},
$$
where
$$
C = \frac{s^2}{t} + \frac{10.5t \log^2 k}{Dk \eta^2} - s \( \( \frac{1}{\eta}
+h \) \( 1 - \frac{1}{h} \)^{s/t} - h \) \log\pfrac{1}{10\eta}.
$$
\endproclaim

Sections 2, 3 and 4 are dedicated to proving explicit bounds for
$J_{s,k}(P)$ (Theorem 3) and $J_{s,k,h}(\CPR)$ (Theorem 4).
In \S 5, we use Vinogradov's
method and Theorems 3 and 4 to prove Theorem 2 for large $\lambda$.
For smaller $\lambda$ we use older methods (\S 6), which give better results.
  This is then applied to the
problem of bounding $|\zeta(s)|$ and $|\zeta(s,u)|$ in \S 7, where Theorem 1
is proved.  Lastly, in \S 8 we discuss the limit of our method,
and briefly indicate some ways in which the
constant $B$ may be improved a little.  

{\bf Acknowledgements}  The author wishes to thank the following people:
 Y. Cheng for several reprints
and preprints of his work, and for helpful discussions concerning the
proof of Theorem 2 for small $\lambda$; A. Ivi\'c for helpful discussions
concering the applications (iv) and (v) above;
D. Meade for help with Maple code;  K. Oskolkov for help with the
Fourier analysis connected with the functions $\ell(x;w)$ in \S 5.


%
%
%
 
\head 2. Preliminary Lemmata. \endhead
\define\iuk{\int_{\Bbb U^k}}

First, we detail some notational conventions.  Let $\Bbb U = [0,1]$, let
$\lfloor x \rfloor$ be the greatest integer $\le x$, let $\lceil x \rceil$
be the smallest integer $\ge x$, write
$e(z)$ for $e^{2\pi i z}$ and let $\| x \|$ be the
distance from $x$ to the nearest integer.  Let $\CPR$ be the set of 
positive integers $n\le P$, all of whose prime factors are in $(\sqrt{R},R]$.
The functions $\omega(n)$ is the
number of distinct prime factors of $n$, $\Omega(n)$ is the number
of prime power divisors of $n$, $\tau(n)$ is the number of
positive divisors of $n$, and $s_0(n)$ is the product of the
distinct primes dividing $n$ (the ``square-free kernel'' of $n$).
Variables in boldface type 
always indicate vector quantities with the components using the same
letter (e.g. $\zz = (z_1,z_2,\ldots)$).

\proclaim{Lemma 2.1} If $N>20$ and $x \ge 2N\log N$, there are
at least $N$ primes in the interval $(x,2x]$.  If $0<\del \le \frac12$,
$\frac{N}{\log N} \ge \frac{6}{\del}$, $x\ge e^{1.5+1.5/\del}$ and
$x\ge \frac{6}{\del}N\log N$, then there are at least $N$ primes in the
interval $(x,x+\del x]$.
\endproclaim

\demo{Proof} This comes directly from the following inequality due to
Rosser and Schoenfeld ([23], Theorems 1 and 2).
 Let $\pi(x)$ be the number of primes $\le x$.
Then for $x>67$ we have
$$
\frac{x}{\log x - 1/2} < \pi(x) < \frac{x}{\log x} \( 1 + \frac{3}{2\log x}\).
\tag{2.1}
$$
Thus for $x\ge 1200$, we have $\pi(2x)-\pi(x) \ge 0.735 \frac{x}{\log x}$.
Taking $x=2N\log N$ proves the first part of the
lemma for $N>130$.  For smaller $N$ we use a short computation.
For the second part, from (2.1) we obtain
$$
\pi(x+\del x)-\pi(x) \ge \frac{x(1+\del)}{\log x} \( 1 + \frac{1/2 - \log(1
+\del)}{\log x} \) - \frac{x}{\log x} - \frac{3x}{2\log^2 x}.
$$
Since $(1+\del)\log(1+\del) \le \del+\frac12 \del^2$, we have
$$
\split
\pi(x+\del x)-\pi(x) &\ge \frac{x}{\log x} \left[ \del - \frac{3/2-(1+\del)
(1/2-\log(1+\del))}{\log x} \right] \\
&\ge \frac{x}{\log x} \left[ \del - \frac{1+\del}{\log x} \right].
\endsplit
$$
Using the lower bounds for $x$ gives
$$
\pi(x+\del x)-\pi(x) \ge \frac{\del x}{3\log x} \ge \frac{2N\log N}
{\log N + \log(\tfrac{6}{\del}\log N)} \ge N. \quad\text{\qed}
$$
\enddemo

\proclaim{Lemma 2.2}  If $0\le \del \le \frac1{10}$, $u\ge 2-3\del$ and
$R \ge 6^{1/\del}$, then
$$
|\CC(R^u,R)| \ge \frac{\del^w}{(w+1)!} \frac{R^u}{\log R},
\quad w = \bigfl{\frac{u}{1-\del}}.
$$
\endproclaim

\demo{Proof}  Let $N_d(x,R) = | \{ n\in \CC(x,R): \Omega(n) \le d \}|$.
We show by induction on $d$ that
$$
N_d(R^u,R) \ge \frac{\del^{d-1}}{d!} \frac{R^u}{\log R} \qquad
(2-3\del \le u < d(1-\del), R\ge 6^{1/\del}). \tag{2.2}
$$
The proof uses another inequality due to Rosser and Schoenfeld ([23], Theorem
5), which states that for some constant $B$ and $x\ge 286$,
$$
\left| \sum_{p\le x} \frac{1}{p} - \log\log x - B \right| \le \frac{1}
{2\log^2 x}. \tag{2.3}
$$
In our applications, $x \ge 6^{1/(2\del)} \ge 6^5 > 286$.  First we 
establish (2.2) when $d=2$ and $d=3$.  Suppose $d=2$ and $2-3\del \le u
< 2-2\del$.  Then $N_2(R^u,R)$ is at least $\frac12$ of the number of
pairs of primes $(p_1,p_2)$ with $R^{u-1} < p_1 \le R$, $\sqrt{R} < p_2
\le R^u/p_1$.  Using $R\ge 6^{1/\del} \ge 6^{10}$, $R^u/p_1 \ge R^{0.7}$,
and (2.1), we have
$$
\split
N_2(R^u,R) &\ge \frac12 \sum_{R^{u-1}<p\le R} \( \pi\pfrac{R^u}{p} -
 \pi(\sqrt{R}) \) \\
&\ge \frac12  \sum_{R^{u-1}<p\le R} \frac{R^u/p}{\log R} \( 1 - 2R^{-0.2}
 \( 1 + \frac{3}{\log R} \) \) \\
&\ge \frac{0.46 R^u}{\log R} \sum_{R^{u-1}<p\le R} \frac{1}{p}.
\endsplit
$$
By (2.3), the last sum is
$$
\ge \log\pfrac{1}{u-1} - \frac{1}{2\log^2 R} \(1 + \frac{1}{(u-1)^2} \)
\ge \log\pfrac{1}{1-2\del}-\frac{\del^2}{2} \ge 2\del,
$$
and (2.2) follows when $d=2$.   Next, let $d=3$.  When $2-3\del \le 
u <2-2\del$, (2.2) follows from the $d=2$ case.
If $2-2\del \le u< 3-3\del$, define
$$
a_1 = \max\( \frac{u-1}{2}, \frac12 \), \qquad a_2 = \min\( 1, \frac{u-
1/2-\del}{2} \).
$$
Then
$$
N_3(R^u,R) \ge \frac16 \sum_{p_1,p_2 \in (R^{a_1},R^{a_2}]} \( \pi\pfrac
{R^u}{p_1p_2} - \pi(\sqrt{R}) \).
$$
For every $p_1,p_2$, 
$$
R \ge R^u/p_1p_2 \ge R^{u-2a_2} \ge R^{1/2+\del}.
$$
By (2.1),
$$ 
\split
\pi\pfrac{R^u}{p_1p_2} - \pi(\sqrt{R}) &\ge \frac{R^u/(p_1p_2)}{\log R} \(
1 - 2R^{-\del} \( 1 + \tfrac{3}{\log R} \) \) \\
&\ge \frac{0.61 R^u}{p_1p_2 \log R},
\endsplit
$$
whence
$$
N_3(R^u,R) \ge \frac{R^u}{10\log R} \biggl( \sum_{R^{a_1}<p\le R^{a_2}} 
\frac1{p} \biggr)^2.
$$
By (2.3),
$$
\sum_{R^{a_1} < p \le R^{a_2}} \frac{1}{p}
\ge \log\pfrac{a_2}{a_1} - \frac{1}{a_1^2 \log^2 R} \ge  \log\pfrac{a_2}{a_1}
 - 1.25 \del^2.
$$
We claim that $\log (a_2/a_1) \ge 1.5\del$, from which (2.2) follows in
the case $d=3$.  Let $I_1=[2-2\del,2)$, $I_2=[2,2.5+\del)$, 
$I_3=[2.5+\del,3-3\del)$.  Then
$$
\log\pfrac{a_2}{a_1} = \cases
\log(u-1/2-\del) \ge \log(1.5-3\del) \ge \log(1+2\del)\ge 1.5\del &
(u\in I_1) \\
\log\pfrac{u-1/2-\del}{u-1} \ge \log\pfrac{2}{1.5+\del} \ge \log(1.25)
 \ge 1.5\del & (u\in I_2) \\
\log\pfrac{2}{u-1} \ge\log\pfrac{2}{2-3\del} \ge 1.5\del & (u\in I_3).
\endcases
$$
Next, let $d\ge 3$ and suppose (2.2) holds.
When $2-3\del \le u < d(1-\del)$, (2.2) follows for all larger $d$ as well.
Suppose $d(1-\del) \le u
\le (d+1)(1-\del)$.  If $p\in (R^{1-\del},R]$, then 
$R^u/p \in (R^{2-3\del},R^{d(1-\del)}]$, and thus
$$
N_d(R^u/p,R) \ge \frac{\del^{d-1}}{d!} \frac{R^u/p}{\log R}.
$$
Summing over primes $p$, each number $pn$ with $n$ counted by
$N_d(R^u/p,R)$ is counted at most $d+1$ times.  Hence
$$
N_{d+1}(R^u,R) \ge \frac{1}{d+1} \sum_{R^{1-\del}<p\le R} N_d(R^u/p,R)
\ge \frac{\del^{d-1}}{(d+1)!} \frac{R^u}{\log R}  \sum_{R^{1-\del}<p\le R}
\frac{1}{p}.
$$
Again using (2.3), the last sum is
$$
\ge \log\pfrac{1}{1-\del}-\frac{1}{(1-\del)^2 \log^2 R} \ge
\del+\frac{\del^2}{2} - 0.4\del^2 > \del,
$$
and (2.2) follows with $d$ replaced by $d+1$.
\qed\enddemo

\proclaim{Lemma 2.3}
Suppose $R \ge (2u)^3 \ge 90000$.  Then $|\CC(R^u,R)| \le R^u (2/u)^u$.
\endproclaim

\demo{Proof}  Suppose $\frac23 \le \beta <1$ and put $P=R^u$.  Then
$$
\split
|\CPR| &\le P^\beta \sum_{n\in \CPR} n^{-\beta} \le P^\beta 
  \prod_{\sqrt{R} < p \le R}\( 1 + p^{-\beta} + p^{-2\beta} + \cdots \) \\
&\le P^\beta \exp \biggl\{ \sum_{\sqrt{R}<p\le R} \frac{1}{p^\beta} + 
  \frac{1}{p^\beta(p^\beta-1)} \biggr\} \\
&\le P^\beta \exp \biggl\{ R^{1-\beta}  \sum_{\sqrt{R}<p\le R} \frac{1}{p}
  + 1.03 \sum_{p> \sqrt{R}} p^{-4/3} \biggr\}.
\endsplit
$$
Since $\sqrt{R} \ge 300$, by (2.3)
$$
\sum_{\sqrt{R}<p\le R} \frac{1}{p} \le \log 2 + \frac{2.5}{\log^2 R} \le 0.713.
$$
Also,
$$
 \sum_{p> \sqrt{R}} p^{-4/3} \le \int_{\sqrt{R}-1}^\infty t^{-4/3}\, dt \le
0.45,
$$
so that
$$
|\CPR| \le P^\beta \exp \{ 0.713 R^{1-\beta}+0.47 \}.
$$
Take $\beta=1-\frac{\log(u/0.713)}{\log R} \ge \frac23$.  Then
$$
|\CPR| \le P \exp \{ -u\log(u/0.713)+u+0.47 \} = P \pfrac{0.713e}{u}^u
 e^{0.47}.
$$
Lastly, $u\ge 22$ and thus $(\frac{0.713e}{2})^u e^{0.47} < 1$.
\qed\enddemo

The next lemma is due to Wooley ([33]), and gives a bound for the number of
non-singular solutions of a system of congruences.  This greatly generalizes
a lemma due to Linnik [14].

\proclaim{Lemma 2.4}  Let $f_1, \ldots,f_d$ be polynomials in 
$\Bbb Z[x_1,\ldots,x_d]$ with respective degrees $k_1,\ldots,k_d$, and write
$$
J(\boldkey f;\xx)=\det \( \frac{\partial f_j(\xx)}{\partial x_i} \)_{1\le i,j
\le d}.
$$
Also, let $p$ be a prime number and $s$ be a natural number.  Then the number,
$N$, of solutions of the simultaneous congruences
$$
f_j(x_1,\ldots,x_d) \equiv 0 \pmod{p^s} \quad (1\le j\le d)
$$
with $1\le x_i\le p^s$ $(1\le i\le d)$ and $(J(\boldkey f;\xx),p)=1$,
satisfies $N \le k_1\cdots k_d$.
\endproclaim

Lastly, we present a general inequality on the number of solutions of
``symmetric'' systems of equations.

\def\zzero{\text{\bf 0}}
\proclaim{Proposition ZRD (Zero Representation Dominates)}  Suppose 
$f_1,\ldots,f_n$ are functions from $\Bbb Z^m$ to $\Bbb Z$ and
$\BB$ is a finite subset of $\Bbb Z^m$.
 Let $I(\ff; \ww;\BB)$ be the number of solutions of
the simultaneaous Diophantine equations
$$
f_j(\xx)-f_j(\yy) = w_j \qquad (1\le j\le n)
$$
with $\xx,\yy \in \BB$.  Then $I(\ff; \ww;\BB) \le I(\ff;\zzero;\BB)$, where
$\zzero=(0,0,\ldots,0)$.
\endproclaim

\demo{Proof}
For $\baa=(\a_1,\ldots,\a_n)$, let
$$
g(\baa) = \sum_{\xx \in \BB} e(\a_1 f_1(\xx) + \cdots + \a_n f_n(\xx)).
$$
Then
$$
I(\ff; \ww;\BB) = \int_{\Bbb U^n} |g(\baa)|^2 e(-\a_1 w_1 - \cdots - \a_k w_k)
\, d\baa \le  I(\ff;\zzero;\BB). 
$$
Alternatively, for $\vv=(v_1,\cdots,v_n)$, let $n(\vv)$ be the number of
solutions of $f_j(\xx) = v_j \; (1\le j\le n)$ with $\xx\in\BB$.  By
the Cauchy-Schwarz inequality,
$$
\split
I(\ff;\ww;\BB) &= \sum_{\Sb \vv, \vv' \\ v_j-v_j'=w_j \endSb} n(\vv) n(\vv')\\
&\le \biggl( \sum_{\Sb \vv, \vv' \\ v_j-v_j'=w_j \endSb} n(\vv)^2 \biggr)^{1/2}
 \biggl( \sum_{\Sb \vv, \vv' \\ v_j-v_j'=w_j \endSb} n(\vv')^2 \biggr)^{1/2}
=I(\ff;\zzero;\BB). \quad\text{\qed}
\endsplit
$$
\enddemo

\vfil\eject
%
%
%
%
%
%
\head 3. Vinogradov's Integral: Complete systems \endhead
%
%
%
%
%
%
In this section, we derive bounds for $J_{s,k}(P)$ using the
iterative methods of Wooley [30], modified using an idea of Arkhipov and
Karatsuba [1] (the introduction of the parameter $r$).
It should be noted that using the method of Tyrina [26] when
$\frac49 k^2 \le \Delta(k,s) \le \frac12 k^2$ gives slightly better
values for $\Delta(k,s)$, but only enough to improve the constant $B$
in Theorem 1 by $0.01$ or less.

The next definition is slightly different from that given in [30].

\proclaim{Definition}  Suppose $0\le d \le k-1$ and $T$ is a positive integer.
We say the $k$-tuple of poynomials $\bpsi=
(\Psi_1,\ldots,\Psi_k) \in \Bbb Z[x]^k$ 
is of type $(d,T)$ if $\Psi_j$ is identically zero for $j\le d$, and for
some integer $m\ge 0$, when
$j> d$, $\Psi_j$ has degree $j-d$ with leading coefficient
$\frac{j!}{(j-d)!}2^m T$.
\endproclaim

\proclaim{Lemma 3.1} Suppose $\bpsi$ is of type $(d,T)$, and 
$z_1,\ldots,z_{k-d}$ are integers.  Then
$$
\split
J_{k-d}(\zz;\bpsi) &:=  
\det \bigl( \Psi_j'(z_i) \bigr)_{1\le i\le k-d \atop d+1\le j\le k} \\
 &= (2^m T)^{k-d} \prod_{j=d+1}^k \frac{j!}{(j-d-1)!}
\prod_{1\le i<j \le k-d} (z_i-z_j),
\endsplit
$$
\endproclaim

\demo{Proof}  This follows by elementary row operations.
\qed\enddemo

The argument will begin with $\Psi_j(z)=z^j$ ($1\le j\le k$), which is of
type $(0,1)$.  At the $d$th iterative stage ($d\ge 0$), 
the system will be transformed
from one of type $(d,T)$ to one of type $(d+1,T')$ in two steps.
First, for some constant $c$ we will take
$$
\Phi_j(z)=\sum_{\ell=0}^j \binom{j}{\ell} \Psi_\ell(z) c^{j-\ell},
$$
which is also  a system of type $(d,T)$.  Then, for a constant $y$ we take
$$
\Upsilon_j(z) = \Phi_j(z+y)-\Phi_j(z) \qquad (1\le j\le k),
$$
which is of type $(d+1,yT)$.

Fix $k$ and suppose $1\le r\le k$.
If $\bpsi=(\Psi_1,\ldots\Psi_k)$ is a system of polynomials, let $K_{s}(P,Q;
\bpsi;q)$ be the number of solutions of the simultaneous equations
$$
\aligned
&\sum_{i=1}^k (\Psi_j(z_i)-\Psi_j(w_i)) + q^j 
\sum_{i=1}^s (x_i^j-y_i^j)=0
\qquad (1\le j\le k), \\
&1 \le z_i,w_i \le P; \quad 1 \le x_i,y_i \le Q.
\endaligned\tag{3.1}
$$
Here the inequalities on the variables $z_i,w_i,x_i,y_i$ hold for every $i$.
For prime $p$, let $L_{s}(P,Q;\bpsi;p,q,r)$
be the number of solutions of
$$
\aligned
&\sum_{i=1}^k
(\Psi_j(z_i)-\Psi_j(w_i)) + (pq)^j \sum_{i=1}^s (u_i^j-v_i^j)=0
\qquad (1\le j\le k), \\
&1 \le z_i,w_i \le  P; \;  z_i \equiv w_i \pmod{p^r}; \;
 1 \le u_i,v_i \le Q.
\endaligned\tag{3.2}
$$
Define the exponential sums
$$
\split
f(\baa)&=f(\baa;Q;q)=\sum_{x\le Q} e(\a_1 q x + \cdots + \a_k q^k x^k), \\
F(\baa)&=F(\baa;P;\bpsi) = \sum_{x\le P} e(\a_1 \Psi_1(x) + \cdots
+\a_k \Psi_k(x)).
\endsplit
$$
Then
$$
K_{s}(P,Q;\bpsi;q) = \iuk |F(\baa)^{2k} f(\baa)^{2s}|
\, d\baa.
$$
The next result relates $K_s$ and $L_s$, and
is a generalization of the ``fundamental lemma'' 
of Wooley ([30], Lemma 3.1).

%
%
%

\proclaim{Lemma 3.2}
Suppose $k$, $r$, $d$ and $s$ are integers with
$$
k\ge 4, \; 2 \le r \le k; \; 0\le d \le r-1; \; s\ge d+1.
$$
Let $M$, $P$ and $Q$ be real numbers with
$$
P^{\frac{1}{k+1}} \le M \le P^{\frac1{r}}; \quad
 32s^2 M < Q \le P;  \quad M\ge k.
$$
Suppose $q$ is a positive integer and
$\bpsi$ is a system of polynomials of type $(d,T)$ with $T\le P^d$.
Denote by $\curly P$ the set of the $k^3$ smallest primes $>M$,
and suppose $\curly P \subset (M,2M]$.
Then there is a
system of polynomials $\bphi$ of type $(d,T)$ and a prime $p \in 
\curly P$ such that
$$
K_{s}(P,Q;\bpsi;q) \le
4k^3 k! p^{2s+\frac12(r^2-r+d^2-d)} L_{s}(P,\tfrac{Q}{p};\bphi;p,q,r).
$$
\endproclaim

\demo{Proof}  Let $W$ be the set of systems of polynomials of type $(d,T)$
with $T\le P^d$.  Since $K_s(P,Q;\bpsi;q) \le P^{2k} Q^{2s}$ trivially,
there is a system $\bpsi_0 \in W$ so that
$$
K_s(P,Q;\bpsi_0;q) = \max_{\bpsi\in W} K_s(P,Q;\bpsi;q).
$$
We therefore assume without loss of generality that $\bpsi=\bpsi_0$.
For brevity, write $K$ for $K_{s}(P,Q;\bpsi;q)$.  We divide the solutions of
(3.1) into two classes: $S_2$ is the number of solutions with $z_i=z_j$
or $w_i=w_j$ for some $i\ne j$;
$S_1$ is the number of remaining solutions.
Clearly $K \le 2\max(S_1,S_2)$.  Suppose first that $S_2\ge S_1$.
By \Ho,
$$
\split
K &\le 2S_2 \le 4 \binom{k}{2} \iuk |F(\baa)^{2k-2} F(2\baa) f(\baa)^{2s}|
\, d\baa \\
&< 2k^2 \( \iuk |F(\baa)^{2k} f(\baa)^{2s}|\, d\baa \)^{1-\frac{1}{k}}
\( \iuk |f(\baa)|^{2s}\, d\baa \)^{\frac{1}{2k}} \( \iuk |F(2\baa)^{2k}
f(\baa)^{2s}|\, d\baa \)^{\frac{1}{2k}} \\
&=   2k^2  K^{1-1/k} (J_{s,k}(Q))^{1/2k} K_s(P,Q;2\bpsi;q) \\
&\le 2k^2  K^{1-1/2k}   (J_{s,k}(Q))^{1/2k}.
\endsplit
$$
Here $2\bpsi = (2\Psi_1(z), \ldots, 2\Psi_k(z))$ is also of type $(d,T)$,
which justifies the last inequality above.  This is the reason for the
introduction of the parameter $m$ in the definition of a system of
polynomials of type $(d,T)$.
Therefore $K\le (2k^2)^{2k} J_{s,k}(Q)$.  On the other hand, counting the
solutions of (3.1) with $z_i=w_i$ for each $i$ produces the lower
bound $K \ge (P-1)^k J_{s,k}(Q)$.  The hypothesis $\curly P \subset (M,2M]$
gives $M \ge k^3-1$ and so $P-1 \ge (k^3-1)^2-1 > 4k^4$.  We have a 
contradiction, therefore $K \le 2S_1$.  To bound $S_1$, we follow the
procedure from Wooley [30].
Consider a solution of (3.1) counted by $S_1$.  By Lemma 3.1,
for some integer $m\ge 0$ we have
$$
\split
J_{k-d}(\zz;\bpsi) J_{k-d}(\ww;\bpsi) &= (2^m T)^{2k-2d} \!
\prod_{j=d+1}^k \! \( \! \tfrac{j!}{(j-d-1)!}
\)^2\!\! \prod_{1\le i<j \le k-d} \! (z_i-z_j)(w_i-w_j) \\
&\ne 0.
\endsplit
$$
By hypothesis, if $p\in \curly P$ then $p>M\ge k$.  Also,
$$
\left| T  \prod_{1\le i<j \le k-d} (z_i-z_j)(w_i-w_j) \right| <
P^{d+(k-d)(k-d-1)} \le P^{k^2-k} < \prod_{p\in \curly P} p.
$$
Thus, for each solution counted by $S_1$,
there is some $p\in \curly P$ which does not divide $J_{k-d}(\zz;\bpsi)
J_{k-d}(\ww;\bpsi)$.
Hence
$$
K \le 2 k^3 \max_{p\in \curly P} S_3(p), \tag{3.3}
$$
where $S_3(p)$ is the number of solutions of (3.1) with $(p,J_{k-d}(\zz;\bpsi)
J_{k-d}(\ww;\bpsi))=1$.
With $p$ fixed, let
$$
\split
g(\baa;b) &= \sum_{x\le Q \atop x\equiv b\pmod{p}} e(\a_1 q x + \cdots +
\a_k q^k x^k), \\
\FT(\baa) &= \sum_{\Sb z_1,\ldots,z_k \\ (J_{k-d}(\zz;\bpsi),p)=1\endSb}
e\( \sum_{j=1}^k \a_j(\Psi_j(z_1)+\cdots +\Psi_j(z_k)) \).
\endsplit
$$
Since $\bpsi$ is of type $(d,T)$, for any solution of (3.1) we have
$$
\sum_{i=1}^s (x_i^j-y_i^j)=0 \qquad (1\le j\le d).
$$
Let $\BB_s(\ww)$ denote the set of solutions (with $0\le c_i \le p-1$ for each
$i$) of the system of congruences
$$
\sum_{i=1}^s c_i^j \equiv w_j \pmod{p} \qquad (1\le j\le d).
$$
Consequently,
$$
S_3(p) \le \iuk |\FT(\baa)|^2 
\sum_{\ww \atop 1\le w_j\le p} |U(\baa;\ww)|^2\, d\baa,
$$
where
$$
U(\baa;\ww) = \sum_{\cc \in \BB_s(\ww)} g(\baa;c_1) \cdots g(\baa;c_s).
$$
By first fixing $c_{d+1},\ldots,c_s$, we have $|\BB_s(\ww)| \le p^{s-d}
 \max_{\vv}|\BB_d(\vv)|$.
  Suppose $\cc$ and $\cc'$ are two solutions counted in
$\BB_d(\vv)$.  Let $q(t) = (t-c_1)\cdots (t-c_d)$.  By Newton's formulas
connecting the sums of the powers of the roots of a polynomial with its 
coefficients, $q(t) \equiv (t-c_1') \cdots (t-c_d') \pmod{p}$.  
Thus, $\cc'$ is
a permutation of $\cc$, whence $|\BB_d(\vv)| \le d!$ and 
$$
|\BB_s(\ww)| \le d! p^{s-d}.
$$
By \CS, followed by an application of the arithmetic mean-geometric mean
inequality, we have
$$
\split
|U(\baa;\ww)|^2 &\le |\BB_s(\ww)| \sum_{\cc \in \BB_s(\ww)} 
|g(\baa;c_1) \cdots g(\baa;c_s)|^2 \\
&\le \frac{d!}{s} p^{s-d}  \sum_{\cc \in \BB_s(\ww)} \sum_{i=1}^s
|g(\baa;c_i)|^{2s}.
\endsplit
$$
We then have
$$
\aligned
S_3(p) &\le d! p^{s-d} \sum_{\cc} \max_{1\le i\le s} \iuk |\FT(\baa)|^2
|g(\baa;c_i)|^{2s}\, d\baa \\
&\le d! p^{2s-d} \max_{0\le c\le p-1} S_4(c,p),
\endaligned\tag{3.4}
$$
where
$$
S_4(c,p) = \iuk |\FT(\baa)^2 g(\baa;c)^{2s}|\, d\baa
$$
is the number of solutions of
$$
\aligned
&\sum_{i=1}^k
(\Psi_j(z_i)-\Psi_j(w_i)) + q^j \sum_{i=1}^s ((pu_i-c)^j-(pv_i-c)^j)
=0 \quad (1\le j\le k), \\
&1\le z_i,w_i\le P; \quad (p,J_{k-d}(\zz;\bpsi)J_{k-d}(\ww;\bpsi))=1;
 \quad 1\le u_i,v_i \le (Q+c)/p.
\endaligned\tag{3.5}
$$
Let $S_5(c,p)$ denote the number of solutions of (3.5) with
$u_i > Q/p$ or $v_i > Q/p$ for some $i$, and let $S_6(c,p)$ denote
the number remaining solutions.  Suppose first that $S_5(c,p) \ge S_6(c,p)$.
By \Ho,
$$
\split
S_4(c,p) &\le 2S_5(c,p) \le 4s \iuk |\FT(\baa)^2 g(\baa;c)^{2s-1}|\, d\baa \\
&\le 4s \( \iuk  |\FT(\baa)^2 g(\baa;c)^{2s}|\, d\baa \)^{1-\frac{1}{2s}} \(
\iuk |\FT(\baa)|^2\,d\baa \)^{\frac{1}{2s}} \\
&= 4s \( S_4(c,p) \)^{1-\frac{1}{2s}} \(
\iuk |\FT(\baa)|^2\,d\baa \)^{\frac{1}{2s}}.
\endsplit
$$
Therefore,
$$
S_4(c,p) \le (4s)^{2s} \iuk |\FT(\baa)|^2\,d\baa.
$$
Note that $\fl{(Q+c)/p} > Q/p$ in this case.
Thus, counting only the solutions of (3.5) with $u_i=v_i$ for every $i$
gives
$$
S_4(c,p) \ge (Q/p)^{s} \iuk |\FT(\baa)|^2\,d\baa.
$$
By our assumed lower bound on $Q$, this is impossible.  Therefore,
$S_4(c,p) \le 2S_6(c,p)$.  By the binomial theorem,
$$
(py)^j = \sum_{\ell=0}^j \binom{j}{\ell} (py-c)^\ell c^{j-\ell}.
$$
Thus, $S_6(c,p)$ is the number of solutions of
$$
\aligned
&\sum_{i=1}^k (\Phi_j(z_i)-\Phi_j(w_i)) + (pq)^j \sum_{i=1}^s (u_i^j-v_i^j)
=0 \qquad (1\le j\le k), \\
&1\le z_i,w_i\le P; \quad (p,J_{k-d}(\zz;\bpsi)J_{k-d}(\ww;\bpsi))=1;
 \quad 1\le u_i,v_i \le Q/p,
\endaligned\tag{3.6}
$$
where, for $1\le j\le k$,
$$
\Phi_j(z) = \sum_{\ell=0}^j \binom{j}{\ell} \Psi_\ell(z) c^{j-\ell}.
$$
The leading coefficients of $\Phi_j$ and $\Psi_j$ are equal, hence
$\bphi$ is also of type $(d,T)$ (with the same value of $m$).
By Lemma 3.1, $J_{k-d}(\zz;\bpsi)=J_{k-d}(\zz;\bphi)$,
so $(p,J_{k-d}(\zz;\bphi)J_{k-d}(\ww;\bphi))=1$ in (3.6).

Lastly, we introduce the congruence condition on $z_i,w_i$.
By (3.6),
$$
\sum_{i=1}^k (\Phi_j(z_i)-\Phi_j(w_i)) \equiv 0 \pmod{p^j} \qquad
(1\le j\le k).
$$
We shall only work with the congruences corresponding to $d+1 \le j\le k$,
since the left side of the above congruence is identically zero when
$j\le d$.
Let $\BB^*(\mm)$ be the set of $\zz$ with $1\le z_i\le p^r$ for each $i$,
$(J_{k-d}(\zz;\bphi),p)=1$ and
$$
\sum_{i=1}^k \Phi_j(z_i) \equiv m_j \pmod{p^{\min(j,r)}} \qquad 
(d+1\le j\le k).
$$
By hypothesis, $d+1 \le r$.  To bound $|\BB^*(\mm)|$,
first fix $z_{k-d+1},\ldots,z_k$ (there are $p^{rd}$ such choices). 
For each $j$, there are $p^{\max(0,r-j)}$ possibilities for  $m_j$ modulo
$p^r$, and with the $m_j$ fixed modulo $p^r$,
 Lemma 2.4 implies that there are at most
$(k-d)!$ solutions $z_1,\ldots,z_{k-d}$ modulo $p^r$.  Therefore,
$$
|\BB^*(\mm)| \le (k-d)! p^{\frac12(r-d-1)(r-d)+rd}.
$$
Define
$$
H(\baa;\zz) = \sum_{\Sb \ww \\ 1\le w_i \le P \\ w_i\equiv z_i\pmod{p^r} 
\endSb}
e \( \sum_{j=1}^k \a_j (\Phi_j(w_1) + \cdots + \Phi_j(w_k))\).
$$
Then, by \CS,
$$
\split
S_6(c,p) &\le  \iuk \sum_{\mm} \biggl| \sum_{\zz \in \BB^*(\mm)}
H(\baa;\zz) \biggr|^2 |f(\baa;Q/p;pq)|^{2s}\, d\baa \\
&\le  \sum_{\mm} |\BB^*(\mm)| \iuk  \sum_{\zz \in \BB^*(\mm)}
|H(\baa;\zz) |^2 |f(\baa;Q/p;pq)|^{2s}\, d\baa \\
&\le (k-d)! p^{\frac12(r-d-1)(r-d)+rd}
 L_{s}(P,Q/p;\bphi;p,q,r).
\endsplit
$$
By (3.4) and the inequality $d!(k-d)! \le k!$,
$$
S_3(p) \le 2k! p^{2s-d+\frac12(r-d-1)(r-d)+rd} L_{s}(P,Q/p;\bphi;p,q,r).
\tag{3.7}
$$
The lemma now follows from (3.3).
\qed\enddemo

%
%
%

\proclaim{Lemma 3.3}  Suppose that $s\ge d$, $k\ge r\ge 2$, $d\le k-2$,
$q\ge 1$, $p$ is a prime and $\bphi$ is a system of polynomials of type
$(d,T)$.  Then there is a system of polynomials $\bups$ of type
$(d+1,T')$ with $T\le T' \le PT$  such that
$$
L_{s}(P;Q;\bphi;p,q,r) \le (2P)^k\max\bigl[k^k J_{s,k}(Q), 2p^{-rk}  \left\{
J_{s,k}(Q) K_{s}(P,Q;\bups;pq) \right\}^{1/2}\bigr].
$$
\endproclaim

\demo{Proof}  For short, write $L$ for $L_{s}(P;Q;\bphi;p,q,r)$.  Then
$L\le 2\max(U_0,U_1)$, where $U_0$ is the number of solutions of (3.2)
with $w_i=z_i$ for some $i$, and $U_1$ is the number of solutions of (3.2)
with $w_i\ne z_i$ for every $i$.  First write $f(\baa)$ for $f(\baa;Q;pq)$ and
$$
I(\baa) = \sum_{1\le c\le P} \biggl| \sum_{\Sb 1\le w\le P 
\\ w\equiv c\pmod{p^r}
\endSb} e(\a_1 \Phi_1(w) + \cdots + \a_k \Phi_k(w)) \biggr|^2,
$$
so that
$$
L = \iuk I(\baa)^k |f(\baa)|^{2s}\, d\baa.
$$
Suppose first that $U_0 \ge U_1$.  By \Ho,
$$
\split
L \le 2U_0 &\le 2kP \iuk I(\baa)^{k-1} |f(\baa)|^{2s}\, d\baa \\
&\le 2kP\( \iuk |I(\baa)^k f(\baa)^{2s}|\, d\baa \)^{1-1/k}
\( \iuk  |f(\baa)|^{2s}\, d\baa \)^{1/k} \\
&= 2kP L^{1-1/k} J_{s,k}(Q)^{1/k},
\endsplit
$$
and the lemma follows in this case.  If $U_1\ge U_0$, for each $i$ we may
write $w_i = z_i + h_i p^r$, where $1 \le |h_i| \le P/p^r$.
We may assume that $P/p^r\ge 1$, else $U_1=0$.  Let
$$
g(\baa;h) = \sum_{1\le z\le P} e \( \sum_{j=1}^k \a_j(\Phi_j(z+hp^r) -
\Phi_j(z)) \).
$$
There are $2^k$
choices for the signs of $w_i-z_i$ ($1\le i\le k$), so
$$
L \le 2 \sum_{\Sb \eta_1,\ldots,\eta_k \\ \eta_i\in \{-1,+1\} \endSb}
\iuk \sum_{\Sb \hh \\ 1\le h_i \le P/p^r \endSb} g(\eta_1\baa;h_1)
\cdots g(\eta_k \baa;h_k) | f(\baa)|^{2s}\, d\baa.
$$
Since $|g(\baa;h)| = |g(-\baa;h)|$,
$$
\sum_{\Sb \hh \\ 1\le h_i \le P/p^r \endSb} |g(\eta_1 \baa;h_1)
\cdots g(\eta_k \baa;h_k)| \le (P/p^r)^k \max_{1\le h\le P/p^r} |g(\baa;h)|^k.
$$
Then, by \CS,
$$
\split
L &\le 2^{k+1}  (P/p^r)^k \max_{1\le h\le P/p^r} \iuk |g(\baa;h)|^k |f(\baa)|
^{2s}\, d\baa \\
&\le 2^{k+1}  (P/p^r)^k \max_{1\le h\le P/p^r} \( \iuk |g(\baa;h)|^{2k}
|f(\baa)|^{2s}\, d\baa \)^{1/2} \( \iuk |f(\baa)|^{2s}\, d\baa 
\)^{1/2} \\
&= 2^{k+1} (P/p^r)^k  \max_{1\le h\le P/p^r} \bigl( K_s(P,Q;\bups;pq) J_{s,k}
(Q) \bigr)^{1/2},
\endsplit
$$
where $\Upsilon_j(z) = \Phi_j(z+hp^r)-\Phi_j(z)$  for $j\ge d+2$ and
$\Upsilon_j(z) \equiv 0$ for $j\le d+1$.
For some integer $m\ge 0$ and
 $j\ge d+2$, $\Upsilon_j$ has degree $j-d-1$ and leading coefficient
$\frac{j!}{(j-d-1)!} hp^r2^m T$,
thus the system $\bups$ is of type $(d+1,Thp^r)$.
\qed
\enddemo

%
%
%
%
%

Next, we iterate Lemmas 3.2 and 3.3 to produce a bound for $J_{s+k,k}(P)$
in terms of the bounds for $J_{s,k}(Q)$.

\proclaim{Lemma 3.4}  Suppose $k\ge 26$, $4\le r\le k$, $k\le s\le k^3$ and
$$
J_{s,k}(Q) \le C Q^{2s-\frac12k(k+1)+\Delta} \qquad (Q\ge 1).
$$
Let $j$ be an integer satisfying
$$
2 \le j \le \tfrac{9r}{10}, \quad (j-1)(j-2) \le 2\Delta-(k-r)(k-r+1).
\tag{3.8}
$$
Define
$$
\phi_j = \frac1{r}, \quad \phi_{J} = \frac1{2r} + 
\frac{k^2+k+r^2-r+J^2-J-2\Delta}{4kr}\, \phi_{J+1} \quad (1\le J\le j-1),
$$
and suppose $r$ and $j$ are chosen so that $\phi_i \ge \frac{1}{k+1}$ for
every $i$.  Suppose 
$$
\frac{1}{3\log k} \le \om \le \frac12, \quad \eta=1+\om, \quad
V=\max\( e^{1.5+1.5/\om}, \frac{18}{\om}k^3\log k\).
$$
If $P \ge V^{k+1}$, then
$$
J_{s+k,k}(P) \le k^{3k} \eta^{4s+k^2} C P^{2(s+k)-\frac12k(k+1)+\Delta'},
$$
where $\Delta' = \Delta(1-\phi_1) - k + \frac{\phi_1}{2}(k^2+k+r^2-r)$.
\endproclaim

\demo{Proof} Let $Q_0=P$ and for $1\le i\le j$ define
$$
M_i = P^{\phi_i}, \qquad Q_i = P^{1-(\phi_1 +\cdots+\phi_i)}.
$$
Let $\curly P_i$ be the set of $k^3$ smallest primes $>M_i$.
By hypothesis, $M_i \ge V$, and by the definition of $\eta$ and $V$,
Lemma 2.1 implies that
$\curly P_i \subset (M_i,\eta M_i]$.
By (3.8), $\phi_i \le \frac{1}{r}$ for each $i$, and for $i\le j-1$
$$
\aligned
Q_i \ge Q_{j-1} &\ge P^{1-(j-1)/r} \ge P^{1/10+1/r} \\
&\ge V^{k/10} P^{\phi_{i+1}} > k^8 P^{\phi_{i+1}}
> 32 s^2 P^{\phi_{i+1}}.
\endaligned\tag{3.9}
$$
Let $\lambda = 2s-\frac12k(k+1)+\Delta$.
We shall show by induction on $J$ that for every system $\bphi$ of
type $(J,T)$ with $1\le T\le P^J$,
every prime $p\in \curly P_{J+1}$ and every positive integer $q$,
$$
L_s(P,Q_{J+1};\bphi;p,q,r) \le E_J C P^k Q_{J+1}^\lambda, \tag{3.10}
$$
where
$$
E_{j-1}=1, \quad E_{J-1} = k^k \eta^{s+\frac14(k^2-k+J^2-J)} E_J^{1/2}
\;\; (1\le J\le j-1).
$$
First, when $J=j-1$, we have $p^r > M_{j-1}^r \ge P$, 
so that in (3.2), $w_i=z_i$ for every $i$.  This gives
$$
L_s(P,Q_{j};\bphi;p,q,r) \le P^k J_{s,k}(Q_j), 
$$
which gives (3.10) for $J=j-1$.  Now suppose $1\le J \le j-1$
and (3.10) holds.  Let $\bpsi$ be a system of polynomials of type
$(J,T)$ with $1\le T\le P^J$, and let $q'$ be any positive integer.
By (3.9), (3.10) and the fact that $L_s(P,Q;\bphi;p,q,r)$ is a
non-decreasing function of $Q$, we find from Lemma 3.2 that
$$
K_s(P,Q_J; \bpsi;q') \le 4k^3 k! 
(\eta M_{J+1})^{2s+\frac12 (r^2-r+J^2-J)} E_J CP^k  Q_{J+1}^{\lambda}.
$$
By Lemma 3.3, for every system of polynomials $\bphi$ of type
$(J-1,T)$ with $1\le T\le P^{J-1}$,
prime $p\in  \curly P_J$ and integer $q$, there is a system $\bpsi$ of
polynomials of type $(J,T')$ with $T' \le P^J$ such that
$$
\split
L_s&(P,Q_J; \bphi;p,q,r) \le (2P)^k \max \bigl[k^k CQ_J^\lambda,
2 P^{-kr\phi_J}
  \( CQ_{J}^\lambda K_s(P,Q_J; \bpsi;pq) \)^{\frac12}\bigr] \\
&\le CQ_J^\lambda (2P)^k \max \left[ k^k, 4(k^3 k!)^{\frac12} 
  E_J^{\frac12} P^{\frac{k}{2}-kr\phi_J} 
  M_{J+1}^{-\frac{\lambda}2} (\eta M_{J+1})^{s+\frac14(r^2-r+J^2-J)} \right].
\endsplit
$$
By the definition of $\phi_i$,
$$
\frac{k}{2} - kr\phi_J  + \frac12 \( \frac{k(k+1)}2 - \Delta + \frac12(r^2-r
+J^2-J) \) \phi_{J+1} = 0,
$$
i.e.,
$$
P^{k/2-kr\phi_J} M_{J+1}^{s-\lam/2+\frac14(r^2-r+J^2-J)} = 1.
$$
Since $r\le k$ and $4(k^3 k!)^{1/2} \le 2^{-k}k^k$ for $k\ge 8$,
this implies
$$
L_s(P,Q_J; \bphi;p,q,r) \le CQ_J^{\lam} (kP)^k \max \( 2^k, 
E_J^{1/2} \eta^{s+\frac14(k^2-k+J^2-J)} \).
$$
Next, $E_J \ge 1$ and
$$
\eta^{s+\frac14(k^2-k)} \ge \( \( 1 + \frac{1}{3\log k} \)^{\frac{k+3}{4}} \)^k
 \ge 2^{k} \quad (k\ge 26).
$$
Therefore, by the definition of $E_{J-1}$,
$$
L_s(P,Q_J; \bphi;p,q,r) \le CE_{J-1} P^k Q_{J}^\lambda,
$$
i.e., (3.10) follows with $J$ replaced by $J-1$.  Finally, taking (3.10)
with $J=0$ and applying Lemma 3.2 with $\Psi_j(x)=x^j$ for each $j$ gives
$$
\split
K_s(P,P;\bpsi;1) &\le 4k^3 k! (\eta M_1)^
{2s+\frac12(r^2-r)} E_0 C P^k Q_1^\lambda\\
&\le CP^{\lambda+k} 4k^3 k! \eta^{2s+\frac12(k^2-k)}
 E_0 M_1^{\frac12(k^2+k +r^2-r)-\Delta}.
\endsplit
$$
From the definition of $E_J$, we have
$$
\split
E_0 &= \prod_{J=1}^{j-1} \pfrac{E_{J-1}}{\sqrt{E_J}}^{2^{1-J}} 
E_{j-1}^{2^{1-j}} \\
&\le\prod_{J=1}^\infty  \(k^k \eta^{s+\frac14(
k^2-k+J^2-J)} \)^{2^{1-J}} \\
&= k^{2k} \eta^{2s+\frac12 k^2 - \frac12 k+2}.
\endsplit
$$
Lastly, $4k^3k! \le k^k$ for $k\ge 11$.  Therefore
$$
J_{s+k,k}(P) = K_s(P,P;\bpsi;1) \le k^{3k}\eta^{4s+k^2} 
C P^{2(s+k)-\frac12 k(k+1) +\Delta'}. \quad \text{\qed}
$$
\enddemo

For a given $k,r,\Delta$, we let $\d_0(k,r,\Delta)$ be
the value of $\Delta'$ coming from Lemma 3.4, where we take $j$ maximal
satisfying (3.8).  The optimal value of $r$ is about $\sqrt{k^2+k-2\Delta}$,
but leads to very messy analysis.  Making the choice
$r \approx k(1-\Delta/k^2)$ simplifies matters and ultimately
increases the value of $B$ in Theorem 1 by only about $0.0074$.

%
%
%

\proclaim{Lemma 3.5} Let $k\ge 26$ and let $\om$, $\eta$ and $V$ be 
as in Lemma 3.4.
Let $\Delta_1=\frac12 k^2(1-1/k)$ and 
for $n\ge 1$, let $r_n$ be an integer in $[4,k]$ satisfying
$$
\phi^*(k,r_n,\Delta_n) :=
\frac{2k}{2r_n k+2\Delta_n - (k-r_n)(k-r_n+1)} \ge \frac{1}{k+1},
\tag{3.11}
$$
then set $\Delta_{n+1} = \del_0(k,r_n,\Delta_n)$.  If $n\le k^2$, then
$$
J_{nk,k}(P) \le C_n P^{2nk-\frac12 k(k+1)+\Delta_n} \qquad (P\ge 1),
$$
where $C_1=k!$ and for $n\ge 2$
$$
C_n = C_{n-1} \max \left[ k^{3k} \eta^{4k(n-1)+k^2}, 
V^{(k+1)(\Delta_{n-1}-\Delta_n)} \right].
$$
\endproclaim

\demo{Proof} Defining $\phi_i$ as in
Lemma 3.4, we must ensure that $\phi_i \ge \frac{1}{k+1}$
for each $i$.
To this end, let $r=r_n$, $\Delta=\Delta_n$,
 $\phi^*=\phi^*(k,r,\Delta)$ and
$y=2\Delta-(k-r)(k-r+1)$.  For $i\ge 1$ let
$\theta_i = \phi_i-\phi^*$.  By (3.8),
$y-(j-1)(j-2) \ge 0$, so $\theta_j = 1/r - \phi^* \ge 0$.  Also,
$$
\theta_J = \frac{\theta_{J+1}}{4kr} (2rk+J^2-J-y) + \frac{J^2-J}{4kr}\phi^*
\quad(1\le J\le j-1).
$$
Since $2\Delta \le k^2-k$, $0 \le 2rk+J^2-J-y \le 2rk$.  It 
follows that for $J\le j-1$,
$$
0 \le \theta_J \le \frac{\theta_{J+1}}{2} + \frac{J^2-J}{4kr}\phi^*.
\tag{3.12}
$$
Thus, (3.11) and (3.12) imply that 
$\phi_i \ge \phi^* \ge \frac{1}{k+1}$ for every $i$.
We now proceed by induction, noting that the lemma holds with $n=1$
by the inequality $J_{k,k}(P) \le k! P^k$.
Assume now that $m\ge 2$ and the lemma holds for $n\le m-1$.  By Lemma 3.4,
$$
J_{mk,k}(P) \le C_{m-1} k^{3k} \eta^{4k(m-1)+k^2}
P^{2mk-\frac12 k(k+1) + \Delta_{m}} \qquad (P\ge V^{k+1}).
$$
For $P<V^{k+1}$, we have trivially
$$
\split
J_{mk,k}(P) &\le P^{2k} J_{(m-1)k,k}(P) \le
C_{m-1} P^{2mk-\frac12 k(k+1)+\Delta_{m-1}} \\
&\le C_{m-1} V^{(k+1)(\Delta_{m-1}-\Delta_m)} P^{2mk-\frac12 k(k+1)+
\Delta_{m}}.
\endsplit
$$
This completes the proof.
\qed\enddemo

For a particular choice of $r_1, r_2, \cdots$, 
the next lemma gives clean upper bounds on $\Delta_n$
and $C_n$ for large $k$. 

%
%
%
%
%

\proclaim{Lemma 3.6}  Suppose that $k\ge 1000$.  For
$$
2k \le n \le \frac{k}{2} \( \frac12 +\log\pfrac{3k}{8} \) + 1,
$$
we have
$$
J_{nk,k}(P) \le C_n P^{2nk-\kk+\Delta_n} \qquad (P\ge 1),
$$
where
$$
\split
\Delta_n &\le \frac{3}{8} k^2 e^{1/2 - 2n/k+1.69/k}, \\
C_n &\le k^{2.055k^3-5.91k^2+3nk} 1.06^{nk^2+2k(n^2-n)-9.7278k^3}.
\endsplit
$$
\endproclaim

\def\del{\delta}
\demo{Proof}  We shall take $r_n = \fl{k-\Delta_n/k+1}$ in Lemma 3.5.
For each $n$ write $\del_n=\Delta_n/k^2$.  Fix $n\ge 2$
and write $\del = \del_{n-1}$, $\del'=\del_n$, $\Delta=\Delta_{n-1}$,
$\Delta'=\Delta_n$, $r=r_{n-1}$.
If $\Delta_{n-1} \le k$, the upper bound for $\Delta_n$ in the
lemma follows from the upper bound on $n$, so from
now on assume that
$$
\Delta_{n-1} > k. \tag{3.13}
$$
We first show that
$$
\del' \le \del \( 1 - \frac{2-\del}{2-\del^2}\(  
\frac{2}{k} - \frac{32}{21k^2} - \frac{16}{7\delta k^3} \) \).
\tag{3.14}
$$
Let
$$
y = 2\Delta-(k-r)(k-r+1), \quad \phi^* = \phi^*(k,r,\Delta) = \frac{2k}{2rk+y}.
$$
By the definition of $r_n$,
$$
k\del(2k-k\del-1) \le y \le k\del(2k-k\del+1).
$$
Hence
$$
\phi^* \ge \frac{2k}{2k(k-k\del+1)+2\del k^2 - k\del(k\del-1)}
= \frac{2}{(2-\del^2)k+2+\del} \ge \frac{1}{k+1},
$$
so (3.11) holds.  Iterating (3.12) gives
$$
\theta_1 \le 2^{1-j}\theta_j + \sum_{h=1}^{j-1} 2^{1-h} (h^2-h) \frac{\phi^*}
{4kr}
\le 2^{1-j}\theta_j + \frac{2\phi^*}{kr}
\le \frac{2^{1-j}}{r} +  \frac{2\phi^*}{kr}.
$$
Next, (3.13) implies
$y\ge 2k-2$.  Since
$\sqrt{2k-2} \le k/3$, we always have $j\ge \sqrt{2k-2}$ (since $j$ is
maximal satisfying (3.8)) and so
for $k\ge 1000$
$$
\frac{2^{1-j}}{r} \le \frac{2^{1-\sqrt{2k-2}}}{r} \le \frac{0.071}
{k^4 r}.
$$
Also, $\del \le \frac12 (1-1/k)$ implies
$$
\phi^* \le \frac{2}{(2-\del^2)k-\del} \le \frac{8}{7k+1/k} < \frac{8}{7k} -
\frac{0.16}{k^3},
$$
and thus
$$
\theta_1 \le \frac{0.071}{k^4 r} + \frac{16}{7k^2 r} - \frac{0.32}{k^4 r}
\le \frac{16}{7k^2 r}.
$$
Since $\Delta \ge k$,
$$
k^2+k-2\Delta = (k-\Delta/k)^2+k-(\Delta/k)^2 \le (k-\del k)(k-\del k+1).
$$
Therefore, from $k-\del k \le r \le k-\del k+1$ 
and the upper bound on $\theta_1$,
$$
\aligned
\Delta' &= \Delta - k + \frac{\phi^*+\theta_1}{2} ( 2kr-y) \\
&\le \Delta - k + \frac{\phi^*}{2} (2kr-y) + \frac{8}{7k^2} \( r-1 + 
\frac{k^2+k-2\Delta}{r} \) \\
&\le \Delta-2k+4k^2 \, \frac{r}{2rk+y}+ \frac{16(1-\del)}{7k}.
\endaligned\tag{3.16}
$$
Next we establish
$$
\frac{1-\del}{2k} \le \frac{r}{2rk+y} \le 
\frac{1-\del}{(2-\del^2)k} + \frac{\del}{(2-\del^2)^2k^2}.
\tag{3.17}
$$
As a function of the real variable $r$, $\frac{r}{2rk+y}$ has positive
second derivative
and a minimum at $r=r_0:=\sqrt{k^2+k-2\Delta}$.  
Therefore, on the interval $[k-k\del,k-k\del+1]$,
the maximum occurs at one of the endpoints.  When $\del > 1/\sqrt{k}$,
$r_0\le k-k\del$, so the minimum occurs at $r=k-k\del$.  When
$\frac1{k} \le \del \le 1/\sqrt{k}$, $k-k\del \le r_0 \le k-k\del+1$,
so the minimum occurs at $r=r_0$.
At $r=k-k\del+1$,
$$
\frac{r}{2rk+y} = \frac{1-\del+\frac{1}{k}}{(2-\del^2)k+2+\del}
=\frac{1-\del}{(2-\del^2)k} \(
1 + \frac{\del}{k(1-\del)(2-\del^2)+2-\del-\del^2} \),
$$
so (3.17) holds for this $r$.  When $r=k-k\del$,
$$
\frac{r}{2rk+y} = \frac{1-\del}{(2-\del^2)k-\del} = 
\frac{1-\del}{(2-\del^2)k} \(1 + \frac{\del}{(2-\del^2)k-\del} \).
$$
Since $(2-\del^2)k-\del > (2-\del^2)k-\del k(2-\del^2)=(2-\del^2)(1-\del)k$,
(3.17) holds for this $r$ as well.  Lastly, when $\frac1{k} \le \del\le
 1/\sqrt{k}$ and $r=r_0$,
$$
\frac{r}{2rk+y} = \frac{1}{4k+1-2\sqrt{k^2+k-2\Delta}}.
$$
Also,
$$
\split
\bigl(k+1/2-(\del+\del^2)k\bigr)^2 &= k^2+k+\tfrac14 - k(2k+1)(\del+\del^2) +
  k^2(\del+\del^2)^2 \\
&\le  k^2+k+\tfrac14  - 2k^2(\del+\del^2) + k^2(\del+\del^2)^2 \\
&= k^2+k-2\del k^2 + \tfrac14 - k^2(\del^2-2\del^3-\del^4) \\
&<  k^2+k-2\del k^2.
\endsplit
$$
Therefore,
$$
\frac{r}{2rk+y} \ge \frac{1}{4k+1-2(k+1/2-(\del+\del^2)k)} =
\frac{1}{2k(1+\del+\del^2)} > \frac{1-\del}{2k}.
$$
This proves (3.17).

By (3.16) and (3.17), plus
the inequality $\frac{(1-\del)(2-\del^2)}{2-\del} \le 1$, we have
$$
\split
\del' &\le \del - \frac{2}{k} + \frac{4-4\del}{(2-\del^2)k} + \frac{4\del}
{(2-\del^2)^2 k^2} + \frac{16(1-\del)}{7k^3} \\
&= \del \( 1 - \frac{4-2\del}{(2-\del^2)k} + \frac{4}{(2-\del^2)^2k^2} \)
+\frac{16(1-\del)}{7k^3} \\
&\le \del\(  1 - \frac{2-\del}{2-\del^2}\(\frac{2}{k} - \frac{32}{21k^2}\)
\) + \frac{16}{7k^3}\, \frac{2-\del}{2-\del^2} \\
&= \del\(  1 - \frac{2-\del}{2-\del^2}\(\frac{2}{k} - \frac{32}{21k^2}
- \frac{16}{7 \del k^3} \)\).
\endsplit
$$
This concludes the proof of (3.14).  
We now use (3.14) to bound $\Delta_n$ and $C_n$.
Let $\beta=\frac2{k} - \frac{32}{21k^2}$, $c=\frac{16}{7k^3}$ and
$\beta'=\beta-c/\del$.
The differential equation analogous to (3.14) is approximately
$dy/dx=-\beta y \frac{2-y}{2-y^2}$, which has the implicit solution
$y+\log y+\log(2-y)=-\beta x + C$ (this serves only as a motivation for the
next inequality).  Let 
$$
\del''=\del \( 1 - \frac{2-\del}{2-\del^2} \beta' \).
$$
Since $y+\log y+\log(2-y)$ is
increasing on $(0,1/2]$, (3.14) gives
$$
\split
\del'+\log \del'&+\log(2-\del') \le \del'' + \log\del'' + \log(2-\del'') \\ 
&=\del+\log \del+\log (2-\del) - \frac{2\del-\del^2}{2-\del^2}\beta'
+ \log\left[\(1 - \frac{2-\del}{2-\del^2}\beta'\)\pfrac{2-\del''}{2-\del}
\right].
\endsplit
$$
Write 
$$
T= - \frac{2\del-\del^2}{2-\del^2}\beta'
+ \log\(1 - \frac{2-\del}{2-\del^2}\beta'\)
+\log\pfrac{2-\del''}{2-\del}.
$$
Using 
$$
\frac{2-\del''}{2-\del} = 1 + \frac{\del \beta'}{2-\del^2} \quad \text{ and }
\quad \log(1+x) \le x-\frac12 x^2 + \frac13 x^3,
$$
we obtain
$$
\split
T &\le -\beta'-\frac{(\beta')^2}{2(2-\del^2)^2} \( (2-\del)^2+\del^2 \) +
\frac{(\beta')^3}{3(2-\del^2)^3} \( -(2-\del)^3+\del^3 \)\\
&\le -\beta' - \frac{2}{5} (\beta')^2 \\
&\le -\beta-\frac25 \beta^2 + \frac{c(1+0.8\beta)}{\del}.
\endsplit
$$
The minimum of $\frac{(2-\del)^2+\del^2}{2(2-\del^2)^2}$ is actually
$0.401\ldots$.  Therefore
$$
\del'+\log \del'+\log(2-\del') \le \del+\log \del+\log(2-\del) - \beta
-0.4 \beta^2 + \frac{c(1+0.8\beta)}{\del}.
$$
Iteration of the above inequality yields
$$
\multline
\del_n+\log \del_n+\log(2-\del_n) \le \del_1+\log \del_1 +\log(2-\del_1)
-(n-1)(\beta+0.4\beta^2) \\
+c(1+1.6/k)\( \frac{1}{\del_1}+ \cdots + \frac{1}{\del_{n-1}} \).
\endmultline
$$
By (3.13) and (3.14),
$$
\del_{i+1} \le \del_i(1-\a), \quad \a = \frac67 (\beta-kc). \tag{3.18}
$$
By (3.13) again, this gives
$$
c(1+1.6/k) \( \frac{1}{\del_1}+ \cdots + \frac{1}{\del_{n-1}} \)
\le \frac{c(1+1.6/k)}{\a \del_{n-1}} \le \frac{1.34}{k}.
$$
Therefore,
$$
\del_n \le \frac{\del_1 (2-\del_1)e^{\del_1}}{(2-\del_n)e^{\del_n}}
e^{-(n-1)(\beta+0.4\beta^2)+1.34/k}. \tag{3.19}
$$
Next,
$$
\beta+0.4\beta^2 \ge \frac{2}{k}-\frac{32}{21k^2} + \frac{0.4}{k^2}
\(2 - \frac{32/21}{1000}\)^2 \ge \frac{2}{k}.
$$
From (3.13) and the inequality $1+x \le e^x$, we have
$$
\split
\del_1(2-\del_1)e^{\del_1} &= \frac{e^{1/2}}{2}\(1-\tfrac{1}{k} \) \( \tfrac32
+\tfrac{1}{2k} \) e^{-1/(2k)} \le \frac34 e^{1/2-7/(6k)}, \\
\frac{e^{-\del_n}}{2-\del_n} &\le \tfrac12 e^{\del_n /(2-\del_n)-\del_n}
\le \frac12 e^{\frac{0.49}{k}}.
\endsplit
$$
Putting these together with (3.19) gives
$$
\del_n \le \frac{3}{8} e^{1/2-2n/k+1.69/k}.
$$

To bound the constants $C_n$, take $\omega=0.06 > 1/(3\log k)$, so that
$$
V^{k+1} = (300k^3\log k)^{k+1} \le k^{4.11k} =: W.
$$
We next prove that
$$
W^{\Delta_{n-1}-\Delta_n} > k^{3k} 1.06^{4k(n-1)+k^2} \qquad (n \le 1.97k+1).
\tag{3.20}
$$
By (3.14),
$$
\del_{n-1}-\del_n \ge \frac{2\del_{n-1}}{k} \( \frac{2-\del_{n-1}}{2-\del_{n-1}
^2} - 0.002 \). \tag{3.21}
$$
By the top line of (3.16) and (3.17),
$$
\split
\del_m &\ge \del_{m-1}-\frac{2}{k} + \frac{4r}{2kr+y} \\
&\ge  \del_{m-1}-\frac{2}{k}
+4\, \frac{1-\del_{m-1}}{2k} = \del_{m-1}
\( 1 - \frac{2}{k} \),
\endsplit
$$
which implies
$$
\del_{n-1} \ge (1-2/k)^{n-2} \del_1 \ge \tfrac12 (1-2/k)^{n-1} \ge \tfrac12
e^{-\frac{2}{k-2}(n-1)} \ge 0.0096476 := \delbar.
$$
The right side of (3.21) is increasing in $\del_{n-1}$, so
$$
\del_{n-1}-\del_n \ge \frac{2\delbar}{k} \( \frac{2-\delbar}{2-\delbar^2}-
0.002 \) \ge \frac{0.01916}{k}.
$$
Therefore, $W^{\Delta_{n-1}-\Delta_n} \ge k^{0.0787 k^2}$.  On the other
hand, 
$$
k^{3k} 1.06^{4k(n-1)+k^2} \le k^{k^2(0.003+8.88\log(1.06)/\log(1000))}
\le k^{0.078k^2}.
$$
This proves (3.20).  Let $n_0=\lfloor 1.97k \rfloor +1$.  By (3.20) and
Lemma 3.5,
$$
C_{n_0} \le W^{\Delta_1 - \Delta_{n_0}} k! \le W^{\frac12 k^2 - \Delta_{n_0}}
$$
and for $n>n_0$
$$
C_n \le k^{3k} 1.06^{4k(n-1)+k^2} W^{\Delta_{n-1}-\Delta_n} C_{n-1}.
$$
Iterating this last inequality gives, for $n>n_0$,
$$
\split
C_n &\le W^{\frac12 k^2} k^{3k(n-n_0)} 1.06^{(n-n_0)k^2+4k(n_0+\cdots + n-1)}\\
&\le W^{\frac12 k^2} k^{3k(n-1.97k)} 1.06^{(n-1.97k)k^2+2k(n^2-n-(1.97k)^2+
  1.97k)}\\
&\le k^{2.055k^3-5.91k^2+3nk} 1.06^{nk^2+2(n^2-n)k-9.7278k^3}.
\endsplit
$$
This finishes the proof of Lemma 3.6.
\qed\enddemo
\bigskip

\demo{Proof of Theorem 3}  Suppose first that $k\ge 1000$.
Every
permissible $s$ can be written as $s=nk+u$ where $0\le u \le k$ and
$n\le \frac{k}{2} ( \frac12+\log \frac{3k}{8})$.  By Lemma 3.6 and \Ho,
$$
J_{s,k}(P) \le  k^{2.055k^3-5.91k^2+3s} 1.06^{sk+2s^2/k-9.7278k^3}
P^{2s-\frac12 k(k+1)+\Delta},
$$
where
$$
\Delta = \frac38 k^2 e^{1/2-2n/k+1.69/k} \left[ 1 - \tfrac{u}{k}+\tfrac{u}{k}
e^{-2/k} \right].
$$
Lastly,
$$
1 - \tfrac{u}{k}+\tfrac{u}{k} e^{-2/k} \le 1-\tfrac{2u}{k^2}+\tfrac{2u}{k^3}
\le e^{-2u/k^2+2u/k^3},
$$
thus $\Delta \le \frac38 k^2 e^{1/2-2s/k^2+1.7/k}$.

Next, suppose $129 \le k\le 1001$.  Start with $\Delta_1=\frac12 k^2(1-1/k)$,
successively choose $r_n$ near $\sqrt{k^2+k-2\Delta_n}$ satisfying (3.11),
and set $\Delta_{n+1} = \del_0(k,r_n,\Delta_n)$.  Also take $C_n$ as in
Lemma 3.5, where we define $\omega$ by
$$
\frac{1}{3\log k} \le \omega \le \frac{1}{2\log k + (4/3)\log\log k},
\quad e^{1.5+1.5/\omega} = \frac{18}{\omega} k^3\log k
$$
and take $\eta=1+\omega$.  To see that $\omega$ is well-defined, let
$h(\omega)=e^{1.5+1.5/\omega}-\frac{18}{\omega}k^3\log k$, $\omega_0=
\frac{1}{3\log k}$ and $\omega_1=1/(2\log k + (4/3)\log\log k)$.
It is easy to verify that $h(\omega_0)>0$, $h(\omega_1)<0$ and $h'(\omega)<0$
for $\omega\in [\omega_0,\omega_1]$.

If $\Delta_{n+1} \le \frac{k^2}{1000} \le \Delta_n$, take
$$
s = \left\lceil \( n + \frac{\Delta_{n} - k^2/1000}{\Delta_{n}-\Delta_{n+1}}
\) k \right\rceil.
$$
By H\"older's inequality,
$$
J_{s,k}(P) \le C_{n+1} P^{2s-\frac12 k(k+1)+0.001 k^2}.
$$
A straightforward computer computation verifies the claimed bounds on $s$
and $C_n$.  The program is listed in the Appendix.
\qed\enddemo 

{\bf Remarks.}  One can obtain slightly better values for $\Delta_n$ 
using a variant of the iterative scheme embodied in Lemmas 3.2 and 3.3.
For example, this alternate method would produce bounds valid with
$\rho = 3.20354$ for $129 \le k\le 199$.  The improvement, however,
becomes negligible for large $k$.  Instead of working with 
$K_s(P,Q;\bpsi;q)$, we work on bounding $K_{s,d}(P,Q;\bpsi;q)$, the
number of solutions of
$$
\aligned
&\sum_{i=1}^{k-d} (\Psi_j(z_i)-\Psi_j(w_i)) + q^j 
\sum_{i=1}^s (x_i^j-y_i^j)=0
\qquad (1\le j\le k), \\
&1 \le z_i,w_i \le P; \quad 1 \le x_i,y_i \le Q.
\endaligned
$$
Define $L_{s,d}(P,Q;\bpsi;p,q,r)$ similarly.  In Lemma 3.2, the variables
$z_{k-d+1},\ldots, z_k$ and $w_{k-d+1}, \ldots,w_k$ are not utilized
in the argument because $\Psi_j(z)=0$ for $j\le d$. 
Following the proof of Lemma 3.2 with the new quantities gives

\proclaim{Lemma 3.2'}  With the same hypotheses as Lemma 3.2,
$$
K_{s,d}(P,Q;\bpsi;q) \le 4k^3 k! p^{2s+\frac12(r-d)(r-d+1)} L_{s,d}(P,Q;
\bphi;p,q,r).
$$
\endproclaim

Likewise, following the proof of Lemma 3.3 and using H\"older's inequality
at the end gives

\proclaim{Lemma 3.3'}  Under the hypotheses of Lemma 3.3,
$$
\multline
L_{s,d}(P;Q;\bphi;p,q,r) \le (2P)^{k-d}\max\bigl[k^{k-d} J_{s,k}(Q), \\
2p^{-r(k-d)} J_{s,k}(Q)^{\frac{k-d-2}{2(k-d-1)}}
 K_{s,d+1}(P,Q;\bups;pq)^{\frac{k-d}{2(k-d-1)}} \bigr].
\endmultline
$$
\endproclaim

In Lemma 3.4, the definition of $\phi_J$ changes to
$$
\phi_J = \frac{1}{2r} + \frac{k^2+k+r^2-r-2\Delta-2rJ}{4r(k-J)}\phi_{J+1}
\quad (1\le J\le j-1),
$$
and this produces slightly smaller values for $\phi_1$.  The only
downside is that the analysis
of the numbers $\delta_n$ (see Lemma 3.6) becomes more complicated.
\vfil\eject

%
%
%
%
%
%
%
\head 4.  Incomplete systems and smooth Weyl sums. \endhead
%
%
%
%
%
%
%
%
\def\hjk{(h \le j \le k)}
\def\iu#1{\int_{\Bbb U^{#1}}}
\def\iut{\iu{t}}
\def\sums{\sum_{i=1}^s}
\def\sumt{\sum_{i=1}^t}

\def\TS{\widetilde S}
\def\GG{\curly G}

The object of this section is to obtain explicit upper bounds on
$J_{s,k,h}(\BB)$, the number of solutions of
$$
\sums (x_i^j-y_i^j) = 0 \qquad \hjk ; x_i,y_i \in \BB, \tag{4.1}
$$
where $\BB=\CPR=\{1\le n\le P: p|n \implies
\sqrt{R} < p\le R\}$.  Suppose $k\ge h\ge 2$ and set $t=k-h+1$.
For a $t$-tuple $\xx=(x_1,\cdots,x_t)$, let
$$
J(\xx) = \det \bigl( jx_i^{j-1} \bigr)_{1\le i\le t \atop h\le j\le k} =
\frac{k!}{(h-1)!} (x_1\cdots x_t)^{h-1} \prod_{1\le i<j\le t} (x_i-x_j)
\tag{4.2}
$$
be the Jacobian of the functions $\sum_{i=1}^t x_i^j$ $\hjk$.
The notation $x\DD(Q) y$ means that there is some $d|x$ with $d\le Q$
and $s_0(x/d)|s_0(y)$.
For $\balph = (\a_h,\cdots,\a_k)$, define the exponential sum
$$
f(\balph) = f(\balph;P,R) = \sum_{x\in \CPR} e(\a_h x^h+ \cdots + 
\a_k x^k)
$$
so that
$$
J_{s,k,h}(\CPR) = \iu{t} |f(\balph)|^{2s}\, d\balph.
$$

Our main lemma is very similar to the
the ``fundamental lemma'' (Lemma 3.1 of [34]).  However,
we do not perform ``repeat efficient differencing'' as in [29], [30], [34],
and Lemma 3.4 of this paper.

%
%
%

\proclaim{Lemma 4.1} Suppose
$$
\aligned
&k \ge h \ge 8, \quad t=k-h+1, \quad
s\ge t+1, \quad h\le r\le k; \\
&P > (8s)^{20}, \qquad R=P^\eta > k^2, \qquad |\CPR| \ge P^{1/2}.
\endaligned\tag{4.3}
$$
Then
$$
\multline
J_{s,k,h}(\CPR) \le \max \biggl[ \( (8s)^2 (22t^2)^{\frac{2}{\eta}}
P^{1/r} \)^{s-t} k^t |\CPR|^{s},
4k^{2t(\frac{1}{r\eta}+1)} \\
\times |\CPR|^t (P^{\frac{1}{r}} R)^{\frac12(r-h)(r-h+1)}
\biggl\{ \sum_{P^{\frac{1}{r}} < q \le P^{\frac1{r}}R} 
J_{s-t,k,h}(\CC(P/q,R))^{\frac{1}{2s-2t}} \biggr\}^{2s-2t} \biggr].
\endmultline
$$
\endproclaim

\demo{Proof}  For short, let $S_0=J_{s,k,h}(\CPR)$,
$\xx=(x_1,\ldots,x_t)$, $\yy=(y_1,\ldots,y_t)$ and
$\baa=(\a_h,\ldots,\a_k)$.  We divide
the solutions of (4.1) into four classes: $S_1$
counts the solutions with $\min(x_i,y_i) \le P^{1/5}$ for some $i$;
$S_2$ counts the solutions with $x_i=x_j$ or $y_i=y_j$ for some
 $1\le i<j\le t$;
$S_3$ counts solutions not counted by $S_1$ or $S_2$, and with 
$x_i\DD(P^{1/r})J(\xx)$ or $y_i\DD(P^{1/r})J(\yy)$ for some $i>t$;
$S_4$ (which will be the main term) counts the solutions not counted by
$S_1$, $S_2$ or $S_3$.

Evidently $S_0 \le 4\max(S_1,S_2,S_3,S_4)$.  If $S_1$ is the largest, then
by a trivial estimate and \Ho,
$$
\split
S_0 \le 4S_1 &\le 8s \iut |f(\baa)^{2s-1} f(\baa;P^{1/5},R)|
 \, d\baa \\
&\le 8s P^{1/5} \( \iut |f(\baa)^{2s}|\, d\baa \)^{1-\frac{1}{2s}} \\
&= 8s S_0^{1-1/2s} P^{1/5}.
\endsplit
$$
Therefore, $S_0 \le (8s P^{1/5})^{2s}$.  However, counting only the trivial
solutions of (4.1) (those with $x_i=y_i$ for every $i$) and using
(4.3) gives
$$
S_0 \ge |\CPR|^s \ge P^{s/2} > (8sP^{1/5})^{2s},\tag{4.4}
$$
giving a contradiction.

If $S_2$ is the largest, then by \Ho,
$$
\split
S_0 &\le 4S_2 \le 8 \binom{t}{2} \iut |f(\baa)^{2s-2} f(2\baa)|\, d\baa \\
&\le 4t^2 \(  \iut |f(\baa)|^{2s} d\baa \)^{1-\frac{1}{s}}
\( \iut  |f(2\baa)|^{2s}\, d\baa \)^{\frac{1}{2s}} \\
&= 4t^2 S_0^{1-\frac{1}{2s}}.
\endsplit
$$
By (4.3), $S_0 \le (4t^2)^{2s} < (8s)^{4s} < P^{s/2}$, contradicting (4.4).
It follows that $S_0 \le 4\max(S_3,S_4)$.

Suppose next that $S_3=\max(S_3,S_4)$.  From (4.2), we have $J(\xx)\ne 0$
and $J(\yy)\ne 0$ for each solution $(x_1,y_1,\ldots,x_s,y_s)$
 of (4.1) counted in $S_3$.
Let
$$
\curly S(\xx) = \{ w\in \CPR : w\DD(P^{1/r})J(\xx) \}
$$
and define
$$
H(\baa) = \sum_{\xx: J(\xx)\ne 0 \atop x_i \in \CPR} \sum_{w\in \curly S(\xx)}
e\( \sum_{j=h}^k \a_j(w^j+x_1^j+\cdots+x_t^j) \).
$$
By \CS,
$$
\split
S_0 \le 4S_3 &\le 8(s-t) \iut |H(\baa) f(\baa)^{2s-t-1}| \, d\baa \\
&\le 8s \( \iut |f(\baa)|^{2s}\,d\baa\)^{1/2} \( \iut |H^2(\baa)
f(\baa)^{2s-2t-2}|\, d\baa \)^{1/2} \\
&= 8s S_0^{1/2} \( \iut |H^2(\baa) f(\baa)^{2s-2t-2}|\, d\baa \)^{1/2}.
\endsplit
$$
Therefore,
$$
S_0 \le (8s)^2 \iut |H(\baa)^2 f(\baa)^{2s-2t-2}|\, d\baa,
$$
and the integral on the right is the number of solutions of
$$
\gather
\sum_{i=1}^{s-1} (x_i^j-y_i^j)+(dw)^j-(ez)^j=0 \qquad
\hjk \\
x_i,y_i \in \CPR; \; d,e\in \CC(P^{1/r},R); \; J(\xx)\ne 0, J(\yy)\ne 0; \\
w\in \CC(P/d,R), z\in \CC(P/e,R); \; s_0(w)|J(\xx), s_0(z)|J(\yy).
\endgather
$$
Writing
$$
\split
G_g(\baa) &= \sum_{\Sb \xx: J(\xx) \ne 0 \\ g|J(\xx) \endSb} e\( \sum_{j=h}^k
\a_j(x_1^j+\cdots+x_t^j) \), \\
\GG(\baa) &= \sum_{\Sb g\in \CPR \\ \mu^2(g)=1 \endSb} G_g(\baa) \sum_{d\in
\CC(P^{1/r},R)} \sum_{\Sb w\in \CC(P/d,R) \\ s_0(w)=g \endSb} e\( \a_h(dw)^h
+\cdots + \a_k(dw)^k \),
\endsplit
$$
it follows that
$$
S_0 \le (8s)^2 \iut | \GG(\baa)^2 f(\baa)^{2s-2t-2}|\, d\baa. \tag{4.5}
$$
By \CS,
$$
|\GG(\baa)|^2 \le \( \sum_g |G_g(\baa)|^2 \) \( \sum_g \biggl| \sum_{d,w} 1
\biggr|^2 \).
$$
Next,
$$
\split
\sum_g \biggl| \sum_{d,w} 1 \biggr|^2 &\le \sum_{g} \( P^{1/r}
 |\{w\le P: s_0(w)=g \}|
\) \sum_{w\in \CPR \atop g|w} \sum_{d\in \CC(P/w,R)} 1 \\
&\le P^{1/r} \sum_g |\{w\le P: s_0(w)=g \}| \sum_{n\in \CPR \atop g|n} d_2(n)
\\
&\le P^{1/r} \max_{g\in \CPR \atop \mu^2(g)=1} \left| \{w\le P:s_0(w)=g\}
\right| \sum_{n\in \CPR} d_2^2(n).
\endsplit
$$
For any $m\in \CPR$, $\tau(m) \le 2^{\Omega(m)} \le 2^{2/\eta}$.
Any $g\in \CPR$ with $\mu^2(g)=1$ can be written as $g=p_1\cdots p_n$,
where $p_1,\ldots,p_n$ are distinct primes each larger than $\sqrt{R}$,
and $0\le n\le 2/\eta$.  Then
$$
\split
 \left| \{w\le P:s_0(w)=g\} \right| &= \left| \{ \uu: u_1\log p_1 +
\cdots + u_n\log p_n \le \log P : u_i\ge 1\; \forall i \} \right| \\
&\le | \{ \uu: u_1+\cdots + u_n \le 2/\eta: u_i\ge 1\; \forall i\}| \\
&= \binom{\fl{2/\eta}}{n} < 2^{2/\eta}.
\endsplit
$$
Therefore,
$$
|\GG(\baa)|^2 \le 2^{6/\eta} P^{1/r} |\CPR| \sum_{g\in \CPR \atop \mu^2(g)=1}
|G_g(\baa)|^2,
$$
whence by (4.5),
$$
S_0 \le (8s)^2 2^{6/\eta} P^{1/r}|\CPR| V, \tag{4.6}
$$
where
$$
V= \iut  \sum_{g\in \CPR \atop \mu^2(g)=1}
\left| G_g(\baa)^2 f(\baa)^{2s-2t-2} \right|\, d\baa.
$$
Here $V$ counts the solutions $(x_1,y_1,\ldots,x_{s-1},y_{s-1},g)$ of
$$
\gather
\sum_{i=1}^{s-1} (x_i^j-y_i^j)=0 \qquad \hjk \\
x_i,y_i,g \in \CPR;\; J(\xx)\ne 0, J(\yy)\ne 0; \;
 \quad \mu^2(g)=1, g|J(\xx), g|J(\yy).
\endgather
$$
Clearly
$$
V \le J_{s-1,k,h}(\CPR) \max_{J(\xx)\ne 0} \left| \{ g\in \CPR, \mu^2(g)=1,
g|J(\xx) \} \right|.
$$
Using (4.2), $\sqrt{R} > k$ and $\mu^2(g)=1$,  $g|J(\xx)$ implies
$g| J^*(\xx)$, where
$$
J^*(\xx) = x_1 \cdots x_t \prod_{1\le i<j\le t} (x_i-x_j).
$$
Since $|J^*(\xx)| < P^{t(t+1)/2}$, $J^*(\xx)$ has at most $t(t+1)/\eta$
distinct prime factors $> \sqrt{R}$.  If $g|J^*(\xx)$, then $g$ is a product of
$n$ of these primes, where $0 \le n \le 2/\eta$.  The number of such $g$
is at most
$$
\sum_{0 \le n\le 2/\eta} \binom{\fl{(t^2+t)/\eta}}{n} \le
 \sum_{0 \le n\le 2/\eta} \frac{(2t^2/\eta)^n}{n!} \le
t^{4/\eta} \sum_{n=0}^\infty \frac{(2/\eta)^n}{n!} =
(et^2)^{2/\eta}.
$$
From (4.6) we conclude that
$$
S_0 \le (8s)^2 (8et^2)^{2/\eta} |\CPR| P^{1/r}  J_{s-1,k,h}(\CPR).
$$
Lastly, applying \Ho, we have 
$$
\split
J_{s-1,k,h}(\CPR) &\le J_{s,k,h}(\CPR)^{1-\frac{1}{s-t}} J_{t,k,h}(\CPR)^
{\frac{1}{s-t}} \\
&= S_0^{1-\frac{1}{s-t}} J_{t,k,h}(\CPR)^{\frac{1}{s-t}}. 
\endsplit
$$
We have $J_{t,k,h}(\CPR) \le k^t |\CPR|^t$, which follows for
instance from Lemma 2.4 (let $p$ be a prime $>tP$,
fix $y_1,\ldots,y_t$ and for each $\uu$ the
number of $\xx$ with $\sum x_i^j \equiv u_j \pmod{p} (h\le j\le k)$
is $\le k^t$).
This proves the lemma in the case $S_3 \ge S_4$.

\def\tp{\widetilde{p}}
\def\tq{\widetilde{q}}
\def\pp{\boldkey p}
\def\qq{\boldkey q}

For the last case, suppose $S_4=\max(S_1,S_2,S_3,S_4)$.  For every solution
of (4.1) counted by $S_4$, each $x_i> P^{1/r}$ and $y_i>P^{1/r}$ and
neither $x_i\DD(P^{1/r})J(\xx)$ nor $y_i\DD(P^{1/r})J(\yy)$ for $i>t$.
  Fix $i>t$
and let $q$ be the greatest divisor of $x_i$ with the property that
$(q,J(\xx))=1$.  If $q\le P^{1/r}$, then $x_i \DD(P^{1/r})J(\xx)$,
a contradiction.  Hence $q>P^{1/r}$, and since every prime divisor of $q$
is $\le R$, there is a divisor $q_i$ of $x_i$ with
$q_i > P^{1/r}$, $q_i \in \CC(P^{1/r}R,R)$ and $(q_i,J(\xx))=1$.
Likewise, each $y_i$  has a divisor $p_i$ with $p_i>P^{1/r}$, $p_i \in
\CC(P^{1/r}R,R)$ and $(p_i,J(\yy))=1$.
Therefore $S_0 \le 4T$, where $T$ is the number of solutions of
$$
\gather
\sumt (x_i^j-y_i^j) + \sum_{i=1}^{s-t} ((q_iu_i)^j-(p_iv_i)^j)=0 \qquad \hjk \\
x_i,y_i \in \CPR; \; u_i\in \CC(P/q_i,R), v_i\in \CC(P/p_i,R); \\
 p_i,q_i \in \CC(P^{1/r}R,R); p_i,q_i > P^{1/r}; (q_i,J(\xx))=(p_i,J(\yy))=1.
\endgather
$$
Let
$$
F_q(\baa) = \sum_{\xx: (q,J(\xx))=1} e\( \sum_{j=h}^k \a_j(x_1^j+\cdots
+x_t^j) \).
$$
 Given $q_1,p_1,\ldots q_{s-t},p_{s-t}$, let
$$
\tp = p_1 \cdots p_{s-t}, \qquad \tq = q_1 \cdots q_{s-t}
$$
and set
$$
\split
X_i(\baa) &= \left| F_{\tq}(\baa)^2 f((q_i^h\a_h,\cdots,q_i^k\a_k);P/q_i,R)
^{2s-2t} \right|, \\
Y_i(\baa) &= \left| F_{\tp}(\baa)^2 f((p_i^h\a_h,\cdots,p_i^k\a_k);P/p_i,R)
^{2s-2t} \right|. \\
\endsplit
$$
Then, by \Ho, we have
$$
\split
S_0 &\le 4 \sum_{\pp,\qq} \iut \prod_{i=1}^{s-t} \( X_i(\baa)Y_i(\baa)\)
^{\frac{1}{2s-2t}} \, d\baa \\
&\le 4\sum_{\pp,\qq} \prod_{i=1}^{s-t} \(\iut X_i(\baa)\,d\baa \)^{\frac{1}
{2s-2t}} \(\iut Y_i(\baa)\,d\baa \)^{\frac{1} {2s-2t}}.
\endsplit
$$
We have $\iut X_i(\baa)\, d\baa \le W(q_i)$ and  
$\iut Y_i(\baa)\, d\baa \le W(p_i)$,  where $W(q)$ is the number of solutions 
of
$$
\aligned
\sumt (x_i^j-y_i^j) + q^j \sum_{i=1}^{s-t} (u_i^j-v_i^j)=0 \qquad \hjk \\
x_i,y_i \in \CPR; \quad u_i,v_i\in \CC(P/q,R); \quad
(q,J(\xx)J(\yy))=1.
\endaligned\tag{4.7}
$$
Thus
$$
S_0 \le 4 \sum_{\pp,\qq} \prod_{i=1}^{s-t} \( W(q_i)W(p_i) 
\)^{\frac{1}{2s-2t}} 
= 4 \biggl( \sum_{\Sb q\in \CC(P^{1/r}R,R) \\ q>P^{1/r} \endSb} W(q)
^{\frac{1}{2s-2t}} \biggr)^{2s-2t}.
\tag{4.8}
$$
Next, by Proposition ZRD, for each possible $2t$-tuple $\xx,\yy$ in
(4.7), the number of $\uu,\vv$ is at most $J_{s-t,k,h}(\CC(P/q,R))$.
By fixing $\yy$, the number of possible $\xx,\yy$ is
$\le |\CPR|^t \max_{\mm} \BB(\mm)$, where $\BB(\mm)$ is the number of
solutions of the simultaneous congruences
$$
\sumt x_i^j \equiv m_j \pmod{q^j} \qquad \hjk
$$
with $1\le x_i\le P$ and $(q,J(\xx))=1$.
  For each $j$, the number of 
possibilities for $m_j$ modulo $q^r$ is $\max(1,q^{r-j})$.  Thus
$$
\BB(\mm) \le q^{(r-h)(r-h+1)/2} \max_{\nn} \BB'(\nn;q^r),
$$
where $\BB'(\nn;q^r)$ is the number of
solutions of
$$
\sumt x_i^j \equiv n_j \pmod{q^r} \qquad \hjk
$$
with $1\le x_i \le q^r$ (recall $q^r \ge P$)
 and $(q,J(\xx))=1$.  By the Chinese Remainder Theorem,
$$
\BB'(\nn;q^r) \le \prod_{p^\ell \| q, p\text{ prime}} \BB'(\nn;p^{r\ell}),
$$
and Lemma 2.4 gives $\BB'(\nn;p^{r\ell}) \le k!/(h-1)! \le k^t$.  
Since $\omega(q)
\le 2/(r\eta)+2$, we have $\BB'(\nn;q^r) \le k^{2t(1+1/(r\eta))}$.  This gives
$$
W(q) \le  k^{2t(1+1/(r\eta))}  q^{(r-h)(r-h+1)/2} |\CPR|^t J_{s-t,k,h}(\CC(P/q,R)).
$$
Together with (4.8), this proves the lemma in the fourth case.
\qed
\enddemo

%
%
%

The optimal choice for $r$ in the above lemma is close to $h$
for the range of $s$ that we are interested in.  The next lemma gives
some bounds achievable with Lemma 4.1.

\proclaim{Lemma 4.2}  Suppose that $k, h$ and $L$ are integers satisfying
$$
k\ge 60, \quad h\le k, \quad t=k-h+1 \le \tfrac{k}{6},
 \quad 1\le L\le h/2. \tag{4.9}
$$
Let $\a=1-1/h$.
Suppose $P,R$ and $\eta$ are real numbers with
$$
0 < \eta \le \frac{2}{3h}, \qquad 
R=P^\eta \ge \pfrac{2}{\eta}^3,\tag{4.10}
$$
and
$$
|\CC(Q,R)| \ge Q^{1/2} \qquad (P^{1/3} \le Q \le P). \tag{4.11}
$$
Then
$$
J_{Lt,k,h}(\CPR) \le (10\eta)^{tL((1/\eta+h) \a^{L-1}-h)} 
C_L (e^2 R)^{\frac{t}{2}L(L-1)}  P^{2Lt-\frac{t}{2}(h+k) + \Delta_L},
$$
where
$$
\split
\Delta_j &= \frac{t(t-1)}{2} + ht(1-1/h)^j \qquad (j\ge 1), \\
C_1 &= k^t,\quad C_\ell = \max_{2\le j\le \ell} e^{t  E_j} \quad (\ell\ge 2),\\
E_j &= \a^{L-j} \left[ \frac{4\log k}{\eta}(j-1) -
  \(j - \frac{j-1}{h} - h + h\a^j \) \log P \right].
\endsplit
$$
\endproclaim

\demo{Proof}  For $1\le j\le L$, define $P_j=P^{\a^{L-j}}$, $M_j=P^{\a^{L-j}}
R^{-h(1-\a^{L-j})}$, $\eta_j=\frac{\log R}{\log M_j}$ and $\eta_j' =
\frac{\log R}{\log P_j} = \a^{j-L}\eta$.
By (4.9) and (4.10), 
$$
M_j \ge P^{\a^L} R^{-h(1-\a^L)} \ge P^{0.6} R^{-0.4h} \ge P^{\frac13}
\ge (2/\eta)^{1/\eta} \ge (3h)^{75} > (8Lt)^{20}. \tag{4.12}
$$
Consequently, $\eta \le \eta_j' < \eta_j \le \eta_1 \le 3\eta$ for every 
$j$.
For $M\ge 1$ let $H_j(M) = J_{tj,k,h}(\CC(M,R))$.  We prove by induction on
$j$ that
$$
H_j(M) \le (10\eta)^{tj/\eta_1} C_j 
  (e^2 R)^{\frac{t}{2}j(j-1)}  M^{2jt-\frac{t}{2}(h+k) + \Delta_j}
\quad (M_j \le M \le P_j). \tag{4.13}
$$
By (4.10), $R \ge (3h)^3 > k^3 > 90000$.  By (4.11) and (4.12), when
$M_1 \le M\le P$ we have $|\CC(M,R)| \ge M^{1/2}$, so all of the hypotheses 
(4.3) of Lemma 4.1 hold (with $M$ in place of $P$).
Also, if $M_1 \le R^u \le P$ then
$R\ge (2/\eta)^3 \ge (2u)^3$, hence the hypotheses of Lemma 2.3 hold.
For $M\ge M_1$, as in the proof of Lemma 4.1 we have
$H_1(M) \le k^t|\CC(M,R)|^t$.  Writing $\nu=\frac{\log M}{\log R}$,
by Lemma 2.3
$$
|\CC(M,R)| \le M (2\nu)^{1/\nu} \le M (6\eta)^{1/\eta_1}, \tag{4.14}
$$
so (4.13) holds for $j=1$.  Next assume $j\ge 2$,
(4.13) holds with $j$ replaced by $j-1$, and assume $M_j \le M\le P_j$.
We will apply Lemma 4.1 with $r=h$ and $P=M$.  By the definition of
$M_j$ and $P_j$, 
$$
M_{j-1} \le M/q \le P_{j-1} \qquad \(P^{1/h} < q \le P^{1/h}R\).
$$
By (4.9) and (4.10),
$$
k(8jt)^2 (22t^2)^{2/\nu} \le \( k^{3/h} 22 t^2 \)^{2/\nu} < (27t^2)^{2/\nu}
< k^{4/\nu} \le k^{4/\eta_j'}.
$$
By Lemma 2.3 and $\nu \le 3\eta \le 2/h$,
$$
4k^{2t(\frac{1}{h\nu}+1)} |\CC(M,R)|^t \le 4M^t e^{\frac{t}{\nu}(\log(2\nu)
+ (2/h+2\nu)\log k)} \le 4M^t e^{\frac{t}{\nu}\log(3.13\nu)}.
$$
Since $e^{\frac{t}{\nu} \log(3.33/3.13)} \ge 4$, it follows that
$$
4k^{2t(\frac{1}{h\nu}+1)} |\CC(M,R)|^t \le M^t (10\eta)^{t/\eta_1}.
$$
By (4.14), Lemma 4.1 and the induction hypothesis,
$$
\split
H_j(M) &\le \max \biggl[(6\eta)^{tj/\eta_1} k^{\frac{4t(j-1)}{\eta_j'}}
  M^{tj+\frac{t(j-1)}{h}}, (10\eta)^{t/\eta_1} M^t \\
&\qquad \times
  \biggl\{ \sum_{M^{\frac1{h}} < q \le M^{\frac{1}{h}} R} H_{j-1}(M/q)
  ^{\frac{1}{2t(j-1)}} \biggr\}^{2t(j-1)} \biggr] \\
&\le (10\eta)^{tj/\eta_1} \max \biggl[
  k^{\frac{4t(j-1)}{\eta_j'}}  M^{tj+\frac{t(j-1)}{h}}, C_{j-1} \\
&\qquad \times (e^2 R)^{\frac{t}{2} (j-1)(j-2)}
  M^{2t(j-1)-\frac{t}{2}(h+k)+\Delta_{j-1}+t} S^{2t(j-1)} \biggr],
\endsplit
$$
where
$$
S =  \sum_{M^{\frac1{h}} < q \le M^{\frac{1}{h}} R}
q^{E}, \quad E= -1 + \frac{(t/2)(h+k) - \Delta_{j-1}}{2t(j-1)} =
-1+\frac{h(1-\a^{j-1})}{2j-2}.
$$
Making use of the inequalities
$$
1-\frac{\ell}{h} \le \a^\ell \le e^{-\ell/h} \le 1 - \frac{\ell}{h} +
\frac{\ell^2}{2h^2}, \tag{4.15}
$$
it follows that $-\frac58 \le E \le -\frac12$.  Thus
$$
S \le \int_1^{RM^{1/h}} x^E \, dx \le \frac{(RM^{1/h})^{E+1}}{E+1}
\le \frac83 (RM^{1/h})^{E+1} \le eR^{1/2} M^{(1-\a^{j-1})/(2j-2)}.
$$
We then obtain
$$
H_j(M) \le (10\eta)^{tj/\eta_1} \max \left[
k^{\frac{4t(j-1)}{\eta_j'}}   M^{tj+\frac{t(j-1)}{h}},
C_{j-1} (e^2 R)^{\frac{t}{2}(j^2-j)} M^{2tj-\frac{t}{2}(h+k)+
\Delta_j} \right].
$$
Write $f_j=j-\frac{j-1}{h} -h(1-\a^j)$, so that $f_1=f_2=0$
and $f_j>0$ for $j>2$.  Then
$$
H_j(M) \le  (10\eta)^{tj/\eta_1} M^{2tj-\frac{t}{2}(h+k)+\Delta_j}
\max \left[k^{\frac{4t(j-1)}{\eta_j'}}  M^{-tf_j}, C_{j-1}
(e^2 R)^{\frac{t}{2}(j^2-j)} \right].
$$
By (4.15),
$$
f_j \le j-\frac{j-1}{h} -h \( \frac{j}{h} -\frac{j^2}{2h^2}\) =
 \frac{j^2-2j+2}{2h} \le \frac{j^2-j}{2h} \quad (j\ge 2).
$$
Since $M \ge M_j \ge R^{-h} P^{\a^{L-j}}$, we have
$$
M^{-tf_j} \le R^{thf_j} P^{-tf_j \a^{L-j}} \le R^{\frac{t}{2}(j^2-j)}
P^{-tf_j\a^{L-j}}.
$$
Recalling the definition of $E_j$ and $\eta_j'$, we conclude that
$$
H_j(M) \le (10\eta)^{tj/\eta_1}
 (e^2 R)^{\frac{t}{2}(j^2-j)} M^{2tj-\frac{t}{2}(h+k)+\Delta_j}
\max\left[ e^{tE_j}, C_{j-1} \right].
$$
Since $e^{tE_2} > C_1$, (4.13) follows at once.  The Lemma then
follows from (4.13) by taking $j=L$.
\qed\enddemo

%
%

\proclaim{Lemma 4.3}  Suppose (4.9), (4.10) and (4.11)
hold, and define $E_j$ as in Lemma 4.2. Suppose that $\log P \ge A$ and
$$
x := \frac{4\log k}{A\eta\a} < 1.
\tag{4.16}
$$
Then
$$
\max_{j \ge 2} E_j \le \frac{4\log k}{\eta}\left[ 1 + h
 \( 1 + \frac{(1-x)\log(1-x)}{x} \) \right].
$$
\endproclaim

\demo{Proof}
We have $E_j \le \max_{z\ge 2} F(z)$, where
$$
F(z) = A(h-1/h-\a x + \a z(x-1) - h \a^z).
$$
By (4.16), $F(z) \to -\infty$ as $z\to \infty$ and $F(z)$ has a 
unique maximum point in $(-\infty,\infty)$. 
Solving $F'(y)=0$, we see that
$$
\a^y = \frac{\a (1-x)}{-h\log\a}. \tag{4.17}
$$
If $y<2$, then 
$$
\max_{z\ge 2} F(z)=F(2)=\frac{4\log k}{\eta}
$$
and the lemma follows in this case, because of the inequality
$(1-x)\log(1-x)\ge -x$.  Now assume $y\ge 2$.
Since $-h\log\a = 1+\frac{1}{2h}+\frac{1}{3h^2}+\cdots$, we have
$$
\frac{1}{1-\frac{1}{2h}} \le -h\log\a \le 1 + \frac{3h-1}{6h(h-1)}
\le 1 + \frac{1}{2h-2}.
$$
Consequently, by (4.17)
$$
x \ge \frac{1}{2h-1}. \tag{4.18}
$$
Also,
$$
\log(-h\log\a) \ge - \log\( 1- \frac{1}{2h} \) \ge \frac{1}{2h}+\frac{1}{8h^2}.
$$
This gives
$$
F(y) = A(h-1)(x+(1-x)V),
$$
where
$$
\split
V &= 1 - \frac{1}{-h\log\a} \(1 +
  \log(-h\log\a)-\log(1 - x) \) \\
&\le  1 - \frac{6h(h-1)}{6h^2-3h-1} \(1 +
  \frac{1}{2h}+\frac{1}{8h^2} \) + \frac{2h-2}{2h-1}\log(1-x) \\
&= \frac{5h+3}{4h(6h^2-3h-1)} + \frac{2h-2}{2h-1}\log(1-x) \\
&\le \frac{1}{4h^2}+\frac{2h-2}{2h-1}\log(1-x).
\endsplit
$$
Using $(1-x)\log(1-x) \ge -x$ again, we obtain
$$
\split
F(y) &\le \frac{(h-1)(1-x)A}{4h^2} + (h-1)Ax + \(1-\frac{1}{2h-1}\)A(h-1)(1-x)
  \log(1-x) \\
&\le \frac{(h-1)A}{4h^2} + (h-1)Ax \( \frac{1}{2h-1}-\frac{1}{4h^2} \) +
  (h-1)A(x+(1-x)\log(1-x)).
\endsplit
$$
By (4.18), we apply $1 \le (2h-1)x$ in the first summand to obtain
$$
\split
F(y) &\le (h-1)Ax \( \frac{2h-1}{4h^2} + \frac{1}{2h-1} - \frac{1}{4h^2}
  +1 + \frac{(1-x)\log(1-x)}{x} \) \\
&\le (h-1)Ax \( \frac{1}{h} +1 + \frac{(1-x)\log(1-x)}{x} \).
\endsplit
$$
The lemma now follows from the definition of $x$ (4.16).
\qed\enddemo

%
%
%
\bigskip
\demo{Proof of Theorem 4}
Let $L$ be an integer, $2 \le L \le h/2$, and put $R=P^\eta$ and $A=Dk^2$.
The hypotheses imply (4.9) and $\eta \le \frac{2}{3h}$.
Next, by (1.10),
$$
R \ge e^{\eta Dk^2} \ge k^{10} > \pfrac{2}{\eta}^3,
$$
so (4.10) holds.   Since $R \ge 6^{11}$, we may apply Lemma 2.2 with $\del=
\frac{1}{11}$.  Suppose $Q=P^\om$ with $\frac13 \le \om \le 1$ and put
$w = \fl{1.1\om/\eta}$.  Since  $m! \le m^m$ and $(w+1)\eta \le 1.1\om+
\eta \le 1.2$,
$$
|\CC(Q,R)| \ge \frac{11^{-w}}{(w+1)w!} \frac{Q}{\log R}
\ge \frac{1}{1.2} \pfrac{1}{11w}^w \frac{Q}{\log P}=Q^{\beta},
$$
where, by (1.10),
$$
\split
\beta &= 1-\frac{\log(1.2\log P)+w\log(11w)}{\log Q} \\
&\ge 1- \frac{3\log(1.2Dk^2)}{Dk^2} - \frac{1.1\log(12.1/\eta)}{\eta Dk^2} \\
&\ge 1-0.001 - 0.03 \ge 0.9.
\endsplit
$$
Thus, (4.11) holds and we may apply Lemmas 4.2 and 4.3.  By (1.10), (4.16)
and the bound $h\ge 54$,
$$
x = \frac{4h \log k}{Dk^2 \eta (h-1)} \in \left[ \frac{18}{k}, 0.408\right],
$$
so that
$$
1 + \frac{(1-x)\log(1-x)}{x} = \frac{x}{2} + \frac{x^2}{6} + \frac{x^3}{12} +
\cdots \le 0.5866 x.
$$
By Lemma 4.3,
$$
\split
\max_{j\ge 2} E_j &\le \frac{4\log k}{\eta} \( 1 + 0.5866hx \) \\
&\le 2.57 \frac{xk\log k}{\eta} \le 10.5 \frac{\log^2 k}{Dk\eta^2}.
\endsplit
$$
Therefore, by Lemma 4.2,
$$
J_{Lt,k,h}(\CPR) \le C_L (e^2 R)^{\frac{t}{2}L(L-1)} 
P^{2Lt-\frac{t}{2}(h+k)+\Delta_L},
$$
where $\Delta_L= \frac{t(t-1)}{2} + ht\a^L$ and
$$
\log C_L =  \frac{10.5 t \log^2 k}{Dk\eta^2} -tL \( \( \frac{1}{\eta}+h \)
\a^{L-1} - h \) \log\pfrac{1}{10\eta}.
$$
By hypothesis,
the number $s$ satisfies $s=Lt+u$, where $0\le u\le t$ and $2 \le L
< L+1 \le h/2$.  By \Ho,
$$
\aligned
J_{s,k,h}(\CPR) &\le \( J_{Lt,k,h}(\CPR) \)^{1-u/t} \( J_{Lt+t,k,h}(\CPR) \)^
{u/t} \\
&\le C_L^{1-u/t} C_{L+1}^{u/t}
 (e^2 R)^{\frac{t}{2}L^2+L(u-t/2)} P^{2s-\frac{t}{2}(h+k) + (1-u/t)\Delta_L
 + (u/t) \Delta_{L+1}}.
\endaligned\tag{4.19}
$$
Next,
$$
(1-u/t)\Delta_L + (u/t) \Delta_{L+1} =\frac{t(t-1)}{2} + ht \a^L
\( 1 - \tfrac{u}{ht} \) < \frac{t(t-1)}{2} + hte^{-s/(ht)}
$$
and 
$$
(e^2R)^{\frac{t}{2}L^2+L(u-t/2)} < (e^2R)^{s^2/(2t)} =e^{s^2/t} 
P^{\eta s^2/(2t)}.
$$
For the constants, we use $\a^{L-1} > \a^L \ge \a^{s/t}$.
Together with (4.19), this proves the theorem.
\qed
\enddemo

\vfil\eject
%
%
%
%

%
%
%
\head 5. Exponential Sums : Theorem 2 for large $\lambda$. \endhead
%
%
%
%
%
%
%
%
%

In this section, we apply Theorems 3 and 4 to  prove Theorem 2 for
large $\lambda$ ($\lam \ge 87$), using
a variant of Vinogradov's method to relate $S(N,t)$
to both $J_{r,k}(P)$ and $J_{s,g,h}(\BB)$.
Korobov's method [11] produces qualitatively
similar bounds, but does not have the seperation of variables
property (the $c_i,d_i$ below in Lemma 5.1), and therefore
one cannot easily modify it to incorporate incomplete systems (1.8).
Rough calculations indicate that Korobov's method, when combined with
Theorem 3, gives $S(N,t) \ll N^{1-1/(866\lam^2)}$.

\def\FM{\fl{M_1}}

\proclaim{Lemma 5.1}  Suppose $k$, $r$ and $s$ are integers $\ge 2$, and
$h$ and $g$ are integers satisfying
$1\le h \le g\le k$.  Let $N$ be a positive integer, and
$M_1$, $M_2$ be real numbers with $1\le M_i\le N$.
Let $\BB$ be a nonempty subset of the
positive integers $\le M_2$.
Then
$$
\multline
S(N,t) \le  2M_1M_2 + \frac{t(M_1M_2)^{k+1}}{kN^{k}} + \\
 N \pfrac{M_2}{|\BB|}^{\frac{1}{r}}\! \(\! (5r)^k M_2^{-2s} 
\FM^{-2r + \kk}
 J_{r,k}(\FM) J_{s,g,h}(\BB) W_h \cdots W_g \)^{\frac{1}{2rs}},
\endmultline
$$
where
$$
W_j = \min \( 2sM_2^j, \frac{2sM_2^j}{r\FM^j} + \frac{stM_2^j}{\pi j N^j}+
\frac{4\pi j (2N)^j} {rt\FM^j} + 2 \)
\qquad (j\ge 1).
$$
\endproclaim

\demo{Proof} For brevity write $M = \fl{M_1}$.
 For $N<R\le 2N$ and $0<u\le 1$, we have
$$
\split
\biggl| \sum_{N<n\le R}& (n+u)^{-it} \biggr|
  = \frac{1}{M |\BB|}
 \biggl| \sum_{\Sb a\le M_1 \\ b\in \BB \endSb}
 \sum_{N<n+ab\le R} (n+ab+u)^{-it} \biggr| \\
&\le  \frac{1}{M|\BB|} \biggl| \sum_{\Sb a\le M_1 \\ b\in \BB \endSb }
\sum_{N<n\le R-1}  (n+ab+u)^{-it} \biggr| + \frac{1}{M|\BB|}
\sum_{\Sb a\le M_1 \\ b\in \BB \endSb } (2ab-1) \\
&\le \frac{N}{M|\BB|}
 \max_{N\le z\le 2N} \biggl| \sum_{\Sb a\le M_1 \\ b\in \BB \endSb}
 e^{-it\log(1+ab/z)} \biggr| + 2M_1M_2.
\endsplit
$$
For $0 \le x \le 1$ we have
$$
\bigl| \log(1+x) - (x-x^2/2 + \cdots + (-1)^{k-1} x^k/k)\bigr| 
\le \frac{x^{k+1}}{k+1}. \tag{5.1}
$$
Also $|e^{iy}-1| \le y$ for real $y$ and $ab/z \le M_1M_2/N$.
Thus, for some $z\in [N,2N]$,
$$
S(N,t) \le \frac{N}{M|\BB|} |U| + 
\frac{t(M_1M_2)^{k+1}}{(k+1) N^k} + 2M_1M_2,
\tag{5.2}
$$
where $U=\sum_{a,b} e(\g_1 (ab) + \cdots + \g_k (ab)^k)$ and
$\gamma_j=(-1)^jt/(2\pi j z^j)$.
By H\"older's inequality,
$$
\split
|U|^r &\le |\BB|^{r-1} \sum_{b\in \BB} \left| \sum_{a\le M_1} 
e(\g_1 (ab) + \cdots + \g_k (ab)^k) \right|^r \\
&= |\BB|^{r-1} \sum_{b\in \BB} \e_b \( \sum_{a\le M_1}
 e(\g_1 (ab) + \cdots + \g_k (ab)^k) \) ^r \\
&= |\BB|^{r-1} \sum_{b\in \BB} \e_b \sum_{c_1,\ldots,c_k} n(\bc)
 e(\g_1 b c_1 + \cdots + \g_k b^k c_k),
\endsplit
$$
where $\e_b$ are complex numbers with $|\e_b|=1$, and for
$\bc=(c_1,\ldots,c_k)$, $n(\bc)$ is the number of solutions of
the simultaneous equations
$c_j=a_1^j+\cdots+a_r^j$ $(1\le j\le k)$ with each $a_i \in [1,M_1]$.
A second application of  H\"older's inequality gives
$$
\aligned
|U|^{2rs} &\le |\BB|^{2rs-2s} \( \sum_{\bc} n(\bc) \)^{2s-2} \( \sum_{\bc}
n(\bc)^2 \) T \\
&=  |\BB|^{2rs-2s} M^{2rs-2r} J_{r,k}(M) T,
\endaligned\tag{5.3}
$$
where
$$
T = \sum_{\bc} \left| \sum_{b\in \BB} \e_b e(\g_1 bc_1 + \cdots + \g_k b^k c_k)
\right|^{2s}.
$$
For $0<w\le \frac12$, let $\ell(x;w)=\max(0,1-\frac{\| x \|}{w})$.
This function has an absolutely and uniformly convergent
Fourier series
$$
\ell(x;w) = \frac{1}{\pi^2 w} \sum_{n=-\infty}^{\infty} \pfrac{\sin \pi n w}{n}
^2 e(nx).
$$
For $1\le j\le k$ define
$$
f_j(x) = \pfrac{rM^j \sin(\pi x/(2rM^j))}{x}^2,
$$
and we note that $f_j(x)\ge 0$ for all $x$ and $f_j(x) \ge 1$
for $1\le x\le rM^j$.  Since $1\le c_j \le rM^j$ for each $j$, we have
$$
\split
T &\le \sum_{\Sb \bc \\ -\infty < c_j < \infty \endSb}
\left| \sum_{b\in \BB} \e_b e(\g_1 bc_1 + \cdots + \g_k b^k c_k)
\right|^{2s}f_1(c_1) \cdots f_k(c_k) \\
&=  \sum_{\Sb \bc \\ -\infty < c_j < \infty \endSb}
\sum_{\Sb b_1,\cdots,b_{2s} \\ b_i\in \BB \endSb}
 \e_{\bb} e(\g_1 d_1 c_1 + \cdots + \g_k d_k c_k)
f_1(c_1)\cdots f_k(c_k),
\endsplit
$$
where $|\e_{\bb}|=1$ and $d_j = b_1^j + \cdots + b_s^j - b_{s+1}^j -\cdots -
b_{2s}^j$ for $1\le j\le k$.  For $\dd=(d_1,\ldots,d_k)$, write
 $J_{s,n,m}(\BB;\dd)$ for the number of
$\bb$ with $b_i\in \BB$ for each $i$ and
 $d_j = b_1^j + \cdots + b_s^j - b_{s+1}^j -\cdots -
b_{2s}^j$ $(m\le j\le n)$.  By Proposition ZRD,
$J_{s,n,m}(\BB;\dd) \le J_{s,n,m}(\BB)$.  Then
$$
\split
T &\le \sum_{d_1,\ldots,d_k} J_{s,k,1}(\BB;\dd) \left| \sum_{\bc} 
  e(\g_1 d_1 c_1 + \cdots + \g_k d_k c_k)f_1(c_1)\cdots f_k(c_k) \right| \\
&=  \sum_{\dd}  J_{s,k,1}(\BB;\dd)  \prod_{j=1}^k \left|
  \sum_{c=-\infty}^{\infty} e(cd_j\g_j) f_j(c) \right| \\
&= \sum_{\dd} J_{s,k,1}(\BB;\dd) \prod_{j=1}^k 
  \( (rM^j)^2 \frac{\pi^2}{2rM^j} \ell(d_j\g_j;\tfrac{1}{2rM^j}) \)  \\
&= (\pi^2 r/2)^k M^{\kk} \sum_{\dd} J_{s,k,1}(\BB;\dd) \prod_{j=1}^k
\ell(d_j\gamma_j; \tfrac{1}{2rM^j}).
\endsplit
$$
Recalling the definition of $\ell(x;w)$, we obtain
$$
T \le (5r)^k M^{\kk} \sum_{d_j\in \DD_j\, \forall j}
J_{s,k,1}(\BB;\dd), \tag{5.4}
$$
where
$$
\DD_j = \{ |d_j| < sM_2^j-1: \| d_j\g_j \| < \tfrac{1}{2rM^j} \}. 
$$
The sum in (5.4) may be interpreted as the number of solutions of the
system of equations
$$
\sum_{i=1}^s (x_i^j-y_i^j)=d_j \quad (1\le j\le k); \;\; x_i,y_i \in \BB;
d_j \in \DD_j. \tag{5.5}
$$
There are now several ways to proceed.  A simple method is to ignore the
equations in (5.5) corresponding to $j>g$ or $j<h$.  Then, by Proposition
ZRD, for each choice of $d_h,\ldots,d_g$, the number of $\xx,\yy$ is
$\le J_{s,g,h}(\BB)$.  Thus, by (5.4),
$$
T \le (5r)^k M^{\kk} J_{s,g,h}(\BB) \prod_{j=h}^g |\DD_j|.
$$
An alternate and slightly better method for bounding the number of solutions
of (5.5) will be given in \S 8.
Lastly, for positive $\delta$, $\gamma$ and $K$, we claim that
$$
| \{ |d| \le K: \| d\gamma \| < \delta \}| \le 4K\delta+2K\gamma+
4\delta/\gamma  +2. \tag{5.6}
$$
Suppose that $\delta<1/2$, else (5.6) is trivial.
The number of intervals of the form $[m-\delta,m+\delta]$ with
integral $m$ which intersect
$[-K\gamma,K\gamma]$ is $\le 2\gamma K+1+2\delta \le 2\gamma K+2$.
Each such interval can contain at most $2\delta/\gamma+1$ points of the
form $d\gamma$, and this proves (5.6).  Putting $K=sM_2^j-1$,
$\gamma=|\gamma_j|$ and $\delta=\frac{1}{2rM^j}$ gives 
$|\DD_j| \le W_j$, hence
$$
T \le  J_{s,g,h}(\BB)  (5r)^k M^{k(k+1)/2} W_h \cdots W_g.
$$
Together with (5.2) and (5.3), this proves the lemma.
\qed\enddemo

{\it Proof of Theorem 2 for $\lam\ge 87$.}
Assume that
$$
\FM \ge M_2 \ge 100g, \quad s\le 2^g, \quad r\ge 13 g, \quad r\ge s,
\quad g\ge h\ge 3.
\tag{5.7}
$$
It turns out that the optimal parameters satisfy (5.7).  By (5.7) and the
definition of $W_j$,
$$
W_j \le 4 + \frac{stM_2^j}{\pi N^j} + \frac{13 g 2^g N^j}{rtM_1^j}
\le 2^{g+1} \max \( 1, \frac{tM_2^j}{N^j}, \frac{N^j}{tM_1^j} \).
$$
Suppose that
$$
M_1 = N^{\mu_1}, \quad M_2 = N^{\mu_2}, \quad \mu_1 > \mu_2. \tag{5.8}
$$
Then, the above bound for $W_j$ is better than the trivial bound
$2sM_2^j$ only when $\lambda < j < \lambda/(1-\mu_1-\mu_2)$.  Let
$$
\phi = g/\lm, \quad \gamma=h/\lm, \qquad 1 \le \gamma \le \frac{1}
{1-\mu_2} < \frac{1}{1-\mu_1} \le \phi \le \frac{1}{1-\mu_1-\mu_2}.
\tag{5.9}
$$
We then have
$$
W_h \cdots W_g \le 2^{g^2} M_2^{h + (h+1) + \cdots + g} N^{-H},
\tag{5.10}
$$
where
$$
H = \sum_{j=h}^g \min \( j\mu_2, j-\lm, \lm-j(1-\mu_1-\mu_2) \). \tag{5.11}
$$
For $i=1,2$, write $\frac{\lm}{1-\mu_i} = m_i + \beta_i$, where $m_i$ is 
an integer and $0\le \beta_i < 1$. Then
$$
\aligned
H &= \sum_{j=h}^{m_2} (j-\lm)+ \sum_{j=m_2+1}^{m_1} j\mu_2 + \sum_{j=m_1+1}^g
   (\lm-j(1-\mu_1-\mu_2)) \\
&= \frac{(m_1^2+m_1)(1-\mu_1)+ (m_2^2+m_2)(1-\mu_2)
  -h^2+h-(1-\mu_1-\mu_2)(g^2+g)}{2} \\
&\qquad +  \lm(h+g-m_1-m_2-1) \\
&= \lm^2 \( \phi+\gamma-\frac{\gamma^2}{2}- \frac{1-\mu_1-\mu_2}{2} \phi^2
  -\frac{2-\mu_1-\mu_2}{2(1-\mu_1)(1-\mu_2)} \) \\
&\qquad  + \lm \( \frac{\gamma}{2}
  -\frac{\phi}{2}(1-\mu_1-\mu_2)\) - \frac{\beta_1(1-\beta_1)(1-\mu_1)+
  \beta_2(1-\beta_2)(1-\mu_2)}{2}.
\endaligned
$$ 
Since $\beta_i(1-\beta_i) \le \frac14$,
$$
\aligned
H &\ge  \lm^2 \( \phi+\gamma-\frac{\gamma^2}{2}- \frac{1-\mu_1-\mu_2}{2} \phi^2
  -\frac{2-\mu_1-\mu_2}{2(1-\mu_1)(1-\mu_2)}  \) \\
&\qquad + \lm \( \frac{\gamma}{2}  -\frac{\phi}{2}(1-\mu_1-\mu_2)\)
 - \frac{2-\mu_1-\mu_2}8 \\
&=: H_2 \lam^2 + H_1\lam - H_0.
\endaligned\tag{5.12}
$$
We shall take the near-optimal choice for the parameters
$$
\aligned
\mu_1 &=0.1905, \quad \mu_2=0.1603, \quad k= \bigfl{\tfrac{\lm}{1-\mu_1
  -\mu_2}+0.000003} \ge 129, \\
r &= \fl{\rho k^2+1}, \qquad \rho \text{ taken from (1.7)},
\endaligned\tag{5.13}
$$
and approximate values (to be specified precisely later)
$$
g \approx 1.2453\lam, \qquad h\approx 1.1818\lam, \qquad 
s\approx 0.3299h(t-1).
$$
With these choices
we quickly deduce that $S(N,t) \ll N^{1-1/(132.31\lm^2)}$
for sufficiently large $\lm$.  By a standard argument (see \S 7), this
implies (1.1) with $B=4.42736$, but only for $1-\sigma$ sufficiently small.
For completely explicit bounds, we pay more attention to the constants,
sacrificing a little bit in $B$ in order to get a fairly small value
for $A$ in Theorem 1.

By (5.13) and Theorem 3, we have
$$
\fl{M_1}^{-2r+\kk} J_{r,k}(\fl{M_1})\le C_1 M_1^{0.001k^2}, \tag{5.14}
$$
where $C_1 = k^{\theta k^3}$ and $\theta$ is taken from (1.7).
Let $Y=300$ and assume that 
$$
N \ge e^{Y \lm^2}, \tag{5.15}
$$
for otheriwse trivially
$$
S(N,t) \le N \le e^{Y/133.66} N^{1-1/(133.66\lam^2)} \le 9.44
N^{1-1/(133.66\lam^2)}.
$$
We shall always choose $g$ so that
$$
106 \le g \le 1.254 \lam. \tag{5.16}
$$
Thus by (5.13) and (5.15), 
$M_2 \ge e^{\mu_2 Y \lam^2} \ge e^{0.1019Y g^2}$.
Let $D=0.1019Y=30.57$ and $\eta=\frac{1}{\xi g^{3/2}}$, where $3\le \xi \le 6$.
By (5.16), (1.10) holds and
hence the hypotheses of Theorem 4 hold (with  $P=M_2$
 and $k=g$).  By Theorem 4,
$$
J_{s,g,h}(\CC(M_2,M_2^\eta)) \le C_2 P^{2s-\frac{t}{2}(h+g) + E_2},
\tag{5.17}
$$
where
$$
\aligned
E_2 &= \frac12 t(t-1) + \frac{\eta s^2}{2t} + ht\exp\{ -\frac{s}{ht} \}, \\
\log C_2&=\frac{s^2}{t} + \frac{10.5\xi^2 t g^2 \log^2 g}{D}
 - s\( (\xi g^{3/2}+h)(1-1/h)^{s/t}-h \) \log(\xi g^{3/2}/10).
\endaligned\tag{5.18}
$$
By (1.10) and (5.16),
$$
R = M_2^{\eta} \ge e^{Dg^2 \eta} \ge g^{10} > 6^{26}.
$$
By Lemma 2.2 (with $\del=\frac{1}{26}$) plus the inequality
$w! \le (w/2.5)^w$ ($w\ge 50$), we have
$$
\split
\frac{M_2}{|\CC(M_2,R)|} &\le (\log R) (1.04\xi g^{3/2}+1) 
\pfrac{27.04\xi g^{3/2}}{2.5}^{1.04\xi g^{3/2}} \\
&\le (\log N) C_3 \le C_3 N^{E_3}, 
\endsplit
$$
where
$$
\aligned
C_3 &= (10.82\xi g^{3/2})^{1.04\xi g^{3/2}}, \\
E_3 &= \frac{\log(Y\lam^2)}{Y\lam^2}.
\endaligned\tag{5.19}
$$
By (5.13),
$$
(5r)^k \le (40\lam^2)^{1.6\lam} \le \lam^{5\lam} \tag{5.20}
$$
and
$$
r \ge 7.509 \lam^2. \tag{5.21}
$$
Consequently
$$
\frac{E_3}{r} \le \frac{\log(Y\lam^2)}{7.5Y\lam^4}.
$$
By Lemma 5.1, (5.10), (5.13), (5.14), (5.17) and (5.20), it follows that
$$
\aligned
S(N,t) &\le \( C_3^{\frac{1}{r}} \( \lam^{5\lam} C_1 C_2
  \)^{\frac{1}{2rs}} \) 
  N^{1+E} +2N^{0.36} + \frac{1}{k} N^{1-0.0000019476}, \\
E &=  \frac{\log(Y\lam^2)}{7.5Y\lam^4} + \frac{1}{2rs} 
\( -H + 0.001\mu_1 k^2 + \mu_2 E_2 \).
\endaligned\tag{5.22} 
$$
We also need bounds on $k/\lam$, which by (5.13) can be written as
$$
k_0 := \tfrac{1}{0.6492} - \tfrac{0.999997}{\lam}
\le \tfrac{k}{\lam} \le \tfrac{1}{0.6492} + \tfrac{0.000003}{\lam} =: k_1.
\tag{5.23}
$$

%
%

\proclaim{Lemma 5.2} When $\lam \ge 220$, we have
$$
S(N,t) \le 7.5 N^{1-1/(133.58\lam^2)} \qquad (N\ge e^{300\lam^2}).
$$
\endproclaim

\demo{Proof}
We take
$$
h=\fl{1.1818\lam+\tfrac12}, \quad g=\fl{1.2453\lam+\tfrac12},
\quad s = \fl{\sigma h(t-1)+1}, \quad \sigma=0.3299. \tag{5.24}
$$
By (5.9) and (5.24), (5.16) holds and also
$$
|\gamma-1.1818| \le \tfrac{1}{2\lam}, \qquad
|\phi-1.2453| \le \tfrac{1}{2\lam}. \tag{5.25}
$$
Further, by (5.13) and (5.24),
$$
g \ge 274, \quad h\ge 260, \quad t\ge 13, \quad k\ge 338, \quad
s\ge 0.02294\lam^2. \tag{5.26}
$$
By (1.7), (5.13) and (5.14),
$$
C_1 = k^{2.3291 k^3} \le e^{9.2 \lam^3\log \lam}. \tag{5.27}
$$
Taking
$$
\xi = 6,
$$
we have by (5.19) and (5.24),
$$
C_3 \le e^{20.31 \lam^{3/2}\log \lam}.
\tag{5.28}
$$
To bound $C_2$, we first note that by (5.24),
$$
(1-1/h)^{s/t} \ge (1-1/h)^{\sigma (h-1)} \ge e^{-\sigma}\ge 0.71899.
$$
This implies
$$
(\xi g^{3/2} + h)(1-1/h)^{s/t}-h \ge 5.9785\lam^{3/2}-0.28101h \ge 
5.956\lam^{3/2}.
$$
By (5.18), (5.24) and (5.26),
$$
\aligned
\log C_2 &\le 0.3907 s\lam + 20.86 t\lam^2 \log^2 \lam - 8.73s \lam^{3/2}
\log \lam \\
&\le 1.52\lam^3\log^2\lam - 8.72 s \lam^{3/2}\log \lam.
\endaligned\tag{5.29}
$$
By (5.21) and (5.26), $2rs\ge 0.3445\lam^4$.
Combining (5.21), (5.27), (5.28) and (5.29), we obtain 
$$
\aligned
C_3^{\frac{1}{r}}  \( \lam^{5\lam} C_1C_2\)^{\frac{1}{2rs}} &\le \exp \biggl\{
\frac{\log \lam}{\lam^{1/2}}\( \frac{20.31-8.72/2}{7.509} \)
 + \frac{\log\lam}{0.3445} \(
\frac{5}{\lam^3}+\frac{9.2+1.52\log \lam}{\lam} \) \biggr\} \\
&\le e^{2.011} \le 7.48.
\endaligned\tag{5.30}
$$
By (5.22), it remains to bound $E$.  Note that $-H+0.001\mu_1 k^2 + \mu_2
E_2 < 0$.
 By (5.9), (5.13), (5.18) and (5.22),
$$
\split
E &\le \frac{\log(Y\lam^2)}{7.5Y\lam^4} + \frac{-H+0.001\mu_1 k^2}
  {2.002\rho \sigma \gamma (\phi-\gamma) \lm^2 k^2} 
  + \frac{\mu_2 E_2}{2\rho k^2 s} \\
&\le \frac{1.52\times 10^{-7}}{\lam^2} + \frac{-\lam^2 H_2 -\lam H_1 + H_0}
  {2.002\rho \sigma \gamma (\phi-\gamma) \lm^2 k^2} +
  \frac{0.001\mu_1}{2.002\rho \sigma \gamma (\phi-\gamma) \lm^2} \\
&\qquad + \frac{\mu_2}{2\rho k^2} \left[ \frac{\phi-\gamma+\tfrac{1}{\lam}}
{2\sigma \gamma} + \frac{\frac{t}{t-1} e^{-\sigma+\sigma/t}}{\sigma}
 +\frac{\sigma h g^{-3/2}}{12} \right].
\endsplit
$$
By (5.26), 
$$
\tfrac{t}{t-1} e^{\sigma/t} \le 1 + \frac{1.33413}{t-1} = 1 + 
\frac{1.33413}{(\phi-\gamma)\lm}.
$$
Therefore
$$
\lm^2 E \le 1.52\times 10^{-7} + \frac{f(\gamma,\phi) +
 G_1/\lam^{1/2} + G_2/\lam}{\rho},
$$
where, by (5.24) and (5.25), 
$$
\split
f(\gamma,\phi) &= \frac{1}{2.002\sigma \gamma} \left[ \frac{0.001 \mu_1}{\phi
  -\gamma} +\frac{1}{k_1^2}\(\frac{-H_2}{\phi-\gamma}+
  1.001\mu_2 \( \frac{\phi-\gamma}{2}+\gamma e^{-\sigma} \)\) \right], \\
G_1 &= \frac{\mu_2 \sigma \gamma \phi^{-3/2}}{24 k_0^2} \le 0.0008, \\
G_2 &= \frac{1}{2.002\sigma(k/\lam)^2} \left[ 
  \frac{-H_1+H_0/\lam+1.33547\mu_2 \gamma e^{-\sigma}}{\gamma(\phi-\gamma)}
  +\frac{1.001\mu_2}{2\gamma} \right].
\endsplit
$$
Let $U$ be the bracketed expression in the definition of $G_2$.
By (5.12) (the definition of $H_1$ and $H_0$), (5.25) and (5.26),
$$
\split
U &\le \frac{0.3246 \phi - 0.34608 \gamma+0.20615/\lam}
  {\gamma(\phi-\gamma)} +\frac{1.001\mu_2}{2\gamma} \\
&= \frac{1.001\mu_2 + 0.6492}{2\gamma} + \frac{-0.02148+\frac{0.20615}{h}}
  {\phi-\gamma} \\
&\le \frac{0.80967}{2.3636-1/\lam}+\frac{-0.02148+0.20615/h}{0.0635+1/\lam} \\
&\le 0.0392.
\endsplit
$$
Thus
$$
G_2 \le  \frac{0.0392}{2.002\sigma \g k_0^2}\le 0.021334.
$$
Then
$$
\aligned
\lam^2 E &\le 1.52\times 10^{-7} + \frac{f(\gamma,\phi)
+ 0.0008 \lam^{-1/2} + 0.021334 \lam^{-1}}{\rho} \\
&\le 0.00004711 + \frac{f(\gamma,\phi)}{\rho}.
\endaligned\tag{5.31}
$$
A short analysis with the aid of Maple shows that in the range
$|\phi-1.2453|\le \frac{1}{440}$, $|\g-1.1818|\le \frac{1}{440}$, we have
$$
f(\gamma,\phi) \le  -0.0242145,
$$
the maximum occuring at $\gamma=1.1818+\frac{1}{440}$, $\phi=1.2453-
\frac{1}{440}$.  By (1.7), (5.13) and (5.31), we conclude that
$$
\lam^2 E \le -0.0074862 \le -\frac{1}{133.58}.
$$
Together with (5.22) and (5.30), this proves the lemma.
\qed\enddemo

\proclaim{Lemma 5.3} When $87 \le \lam \le 220$, we have 
$$
S(N,t) \le 8.4 N^{1-1/(133.66\lam^2)} \qquad (N \ge e^{300\lam^2}).
$$
\endproclaim

\demo{Proof} Here we take 
$$
\xi = 3.6, \qquad s= \fl{\sigma ht}+1, \quad \sigma=0.3299.
$$
We choose $g,h$ satisfying (5.16) and 
$$
g=\bigfl{\tfrac{\lam}{1-\mu_1}}+1+a, \quad h=\bigfl{\tfrac{\lam}{1-\mu_2}}-b,
\quad t=g-h+1, \quad a,b\in\{0,1\}.
$$
To bound the exponent of $N$,
consider $\lam\in I=[\lam_1,\lam_2)$,
a small interval on which each of the quantities 
$m_1=\fl{\frac{\lam}{1-\mu_1}}$,
$m_2=\fl{\frac{\lam}{1-\mu_2}}$ and $k$ (defined in (5.13)) is constant. 
We choose constant values of $a$ and $b$ in $I$, so that $g,h,t,s,r$ are
also fixed.
By the definition of $H$, we  have for $\lam\in I$
$$
\split
H &= Z_0 + Z_1 \lam, \\
Z_0 &= \frac{(m_1^2+m_1)(1-\mu_1)+ (m_2^2+m_2)(1-\mu_2)
  -h^2+h-(1-\mu_1-\mu_2)(g^2+g)}{2}, \\
Z_1 &= h+g-m_1-m_2-1 = a-b \in \{-1,0,1\}.
\endsplit
$$
Therefore,
$$
H \ge H' := Z_0 + \cases \lam_1 & Z_1=1 \\ 0 & Z_1=0 \\ -\lam_2 & Z_1=-1 
\endcases.
$$
By (5.22),
$$
E \le \frac{\log(Y\lam_1^2)}{7.5Y\lam_1^4} - \frac{H'-0.001\mu_1 k^2 -
\mu_2 \( \frac{t(t-1)}{2} + \frac{s^2}{\xi tg^{3/2}}+hte^{-s/(ht)} \)}{2rs}
:= E'.
$$
Then, by (5.22), when $\lam\in I$ we have
$$
S(N,t) \le C N^{1-1/(u\lam^2)} + \frac{1}{k} N^{1-1/(133\lam^2)},
$$
where $u = 1/(E' \lam_1^2)$ and $C=C_3^{1/r} (\lam^{5\lam} C_1 C_2)^{1/(2rs)}$.
A short computer program (Program 2 in the Appendix) is used to compute
$C$ and $u$ in each interval, and to find the
best choice for $a$ and $b$ (the choice which gives the smallest $C$
subject to $u\le 133.66$).  In all cases, $C\le 8.38$.
For most $\lambda$, we take $b=0$ and for $\lam\in [136,220]$ we take
$a=1$.
This concludes the proof.

No choice of parameters $g,h,s$ produced $C < 9.5$ in the range
$86 \le \lam\le 87$.
\qed\enddemo

Together, Lemma 5.2 and 5.3 prove Theorem 2 for $\lam \ge 87$.

\vfil\eject

%
%
%
\head 6. Theorem 2 for small $\lambda$ \endhead
%
%
%
%

We begin with a general inequality derived from the Weyl shifting method.
Suppose $N$ is a positive integer and $M$ is a  real number satisfying
$1 \le M\le N$.
Arguing as in the proof of Lemma 5.1,
for $N<R\le 2N$ and $0<u\le 1$, we have
$$
\split
\biggl| \sum_{N<n\le R}& (n+u)^{-it} \biggr|
  = \frac{1}{\fl{M+1}}
  \left| \sum_{m\le M+1}
  \sum_{N<n+m\le R} (n+m+u)^{-it} \right| \\
&\le  \frac{1}{M} \left| \sum_{m\le M}
  \sum_{N<n\le R-1}  (n+m+u)^{-it} \right| + \frac{N}{M}+\frac{1}{M}
\sum_{m\le M}  (2m-1).
\endsplit
$$
Therefore,
$$
S(N,t) \le \frac{1}{M} \max_{0<u\le 1}
 \sum_{N< n\le 2N-1} \left| \sum_{m\le M}
 e^{-it\log(1+m/(n+u))} \right| + \frac{N}{M}+M. \tag{6.1}
$$

\proclaim{Lemma 6.1} If $|\alpha - p/q| \le 1/q^2$, $(p,q)=1$,
$m$ is a positive integer, and $x\ge 1$ and
$y\ge 2$ are real numbers, then
$$
\sum_{n\le x} \min \( y, \frac{1}{2\| \alpha m n \|}\)
\le \( 1 + \frac{2mx}{q}\) \(2q \log(ey)+ 4y \).
$$
\endproclaim

\demo{Proof} For $0\le j\le 2q-1$ let $I_j$ be the interval
$[\frac{j}{2q},\frac{j+1}{2q})$.  The interval
$[1,x]$ can be partitioned into intervals $B_i$, $1\le i\le 1+2mx/q$, each of
length $\le q/(2m)$.  If $n,n' \in B_i$ and $\{\alpha mn\}, \{ \alpha mn'\}
\in I_j$ then
$$
\left\| \frac{pm}{q}(n-n') \right\| \le \| \alpha mn - \alpha mn' \| + \left|
\frac{m(n-n')}{q^2} \right| < \frac{1}{q},
$$
hence $n=n'$.  So, for $0 \le j\le q-1$,
 there are at most $G = 2 + 4mx/q$ values of
$n$ giving $\| \alpha mn\| \in I_j$.  We take the summand to be
$y$ when $j \le q/y+1$, thus
$$
\sum_{n\le x} \min \( y, \frac{1}{2\| \alpha m n \|}\) \le Gy(q/y+2)+ G
 \sum_{q/y+1< j\le q-1} q/j \le  G (q+2y+q\log y).
 \text{\qed}
$$
\enddemo

%
%
%
%

Next, we use the Weyl method to prove Theorem 2 for $1\le \lam\le 2.6$.
There is much room for improvement here, but the bounds below more than
suffice for our purposes.

\proclaim{Lemma 6.2} 
We have
$$
\split
S(N,t) &\le  5 N^{1-1/20} \qquad (1\le \lam \le 1.9), \\
S(N,t) &\le 30 N^{1-1/83} \qquad(1.9\le \lam \le 2.6).
\endsplit
$$
Consequently, when $1\le \lam\le 2.6$, we have
$$
S(N,t) \le 1.81 N^{1-1/(133\lam^2)}.
$$
\endproclaim

\demo{Proof} Suppose $k\ge 2$.
By (6.1) and (5.1), for some real number $z\in [N,2N]$,
$$
S(N,t) \le \frac{N}{M} |U| + \frac{N}{M}+M +\frac{tM^{k+1}}{(k+1)N^k},
\tag{6.2}
$$
where
$$
U = \sum_{m\le M} e^{-it((m/z)-m^2/(2z^2) + \cdots + (-1)^{k-1} m^k/
(kz^k))}.
$$
By the proof of Weyl's inequality (e.g. Lemma 2.4 of [27]), we have
$$
|U|^{2^{k-1}} \le (2M)^{2^{k-1}-k} \sum_{\Sb h_1,\ldots, h_{k-1} \\
|h_i| \le M-1 \endSb} \min\( M, \frac{1}{2 \| \a h_1 \cdots h_{k-1} k! \|}\),
$$
where $\a = t/(2\pi k z^k)$.  There are at most $(k-1)(2M)^{k-2}$\;
vectors $(h_1,\cdots,h_{k-1})$ with some $h_i=0$, thus
$$
\aligned
|U|^{2^{k-1}} &\le (2M)^{2^{k-1}-k} \biggl( (k-1)M^{k-1}2^{k-2} \\
&\qquad + 2^{k-1} \sum_{1\le h\le M^{k-1}} d_{k-1}(h) \min 
\bigl( M, \frac{1}{2\| \a h k! \|} \bigr) \biggr).
\endaligned\tag{6.3}
$$

Suppose $1\le \lam \le 1.9$.  Let $q=\fl{1/\a}$ and note that
$\frac{(4\pi-1) N^2}{t} \le q \le \frac{16\pi N^2}{t}$.
Assume $M\ge 10000$.   By (6.3) with $k=2$ and Lemma 6.1,
$$
\aligned
|U|^2 &\le 9M + \frac{32M^2}{q}+(16M+4q)\log(eM) \\
&\le \frac{32M^2t}{(4\pi-1) N^2}+ \( 17M+\frac{64\pi N^2}{t} \)\log(eM).
\endaligned\tag{6.4}
$$
We may assume that $N \ge 5^{20}$, otherwise the claimed bound
is trivial.  We shall take $M=N^\mu$, where $\mu=\frac{2.95-\lambda}{3}
\in [0.35,0.65]$, so that $M\ge 5^7 > 10000$.
By (6.2) and (6.4),
$$
\split
S(&N,t) \le N \( \frac{32t}{(4\pi-1) N^2}+\( \frac{17}{M}+\frac{64\pi N^2}
{tM^2}  \)\log(eM) \)^{1/2}+\frac{N}{M}+M+\frac{tM^3}{3N^2} \\
&\le N \( 3N^{\lambda-2}+\( 17 N^{-0.35}+64\pi N^{-0.3} \)\log N
  \)^{1/2} + 2N^{0.65} + \frac{N^{0.95}}{3} \\
&\le N \( 3N^{-0.1} + 205 N^{-0.3} \log N \)^{1/2} + 0.334 N^{0.95} \\
&\le 4.1 N^{0.95}.
\endsplit
$$

When $1.9\le \lam\le 2.6$, we apply (6.3) with $k=3$, obtaining
$$
|U|^4 \le 2M \( 4M^2 + 4 \sum_{1\le h\le M^2} d_2(h) \min \bigl( M,
\frac{1}{2\| 6 \a h \|} \bigr) \).
$$
We shall use a crude upper bound on $\tau(h)$:
$$
\split
\frac{d_2(h)}{h^{1/3}} &= \prod_{p^e \| h} \frac{e+1}{p^{e/3}}
\le \prod_p \max_{e\ge 0} \frac{e+1}{p^{e/3}} \\
&= \prod_{p\le 7}  \max_{e\ge 0} \frac{e+1}{p^{e/3}} =
\frac{24}{315^{1/3}} \le 3.53.
\endsplit
$$
Take $q=\fl{\frac{6\pi z^3}{t}}$, so that $\frac{(6\pi-1) N^3}{t}\le q
\le \frac{48\pi N^3}{t}$.  By Lemma 6.1 (with $m=6$, $x=M^2$, $y=M$), we obtain
$$
\aligned
|U|^4 &\le 8M^3 + 28.24 M^{5/3}  \sum_{1\le h\le M^2} \min \bigl( M,
   \frac{1}{2\| 6 \a h \|} \bigr) \\
&\le 8M^3 + 28.24 M^{5/3} \( (2q+24M^2) \log(eM) +4M +
\tfrac{48M^3}{q} \).
\endaligned\tag{6.5}
$$
We choose $\mu=1-\frac{\lam+1/50}{4}\in [0.345,0.52]$ and put $M=N^\mu$.
Then
$$
3-\lam = 3-4(1-\mu)+\tfrac1{50} = -\tfrac{49}{50} + 4\mu \in \bigl[ \mu+0.055,
\mu(2+\tfrac3{26}) \bigr],
$$
and consequently
$$
M N^{0.055} \le \frac{N^3}{t} \le M^{2+3/26}.
$$
We assume that $N\ge 30^{60}$, otherwise the claimed bound is trivial.
Then $M\ge 30^{20.7}$ and
$$
N^{0.055} \ge 74000, \qquad \frac{\log(eM)}{M^{1/13}} \le 0.318.
$$
Thus,
$$
\split
(2q+24M^2)&\log(eM) + 4M + \frac{48M^3}{q} \le (96\pi M^{2+3/26}+24M^2)
  \log(eM) \\
&\hskip 2.5in +4M+\frac{48}{6\pi-1} \frac{M^2}{N^{0.055}} \\
&\le M^{2+5/26} \( (96\pi+24 M^{-3/26}) (0.318)+\frac{4}{M^{1+5/26}}+\frac{1}
  {26000M^{5/26}} \) \\
&\le 96 M^{2+5/26}.
\endsplit
$$
By (6.5), $|U|^4 \le 2712 M^{4-\frac{11}{78}}$, and (6.2) then gives
$$
S(N,t) \le 4N(7.22M^{-\frac{11}{312}}) + 2 N^{0.655}+\tfrac14 N^{1-\frac
{1}{50}} \le 30 N^{1-1/83}.
$$
This completes the proof of the first part of the lemma.
The last part follows a general inequality: if $\lam$
is fixed and $0<d<c<1$, then
$$
S(N,t) \le C N^{1-c} \quad (N\ge 1) \quad \implies \quad S(N,t) \le
C^{d/c} N^{1-d} \quad (N\ge 1). \tag{6.6}
$$
For the proof, if $N\le C^{1/c}$, then trivially 
$S(N,t)\le N = N^d N^{1-d} \le C^{d/c} N^{1-d}$.  When $N>C^{1/c}$,
the hypothesis of (6.6) implies that
$$
S(N,t) \le C N^{1-c} = CN^{d-c} N^{1-d} \le C \cdot C^{\frac{1}{c}(d-c)}
 N^{1-d} = C^{d/c} N^{1-d}.
$$
For $\lam\in [1,1.9]$, take $c=\frac{1}{20}$,
$d=\frac{1}{133}$ in (6.6) and for $\lam\in [1.9,2.6]$ take $c=\frac{1}{83}$,
$d=\frac{1}{133.66(1.9^2)}$.
\qed
\enddemo

%
%
%
%

For larger $\lam$, we relate $S(N,t)$ to $J_{s,k}(P)$ using
an older method (\S 6.12 of [25]).

%
%

\proclaim{Lemma 6.3}  Suppose $k\ge 2$, $s\ge 2$, $N\ge 1$,
$1\le M\le N t^{-\frac{1}{k+1}}$ and $t\le N^k$.   Then
$$
S(N,t) \le \frac{4 N^{1-\frac{1}{2s}}}{M} \( \pi^k k! k^k W M^{\kk} J_{s,k}(M)
\)^{\frac{1}{2s}} + \frac{N}{M}+M,
$$
where
$$
W = \frac{2^{k+2}N^{k+1}}{k^2 tM^k}+1.
$$
\endproclaim

\demo{Proof}
By (6.1) and \Ho,
$$
S(N,t) \le \max_{0<u\le 1} \frac{N^{1-\frac1{2s}}}{M} \( \sum_{N<n\le 2N-1}
|T(n)|^{2s} \)^{\frac1{2s}} + \frac{N}{M}+M, \tag{6.7}
$$
where
$$
T(n) = \sum_{m\le M} e\( -\frac{t}{2\pi} \log\(1+\frac{m}{n+u}\) \).
$$
With $n$ fixed, let $\g_j = \g_j(n)= (-1)^j \frac{t}{2\pi j(n+u)^j}$ for
$1\le j\le k$. 
Define
$$
\split
S(x;\bbeta) &= \sum_{m\le x}
 e(m\beta_1 + \cdots+ m^k \beta_k), \\
\delta(m;\bbeta)&=  -\frac{t}{2\pi} \log\(1+\frac{m}{n+u}\) - \sum_{j=1}^{k}
\beta_j m^j.
\endsplit
$$
When $0\le w \le M$,
$$
\aligned
|\delta'(w;\bbeta)| &\le \frac{tM^k}{2\pi N^{k+1}}+\sum_{j=1}^k j|\beta_j-\g_j|
M^{j-1} \\
&\le \frac{1}{2\pi M} + \sum_{j=1}^k j|\beta_j-\g_j| M^{j-1}.
\endaligned
\tag{6.8}
$$
Let $\Omega_n$ be the region $\{ \bbeta: |\beta_j-\g_j| \le \frac{1}{2\pi jk
M^j} \; \forall j \}$.   By (6.8),
 for $\bbeta \in \Omega_n$ and $0 \le w\le M$,
$|\delta'(w;\bbeta)| \le \frac{1}{\pi M}$.  For any $\bbeta\in \Omega_n$, 
 partial summation gives
$$
T(n) = S(M;\bbeta)e(\delta(\fl{M};\bbeta))-2\pi i \int_0^M S(w;\bbeta) 
e(\delta(w;\bbeta)) \delta'(w;\bbeta)\, dw,
$$
and thus
$$
|T(n)| \le  |S(M;\bbeta)| + \frac{2}{M}  \int_0^M |S(w;\bbeta)|\, dw =:
 S_0(\bbeta).
$$
Integrating over $\Omega_n$ then gives
$$
|T(n)|^{2s} \le \frac{1}{|\Omega_n|} \int_{\Omega_n} S_0(\bbeta)^{2s}\,
d\beta = \pi^k k! k^k M^{\kk} \int_{\Omega_n} S_0(\bbeta)^{2s}.
$$
For any $\bbeta$, the number of $n$ with $\bbeta \in \Omega_n$ 
is at most the number of $n$ with $|\g_k(n)-\beta_k| 
\le \frac{1}{2\pi k^2 M^k}$.
By hypothesis, $|\g_k(N)-\g_k(2N)| < \frac12$ and by the mean value theorem,
when $N\le n\le 2N-2$,
$$
|\g_k(n)-\g_k(n+1)| \ge \frac{t}{2\pi (2N)^{k+1}}.
$$
Therefore the number of such is $n$ is at most $W$.
Hence
$$
\sum_{N<n\le 2N-1} |T(n)|^{2s} \le \pi^k k! k^k M^{\kk} W \iuk
 S_0(\bbeta)^{2s} \, d\bbeta. \tag{6.9}
$$
By \Ho,
$$
\split
S_0(\bbeta)^{2s} &\le 2^{2s-1} \( |S(M;\bbeta)|^{2s}+ \pfrac{2}{M}^{2s}
\(\int_0^M |S(w;\bbeta)|\, dw \)^{2s} \) \\
&\le 2^{2s-1} |S(M;\bbeta)|^{2s} + \frac{2^{4s-1}}{M} \int_0^M 
 |S(w;\bbeta)|^{2s} \, dw.
\endsplit
$$
Thus
$$
\split
\iuk S_0(\bbeta)^{2s}\, d\bbeta &\le  2^{2s-1} \iuk
 |S(M;\bbeta)|^{2s}\, d\bbeta + \frac{2^{4s-1}}{M} \int_0^M \iuk
 |S(w;\bbeta)|^{2s} \,d\bbeta\, dw \\
&= 2^{2s-1} J_{s,k}(M) + \frac{2^{4s-1}}{M} \int_0^M J_{s,k}(w)\, dw \\
&\le 2^{4s} J_{s,k}(M).
\endsplit
$$
Combined with (6.7) and (6.9), this gives the lemma.
\qed\enddemo

%
%
%

\proclaim{Corollary 6.4}  Suppose $k\ge 4$, $1\le n\le k^2$ and $s=nk$.
Assume that 
$$
J_{s,k}(P) \le C P^{2s-\kk+\Delta} \qquad (P\ge 1).
$$
Then, for $t=N^\lam$ with $k-1\le \lam \le k$, we have
$$
S(N,t) \le 4N^{1-\frac1{2nk}} \( C (2\pi k)^k k! 
N^{\frac{2+2\Delta}{k+1}} \)^{
\frac{1}{2nk}} + N^{\frac{k}{k+1}} + N^{\frac{1}{k+1}}.
$$
\endproclaim

\demo{Proof}  This follows from Lemma 6.3, taking $\mu=1-\frac{\lam}{k+1} \in
[\frac{1}{k+1}, \frac{2}{k+1}]$, $M=N^\mu$ and noting that $W\le 2^kM$.
\qed
\enddemo

Bounds for $J_{s,k}(P)$ with the best exponents of $P$ come from Lemma 3.5,
however the constants are very large.  By using older methods without
``repeat efficient differencing'', we obtain bounds with far better
constants, while sacrificing something in the exponents of $P$.  In fact,
using Corollary 6.4 with the older bounds for $J_{s,k}(P)$ gives
$$
S(N,t) \le C_\lam N^{1-1/16\lam^2} \quad (6\le \lam \le 100),
$$
which is far better than needed for Theorem 2.  Since we will then use
(6.6) to greatly reduce the constant (to $C_\lam^{16/133.66}$),
it is better for us to minimize $C_\lam$ rather than the exponent of $P$.
Lemma 6.5 below comes from using Lemma 3.2 in a non-iterative way.  For
some $s$, even better constants can be obtained using an older
variation of the method (Lemma 6.6), where solutions modulo a single
prime are considered (as opposed to considering a set of $k^3$ primes).

%
%

\proclaim{Lemma 6.5}
Suppose $k\ge 4$ and $1\le n \le k^2$.  Suppose $0 < \omega\le \frac12$
or $\omega=1$, and let $\eta=1+\omega$.  Put
$V(\om)=6k^3\log k$ if $\omega=1$ and 
$V(\om)=\max(e^{1.5+1.5/\omega},\tfrac{18}{\omega} k^3\log k)$ otherwise.  If
$$
J_{nk,k}(P) \le C P^{2nk-\kk+\Delta} \qquad (P\ge 1),
$$
then
$$
J_{nk+k,k}(P) \le C' P^{2n k -\kk + \Delta'},
$$
where $\Delta' = (1-1/k)\Delta$ and
$$
C' = C \max \left[ 4k^3 k! \eta^{k^2-\Delta}, 
V(\om)^{\Delta} \right].
$$
\endproclaim

\demo{Proof} 
This comes from Lemma 3.2 with $Q=P$, $M=P^{1/k}$, $r=k$, $d=0$, 
$T=1$, $s=nk$, $q=1$ and $\psi_j(z)=z^j$ for each $j$.
Lemma 2.1 implies that the interval $[P^{1/k},\eta P^{1/k}]$ contains at
least $k^3$ primes.
Also, $P^{1-1/k} \ge k^{3k-3}\ge 32s^2$, so the
hypotheses of Lemma 3.2 are satisfied.  Together with the inequality
$$
L_s(P,P/p;\bphi;p,1,k) \le P^k J_{s,k}(P/p),
$$
this proves that for $P\ge V(\om)^k$ and
$k\le s \le k^3$, there is a  prime $p\in [P^{1/k},\eta P^{1/k}]$ giving
$$
J_{s+k,k}(P) \le 4k^3 k! p^{2s+\frac12(k^2-k)} P^k J_{s,k}(P/p).
\tag{6.10}
$$
The upper bound on $p$ now gives the lemma.
If $P<V(\om)^k$, trivially
$$
\split
J_{(n+1)k,k}(P) &\le P^{2k} J_{nk,k}(P) \le C P^{2(n+1)k-\kk+\Delta}\\
&\le C V(\om)^{k(\Delta-\Delta')} P^{2s-\kk+\Delta'}.
\endsplit
$$
\qed\enddemo

%
%
%

\proclaim{Lemma 6.6}  Suppose $k\ge 9$, $k\le s\le k^3-k$, $P\ge k^k$ and
$p$ is a prime in $[P^{1/k}, 2P^{1/k}]$.  Then
$$
J_{s+k,k}(P) \le  \max\left[ (ep)^{2k-2}(k-1)^{2s+2}J_{s+k,k}(\tfrac{P}{p}),
\tfrac{32}{k!} (s+k)^{2k} p^{2s+\frac12 (k^2-k)} P^k J_{s,k}(\tfrac{P}{p})
\right].
$$
\endproclaim

\demo{Proof}  Let $S_1$ be the number of solutions of (1.4) (with $s \to s+k$)
with at least $k$ distinct residues modulo $p$ among $x_1,\ldots,x_{s+k}$
or  at least $k$ distinct residues modulo $p$ among $y_1,\ldots,y_{s+k}$.
Let $S_2$ be the number of remaining solutions.  Clearly $J_{s+k,k}(P) \le
2\max(S_1,S_2)$.  Let
$$
f(\baa;Q) = \sum_{x\le Q} e(\a_1 x + \cdots + \a_k x^k).
$$
If $S_2\ge S_1$, for $1\le b\le p$ let
$$
\split
g(\baa;b) &= \sum_{\Sb x\le P \\ x\equiv b\pmod{p} \endSb} e(\a_1 x + \cdots
+ \a_k x^k) \\
&= \sum_{1\le y\le \frac{P+p-b}{p}} e(\a_1 (py+b-p) + \cdots
+ \a_k (py+b-p)^k).
\endsplit
$$
Define $\curly B$ to be the set of $(b_1,\ldots,b_{s+k})$
with $1\le b_i \le p$ for each $i$, and containing at most $k-1$
distinct values.  Then
$$
|\curly B| \le \binom{p}{k-1} (k-1)^{k+s} \le\frac{p^{k-1}}{(k-1)!} 
(k-1)^{k+s} \le \tfrac12 (ep)^{k-1} (k-1)^{s+1}.
\tag{6.11}
$$
By \Ho,
$$
\split
J_{s+k,k}&(P) \le 2 \iuk \left| \sum_{\bb \in \curly B} g(\baa;b_1) \cdots
  g(\baa;b_{s+k}) \right|^2 \, d\baa \\
&\le 2 \iuk \( \sum_{\bb,\bb' \in \curly B} |g(\baa;b_1)|^{2s+2k} \)^{\frac{1}
  {2s+2k}}\!\! \cdots \( \sum_{\bb,\bb' \in \curly B} 
  |g(\baa;b'_{s+k})|^{2s+2k}\)^{\frac{1} {2s+2k}} \!\! d\baa \\
&= 2\sum_{\bb,\bb' \in \curly B} \iuk |g(\baa;b_1)|^{2s+2k}\, d\baa.
\endsplit
$$
By the binomial theorem, the last integral is $J_{s+k,k}(\frac{P+p-b_1}{p})$,
hence
$$
J_{s+k,k}(P) \le 2 |\BB|^2 J_{s+k,k}(P/p+1).
$$
For brevity, write $P_1=P/p+1$.  We have $J_{s+k,k}(P_1)\le 2\max(S_3,S_4)$,
where $S_3$ is the
number of solutions of (1.4) with every $x_i,y_i \le P/p$ and
$S_4$ is the number of remaining solutions.  If $S_4\ge S_3$,  \Ho\, implies
$$
\split
J_{s+k,k}(P_1) &\le 2(2s+2k) \iuk |f(\baa;P_1)|^{2s+2k-1}\, d\baa \\
&\le 4(s+k) \( \iuk |f(\baa;P_1)|^{2s+2k}\,d\baa \)^{1-\frac{1}{2s+2k}} \\
&= (4s+4k) J_{s+k,k}(P_1)^{1-\frac{1}{2s+2k}},
\endsplit
$$
whence $J_{s+k,k}(P_1) \le (4s+4k)^{2s+2k}$.
On the other hand, since $k\ge 9$ and $P\ge k^k$,
counting only trivial solutions gives
$$
J_{s+k,k}(P_1) \ge (P/p)^{s+k} \ge (\tfrac12 P^{1-1/k})^{s+k} >
(4k^3)^{2s+2k} \ge (4s+4k)^{2s+2k},
$$
a contradiction.  Therefore, $J_{s+k,k}(P_1) \le 2J_{s+k,k}(P/p)$,
and by (6.11),
$$
J_{s+k,k}(P) \le 4 |\curly B|^2 J_{s+k,k}(P/p) \le 
 (ep)^{2k-2} (k-1)^{2s+2} J_{s+k,k}(P/p).
$$
This proves the lemma in the case $S_2\ge S_1$.

If $S_1\ge S_2$, then $S_1$ is at most $2\binom{s+k}{k}$ times the 
number of solutions of (1.4) with $x_1,\cdots,x_k$ distinct modulo $p$.
Let $\curly X$ be the set of $k$-tuples $(x_1,\cdots,x_k)$ which are
distinct modulo $p$ and
$$
F(\baa) = \sum_{\xx\in\curly X} e(\a_1 (x_1+\cdots+x_k)+ \cdots + \a_k (
x_1^k+\cdots+x_k^k)).
$$
Then, by \CS,
$$
\split
J_{s+k,k}(P) &\le 2S_1 \le 4\binom{s+k}{k} \iuk | F(\baa) f(\baa;P)^{2s+k}|
\, d\baa \\
&\le  4\binom{s+k}{k} \( \iuk |F(\baa)^2 f(\baa;P)^{2s}|\, d\baa \)^{\frac12}
 \( \iuk |f(\baa;P)^{2s+2k}|\, d\baa \)^{\frac12} \\
&=  4\binom{s+k}{k} \( \iuk |F(\baa)^2 f(\baa;P)^{2s}|\, d\baa \)^{\frac12}
\( J_{s+k,k}(P) \)^{1/2}.
\endsplit
$$
Thus
$$
J_{s+k,k}(P) \le  16\binom{s+k}{k}^2 \iuk |F(\baa)^2 f(\baa;P)^{2s}|\, d\baa
=  16\binom{s+k}{k}^2 S_3(p),
$$
where $S_3(p)$ is defined as in the proof of Lemma 3.2 (with $\Psi_j(x)=x^j$
for $j=1,\ldots,k$).   All the hypotheses
of Lemma 3.2 hold, with $d=0$, $T=1$, $M=P^{1/k}$, $r=k$, $Q=P$ and
$q=1$.  Recalling the definition (3.2) of $L_s(P,Q;\bpsi;p,q,r)$
and using (3.7),
$$
\split
J_{s+k,k}(P) &\le  16\frac{(s+k)^{2k}}{(k!)^2} 2k! p^{2s+\frac12k(k-1)}
 L_s(P,P/p;\bphi;p,1,k)\\
&\le \tfrac{32}{k!}(s+k)^{2k} p^{2s+\frac12k(k-1)} P^k J_{s,k}(P/p),
\endsplit
$$
and the lemma follows in the case $S_1\ge S_2$.
\qed\enddemo

The chief advantage of Lemma 6.6 over Lemma 6.5 is the much smaller
lower bound required for $P$ (see (6.10)).

%
%

\proclaim{Lemma 6.7}  Suppose $k\ge 9$, $1\le n\le k^2$ and
$$
J_{nk,k}(P) \le C P^{2nk-\kk+\Delta} \qquad (P\ge 1).
$$
Suppose that $1<\eta \le 2$ and that
for $x\ge k$, there is a prime in $[x,\eta x]$.
Then
$$
J_{(n+1)k,k}(P) \le C' P^{2nk-\kk+\Delta'}, \tag{6.12}
$$
where $\Delta'=(1-1/k)\Delta$ and 
$$
\split
C' &= C \max \left[ U^{\Delta}, \tfrac{32}{k!} (nk+k)^{2k} 
\eta^{k^2-\Delta} \right], \\
U &= \max \left[ k, \left\{ e^{2k-2} (k-1)^{2kn+2} \right\}^{\frac{1}{2nk-
k(k+1)/2+\Delta'+2}} \right].
\endsplit
$$
\endproclaim

\demo{Proof}  If $P\le U^k$ then as in the proof of Lemma 6.5, we have the 
trivial estimate
$$
J_{(n+1)k,k}(P) \le C U^{k(\Delta-\Delta')} P^{2(n+1)k-\frac12 k(k+1)
+\Delta'}.
$$
Next suppose $P > U^k \ge k^k$.  We prove (6.12) 
by induction on $\fl{P}$, observing that (6.12) for integral
$P=m$ implies (6.12) for $m\le P<m+1$.  Assume (6.12) is true for
$\fl{P} \le Q-1$, where $Q$ is an integer $\ge U^k$, and apply
Lemma 6.6.
If the first term in the maximum in the conclusion on
Lemma 6.6 is largest, (6.12) follows from the bound $p\ge U$ and the
induction hypothesis on $m$.  If the second term in the maximum is largest,
(6.12) follows from the upper bound $p\le \eta P^{1/k}$ and the
upper bound on $J_{rk,k}(P)$.
\qed\enddemo


\midinsert
\vbox{
\offinterlineskip
\hrule
\def\twoptsspace{height1.6pt&\omit&\omit&\omit&\omit&\omit&\omit&&\omit&\omit&
\omit&\omit&\omit&\cr}
\halign to \hsize{\tabskip=0pt
\vrule#&
\hbox{\vrule height6.8pt depth2.8pt width0pt}
#\tabskip=1em plus2em minus0.7em&
\hfil#&\hfil#&\hfil#&\hfil#&\hfil#&\vrule#&
\hfil#&\hfil#&\hfil#&\hfil#&\hfil#&\vrule#
\tabskip=0pt\cr
\twoptsspace
&& $\lam$\hfil & $k$\hfil & $n_0$\hfil & $n$ \hfil & $C$\hfil
 && $\lam$ \hfil & $k$ \hfil & $n_0$ \hfil & $n$ \hfil & $C$ \hfil &\cr
\twoptsspace
\noalign{\hrule}
\twoptsspace
&& 2.6--4  &  4 &   1 &  13 &  2.5543 && 45--46 & 46 &  44 & 365 &  3.5897 &\cr
&&   4--5  &  5 &   1 &  17 &  1.7474 && 46--47 & 47 &  46 & 375 &  3.6728 &\cr
&&   5--6  &  6 &   1 &  22 &  1.7805 && 47--48 & 48 &  48 & 386 &  3.7580 &\cr
&&   6--7  &  7 &   1 &  28 &  1.8406 && 48--49 & 49 &  50 & 397 &  3.8453 &\cr
&&   7--8  &  8 &   1 &  34 &  1.9173 && 49--50 & 50 &  52 & 408 &  3.9348 &\cr
&&   8--9  &  9 &   3 &  40 &  1.6808 && 50--51 & 51 &  54 & 419 &  4.0266 &\cr
&&   9--10 & 10 &   3 &  46 &  1.7062 && 51--52 & 52 &  56 & 430 &  4.1207 &\cr
&&  10--11 & 11 &   3 &  52 &  1.7362 && 52--53 & 53 &  58 & 441 &  4.2171 &\cr
&&  11--12 & 12 &   4 &  59 &  1.7678 && 53--54 & 54 &  60 & 452 &  4.3160 &\cr
&&  12--13 & 13 &   4 &  66 &  1.8021 && 54--55 & 55 &  63 & 465 &  4.4174 &\cr
&&  13--14 & 14 &   5 &  73 &  1.8295 && 55--56 & 56 &  65 & 476 &  4.5214 &\cr
&&  14--15 & 15 &   6 &  81 &  1.8669 && 56--57 & 57 &  67 & 487 &  4.6280 &\cr
&&  15--16 & 16 &   6 &  88 &  1.9057 && 57--58 & 58 &  69 & 498 &  4.7373 &\cr
&&  16--17 & 17 &   7 &  96 &  1.9464 && 58--59 & 59 &  71 & 509 &  4.8494 &\cr
&&  17--18 & 18 &   8 & 104 &  1.9883 && 59--60 & 60 &  74 & 522 &  4.9643 &\cr
&&  18--19 & 19 &   8 & 111 &  2.0317 && 60--61 & 61 &  76 & 533 &  5.0821 &\cr
&&  19--20 & 20 &   9 & 119 &  2.0766 && 61--62 & 62 &  79 & 546 &  5.2030 &\cr
&&  20--21 & 21 &  10 & 127 &  2.1229 && 62--63 & 63 &  81 & 557 &  5.3268 &\cr
&&  21--22 & 22 &  11 & 136 &  2.1706 && 63--64 & 64 &  84 & 569 &  5.4539 &\cr
&&  22--23 & 23 &  11 & 143 &  2.2190 && 64--65 & 65 &  86 & 581 &  5.5841 &\cr
&&  23--24 & 24 &  12 & 152 &  2.2688 && 65--66 & 66 &  89 & 593 &  5.7176 &\cr
&&  24--25 & 25 &  13 & 161 &  2.3201 && 66--67 & 67 &  91 & 605 &  5.8546 &\cr
&&  25--26 & 26 &  14 & 169 &  2.3728 && 67--68 & 68 &  94 & 617 &  5.9950 &\cr
&&  26--27 & 27 &  15 & 178 &  2.4270 && 68--69 & 69 &  96 & 629 &  6.1390 &\cr
&&  27--28 & 28 &  17 & 188 &  2.4826 && 69--70 & 70 &  99 & 642 &  6.2867 &\cr
&&  28--29 & 29 &  17 & 196 &  2.5398 && 70--71 & 71 & 102 & 654 &  6.4381 &\cr
&&  29--30 & 30 &  19 & 206 &  2.5987 && 71--72 & 72 & 104 & 666 &  6.5934 &\cr
&&  30--31 & 31 &  20 & 215 &  2.6590 && 72--73 & 73 & 107 & 679 &  6.7527 &\cr
&&  31--32 & 32 &  21 & 224 &  2.7210 && 73--74 & 74 & 110 & 691 &  6.9160 &\cr
&&  32--33 & 33 &  23 & 233 &  2.6797 && 74--75 & 75 & 113 & 704 &  7.0836 &\cr
&&  33--34 & 34 &  25 & 243 &  2.7396 && 75--76 & 76 & 116 & 717 &  7.2553 &\cr
&&  34--35 & 35 &  26 & 252 &  2.8010 && 76--77 & 77 & 118 & 729 &  7.4315 &\cr
&&  35--36 & 36 &  28 & 263 &  2.8641 && 77--78 & 78 & 121 & 742 &  7.6122 &\cr
&&  36--37 & 37 &  29 & 272 &  2.9287 && 78--79 & 79 & 124 & 754 &  7.7975 &\cr
&&  37--38 & 38 &  31 & 283 &  2.9950 && 79--80 & 80 & 127 & 767 &  7.9876 &\cr
&&  38--39 & 39 &  32 & 292 &  3.0630 && 80--81 & 81 & 130 & 780 &  8.1825 &\cr
&&  39--40 & 40 &  34 & 303 &  3.1327 && 81--82 & 82 & 133 & 793 &  8.3825 &\cr
&&  40--41 & 41 &  36 & 313 &  3.2042 && 82--83 & 83 & 136 & 806 &  8.5876 &\cr
&&  41--42 & 42 &  37 & 323 &  3.2775 && 83--84 & 84 & 139 & 819 &  8.7979 &\cr
&&  42--43 & 43 &  39 & 333 &  3.3526 && 84--85 & 85 & 143 & 833 &  9.0136 &\cr
&&  43--44 & 44 &  41 & 344 &  3.4297 && 85--86 & 86 & 146 & 846 &  9.2350 &\cr
&&  44--45 & 45 &  43 & 355 &  3.5088 && 86--87 & 87 & 149 & 859 &  9.4620 &\cr
\twoptsspace
\noalign{\hrule}
}}
\botcaption{Table 6.1} 
\endcaption
\endinsert


\proclaim{Lemma 6.8}  Theorem 2 holds for $2.6\le \lam \le 87$.  In
particular, for each row of Table 6.1, when $\lam$ is in the stated
range,
$$
S(N,t) \le C N^{1-1/(133.66\lam^2)} \quad (N\ge 1).
$$
\endproclaim

\demo{Proof}  Take $k$, $n$ and $n_0$ from a row of the table. 
For reasons connected with the size of $U$ in Lemma 6.7, it is
advantageous to use a completely trivial bound
$$
J_{nk,k}(P) \le P^{2nk-2k} J_{k,k}(P) \le k! P^{2nk-\frac12 k(k+1) +
\Delta_n}, \quad \Delta_n=\frac12 k^2 (1-1/k)
$$
for $1\le n\le n_0$.  We then proceed iteratively, taking a
bound of the form 
$$
J_{nk,k}(P) \le C_n P^{2nk-\frac12 k(k+1) + \Delta_n} \qquad (P\ge 1) 
$$
and producing a bound
$$
J_{(n+1)k,k}(P) \le C_{n+1} P^{2nk-\frac12 k(k+1) + \Delta_{n+1}}
\qquad (P\ge 1),
$$
where $\Delta_{n+1} = (1-1/k)\Delta_n$ and $C_{n+1}$ is the smaller of
the constants coming from Lemmas 6.7 (only for $k\ge 9$)
or 6.5 (with optimal choice of $\om$).
As for the number $\eta$ in Lemma 6.7, (2.1) implies that
$$
\split
\pi(1.12x)-\pi(x) &\ge \frac{x}{\log x} \left[ 1.12 \( 1 + \frac{1/2-\log 1.12}
{\log x} \) - 1 - \frac{3}{2\log x} \right] \\
&\ge \frac{x}{\log x} \( 0.12 - \frac{1.067}{\log x} \)  > 0 \quad (x\ge 7300).
\endsplit
$$
Using a table of primes $<7300$, we find that the following are admissible
choices for $\eta$:
$$
\eta = \cases 17/13 & 9\le k\le 13 \\ 29/23 & 14\le k\le 32 \\ 53/47 &
k\ge 33. \endcases
$$
The optimal value of $\omega$ in Lemma 6.5 is found by solving 
$$
4k^3 k! (1+\om)^{k^2-\Delta_n} = \max ( e^{1.5+1.5/\om}, \tfrac{18}{\om}
k^3\log k ), \tag{6.13}
$$
obtaining a positive real solution $\om_0$.  The solution is unique since
the left side of (6.13) is increasing in $\om$, while the right side is
decreasing.  If $\om_0 \ge 1$, we take $\om=1$.  If $\frac12 \le \om_0 < 1$
we take $\om$ to be either $\frac12$ or $1$, whichever gives the best
constant $C'$.  Otherwise take $\om=\om_0$.

Having computed admissible sequences $C_n$ and $\Delta_n$,
we turn to Lemma 6.3 and Corollary 6.4 to bound $S(N,t)$.
When $2.6\le \lam \le 4$
($k=4$, $n=12$), take $\mu=1-\frac{\lam}{5}$, $M=N^\mu$ and apply
Lemma 6.3.  We have $W\le 2^k M$ and thus
$$
S(N,t) \le 4 \( 4! (8\pi)^4 C_n \)^{\frac{1}{2nk}} N^{1-c} + 2N^{0.8},
\quad c=\frac{1}{2nk} (1 - 0.48(1+\Delta))<0.8.
$$
Hence
$$
S(N,t) \le \( 4 \( 4! (8\pi)^4 C_n \)^{\frac{1}{2nk}}+2 \) N^{1-c} \quad
(N\ge 1).
$$
Applying (6.6) then gives the claimed inequality.  For $\lam\ge 4$, we use
Corollary 6.4 directly, obtaining
$$
S(N,t) \le  \( 4 \( k! (2\pi k)^k C_n \)^{\frac{1}{2nk}}+2 \) N^{1-c}.
$$
Then (6.6) implies the stated claim.  A short computer program (Program 3 in
the appendix) provided the computations of $C_n$ and $\Delta_n$,
and found the best choices of parameters $n_0$ and $n$.
The values of $C$ listed in the table
have been rounded up in the last displayed decimal place.
\qed\enddemo


\vfil\eject

%
%
%
\head 7. Bounding $\zeta(s)$ and $\zeta(s,u)$ \endhead
%
%
%
%
%
%

We start with a crude bound for $\zeta(s)$ and $\zeta(s,u)$ which takes 
care of $s$ with either $\sigma$ or $t$ small.

\proclaim{Lemma 7.1}  Suppose $\frac12 \le \sigma\le 1$, $0<u\le 1$, $t\ge 3$
and $s=\sigma+it$.  If either
$\sigma \le \frac{15}{16}$ or $t\le 10^{100}$, then
$$
|\zeta(s)|, |\zeta(s,u)-u^{-s}| \le 58.1 t^{4(1-\sigma)^{3/2}} \log^{2/3} t.
$$ 
\endproclaim

\demo{Proof}  Applying integration by parts, when $\sigma>0$ we have
$$
\zeta(s,u)  = \sum_{0 \le n\le N} \frac{1}{(n+u)^{s}}
 + \frac{(N+\tfrac12 + u)^{1-s}}{s-1}
+s\int_{N+1/2}^\infty \frac{1/2-\{ w\}}{(w+u)^{s+1}}\, dw, \tag{7.1}
$$
where $N$ is a positive integer.  We take $N=\lfloor t \rfloor$, and note
that $\frac{d^2}{dn^2} (n+u)^{-\sigma} > 0$.  Therefore,
$$
\split
|\zeta(s,u) &- u^{-s} | \le \int_{1/2+u}
   ^{N+1/2+u} \frac{dw}{w^\sigma} + \frac{(N+\tfrac12+u)^{1-\sigma}}{t} +
   \frac{|s|}{2} \int_{N+1/2}^\infty \frac{dw}{(w+u)^{1+\sigma}}\\
&=\int_{1/2+u}^{N+1/2+u} \frac{dw}{w^\sigma}+ 
  \frac{(N+\tfrac12+u)^{1-\sigma}}{t} +
  \frac{|s|(N+\tfrac12+u)^{-\sigma}}{2\sigma} \\
&\le \int_{1/2+u}^{N+1/2+u} \frac{dw}{w^\sigma} +(1+\tfrac{1}{t}) (t+3/2)
^{1-\sigma}.
\endsplit
$$
If $\sigma<1$, 
$$
\int_{1/2+u}^{N+1/2+u} \frac{dw}{w^\sigma} \le \frac{(N+1/2+u)^{1-\sigma}}
{1-\sigma} \le \frac{(t+3/2)^{1-\sigma}}{1-\sigma}
$$
and for all $\sigma\in (0,1]$, we have
$$
\int_{1/2+u}^{N+1/2+u} \frac{dw}{w^\sigma} \le (N+1/2+u)^{1-\sigma}
\int_{1/2+u}^{N+1/2+u} \frac{dw}{w} \le (t+3/2)^{1-\sigma}\log(2N+1).
$$
Therefore, we obtain the inequality
$$
|\zeta(s,u) - u^{-s}| \le (t+3/2)^{1-\sigma} \( 1+\tfrac{1}{t} +
\min \( \tfrac{1}{1-\sigma}, \log(2t+1) \) \). \tag{7.2}
$$
Consider first the case when $t\ge 3$ and $\frac12 \le \sigma \le 
\frac{15}{16}$.  Here $(1-\sigma) \le 4(1-\sigma)^{3/2}$, so by (7.2)
$$
|\zeta(s,u) - u^{-s}| \le \sqrt{1.5} t^{4(1-\sigma)^{3/2}} \( \tfrac43 + 16
\) \le 21.3 t^{4(1-\sigma)^{3/2}}.
$$
Next, if $\frac{15}{16} \le \sigma \le 1$ and $3\le t\le 10^{100}$, (7.2)
gives
$$
|\zeta(s,u) - u^{-s}| \le (t+3/2)^{1-\sigma} (1+1/t+\log(2t+1)).
$$
If $3\le t\le 10^6$, the right side is $\le 36.8$.  If $t>10^6$, the
right side is
$$
\le 1.123 t^{1-\sigma} \log t = 1.123 \( t^{4(1-\sigma)^{3/2}} \log^{2/3} t
\) \( t^{1-\sigma-4(1-\sigma)^{3/2}} \log^{1/3} t \).
$$
The maximum of $1-\sigma-4(1-\sigma)^{3/2}$ is $\frac{1}{108}$, thus
$$
|\zeta(s,u) - u^{-s}| \le 58.1 t^{4(1-\sigma)^{3/2}} \log^{2/3} t.
$$
Lastly, taking $u\to 0^+$ shows that the lemma holds for $|\zeta(s)|$
as well.
\qed\enddemo

\proclaim{Lemma 7.2}  If $s=\sigma+it$, $\frac{15}{16} \le \sigma \le 1$,
$t \ge 10^{100}$ and $0<u\le 1$, then
$$
\left| \zeta(s,u) - \sum_{0\le n\le t} (n+u)^{-s} \right| \le 10^{-80}.
$$
\endproclaim

\demo{Proof}  Let $E(s,u) = \zeta(s,u) - \sum_{0\le n\le t} (n+u)^{-s}$.  By
(7.1) with $N=\lfloor t \rfloor$,
$$
\split
|E(s,u)| &\le \frac{(t+3/2)^{1-\sigma}}{t} + |s| \biggl|
\int_{\fl{t}+1/2+u}^{\infty} \frac{1/2-\{w\}}{w^{1+s}}\, dw
\biggr| \\
&\le \frac{(t+3/2)^{1-\sigma}}{t} + \frac{3|s|}{4(t-1/2)^{\sigma+1}}
 + |s| \left| \int_t^{t^2} \frac{1/2 - \{ w \}}{w^{s+1}} \, dw \right|
 + \frac{|s| t^{-2\sigma}}{2\sigma} \\
&\le 10^{-81} + (t+1) \left| \int_t^{t^2} \frac{\{ w \} - 1/2}{w^{\sigma+1}}
(\cos(t\log w)-i\sin(t\log w)) \, dw 
\right|.
\endsplit
$$
We bound the intergal using the Fourier expansion $\{ x \} - \frac12 =
-\frac{1}{\pi} \sum_{m=1}^\infty \frac{\sin(2\pi mx)}{m}$, as in
[3].  We also use the trigonometric identities
$$
\sin a \sin b = \frac{\cos(a-b)-\cos(a+b)}{2}, \qquad
\sin a \cos b = \frac{\sin(a+b)+\sin(a-b)}{2}.
$$
Therefore, writing
$$
I_m = \max_{h=\sin, \cos} 
\left| \int_t^{t^2} \frac{h(t\log x + 2\pi mx)}{x^{1+\sigma}}\, dx \right| + 
\left| \int_t^{t^2} \frac{h(t\log x - 2\pi mx)}{x^{1+\sigma}}\, dx \right|
$$
and separating real and imaginary parts, we obtain
$$
|E(s,u)| \le 10^{-81} + \frac{t+1}{\pi} \sum_{m=1}^\infty \frac{I_m}{m}.
\tag{7.3}
$$
To bound $I_m$, let $f(x) = x^{-\sigma}/(t\pm 2\pi mx))$ and
$g(x) = k(t\log x \pm 2\pi mx)$, where $k'(x)=h(x)$ and
$k(x) \in \{ \pm \sin(x),\pm \cos(x)\}$).   Since $f$ is 
monotonic on $[t,t^2]$, we obtain
$$
\split
\left| \int_t^{t^2} \frac{h(t\log x \pm 2\pi mx)}{x^{1+\sigma}}\, dx \right|
&=\left| \int_t^{t^2} f(x) g'(x)\, dx \right| \\
&= \left| f(t^2)g(t^2)-f(t)g(t) -
\int_t^{t^2} g(x)f'(x)\, dx \right| \\
&\le |f(t)g(t)|+|f(t^2)g(t^2)| + \max_{t\le x\le t^2}
|g(x)| \int_t^{t^2} |f'(x)|\, dx \\
&= |f(t)g(t)|+|f(t^2)g(t^2)| +\max_{t\le x\le t^2}
|g(x)| |f(t^2)-f(t)|  \\
&\le \frac{4}{t^{1+\sigma}(2\pi m \pm 1)}.
\endsplit
$$
Therefore,
$$
I_m \le \frac{4}{t^{1+\sigma}(2\pi m +1)} +
\frac{4}{t^{1+\sigma}(2\pi m -1)} = \frac{16\pi m}{t^{1+\sigma} (4\pi^2
m^2-1)} \le \frac{16\pi}{(4\pi^2-1) t^{1+\sigma} m}.
$$
Together with (7.3), this proves the lemma.
\qed\enddemo

\proclaim{Lemma 7.3}  Suppose that $S(N,t) \le C N^{1-1/(D\lambda^2)}$
$(1\le N\le t)$ for positive constants $C$ and $D$, where
$\lambda=\frac{\log t}{\log N}$.
Let $B=\frac29 \sqrt{3D}$.  Then, for $\frac{15}{16} \le \sigma \le 1$,
$t\ge 10^{100}$ and $0<u\le 1$, we have
$$
\split
|\zeta(s)| &\le
\( \frac{C+1+10^{-80}}{\log^{2/3} t}+1.569CD^{1/3} \)  t^{B (1-\sigma)^{3/2} }
 \log^{2/3} t, \\
|\zeta(s,u)-u^{-s}| &\le
\(  \frac{C+1+10^{-80}}{\log^{2/3} t}+1.569CD^{1/3} \)  t^{B (1-\sigma)^{3/2} }
 \log^{2/3} t. \\
\endsplit
$$
\endproclaim

\demo{Proof} Let 
$$
S_1(u) = \sum_{1\le n \le t} (n+u)^{-s}.
$$
By Lemma 7.2, $|\zeta(s,u)-u^{-s}| \le 10^{-80} + S_1(u)$.
Put $r = \lceil\frac{\log t}{\log 2} \rceil$.  By
partial summation,
$$
\split
|S_1(u)| &\le 1+ \sum_{j=0}^{r-1} \left| \sum_{2^j < n \le \min(t,2^{j+1})}
(n+u)^{-\sigma-it} \right| \\
&\le 1+\sum_{j=0}^{r-1} (2^j)^{-\sigma}  S(2^j,t) \\
&\le 1+ C \sum_{j=0}^{r-1} e^{g(j)},
\endsplit
$$
where
$$
g(j) = (1-\sigma)(j\log 2)-\frac{(j\log 2)^3}{D\log^2 t}.
$$
As a function of $x$, $g(x)$ is increasing on $[0,x_0]$ and decreasing
on $[x_0,\infty)$, where $x_0\log 2 = \sqrt{D(1-\sigma)/3} \log t$.  Thus
$$
\split
\frac{|S_1(u)|-1}{C} &\le e^{g(x_0)} + \int_0^{r} e^{g(x)}\, dx \\
&\le t^{B(1-\sigma)^{3/2}} + \frac{D^{1/3} \log^{2/3} t}{\log 2} 
\int_0^\infty e^{3y^2 u-u^3}\, du,
\endsplit
$$
where  $y=\sqrt{(1-\sigma)/3} D^{1/6}\log^{1/3} t$.
To bound the last integral, we make use of the inequality
$$
e^{-2y^3} \int_0^\infty e^{3y^2u-u^3}\, du \le 1.0875034 \qquad
(y \ge 0),
$$
where the maximum occurs near $y=0.710$.  Therefore
$$
\frac{|S_1(u)|-1}{C} \le t^{B(1-\sigma)^{3/2}} \( 1 +
 \frac{1.0875034}{\log 2} D^{1/3}\log^{2/3} t \),
$$
which proves the lemma.
\qed\enddemo

\demo{Proof of Theorem 1}
Apply Lemma 7.3 using $C=9.463$, $D=133.66$ (from Theorem 2).
\qed\enddemo

\vfil\eject
%
%
%
%
%

\head 8. Possible improvements to the constant $B$ \endhead

There are a number of ways in which the constant $B$ in Theorem 1 may be
improved, and we sketch three of them below.  To provide complete details
would involve a substantial lengthening of this paper, and even more
work would be required to obtain a decent constant $A$.
Taken together, the three ideas have the potential 
to reduce the constant $B$ only to about $4.1$.

\medskip
1. As noted in section 3, there are some improvements possible in the method
for bounding $J_{s,k}(P)$.  Tyrina's method could be used for small
$s$ (when $\Delta \ge \frac49 k^2$), and in Lemma 3.5 we could take
$r\approx \sqrt{k^2+k-2\Delta}$ in Lemma 3.5.  The end result is a 
slight reduction in the constant $\frac38$ appearing in the definition
of $\Delta_s$ in Theorem 3.  This can lower $B$ by less than $0.02$.
\medskip

2.  As mentioned in section 4, the use of repeat efficient differencing
(repeatedly forming divided differences of the polynomials $\Psi_j$
as in [34]) produces superior bounds for $J_{s,g,h}(\CPR)$.  
Preliminary computations indicate a potential reduction in $B$ of $4-5\%$,
or $0.2$ at most, making it hardly worth the effort of working out
the details.  There is also the problem of obtaining good explicit
constants (e.g. $e^C$ in Theorem 4).  In particular, when Wooley's methods
are used directly, the constants $C$ are far too large to be of any use
in bounding $\zeta(s)$.  Referring to Lemma 4.1 of [34], relations
(4.9) and (4.10) essentially bound $J_{s,g,h}$ in terms
of $J_{s-1,g,h}$.  When iterated, the constants grow too rapidly with
$s$.  In our Lemma 4.1 above, we avoided this pitfall by an application
of H\"older's inequality at the end of the third case (assuming $S_3=
\max(S_1,S_2,S_3,S_4)$), a tool which is unavailable when using repeat
efficient differencing.  Incidentally, this idea was also used in the
proof of Lemma 6.7 above.  Presumably some clever argument would overcome
this problem.
\medskip

3.  In the estimation of the quantity $T$ in section 5, the number of
solutions of (5.5) may be bounded in a more sophisticated way.  
First we note that when $sM_2^j |\g_j| \le \frac14$ (essentially $j \ge
\frac{\lam}{1-\mu_2}$), $\curly D_j$ is the set of integers in an
interval of the form $[-D_j,D_j]$, where $D_j$ is a non-negative integer.
If in addition $|\g_j| \ge \frac{1}{2rM^j}$ (essentially $\frac{\lam}{1-\mu_2}
\le j\le \frac{\lam}{1-\mu_1}$), in fact $\curly D_j = \{ 0 \}$ (i.e.
$D_j=0$ in this case).

\def\bg{{\tilde{g}}}
\def\DDD{\boldkey D}
Let $h_0$ be the smallest integer with $D_{h_0}=0$ and
let $\bg$ be the largest integer with $|\g_\bg| \ge \frac{1}{2rM^\bg}$.
Assuming $h_0 \le h \le \bg \le g \le k$,
the number of solutions of (5.5) is at most $J^*_{s,g,h}(\BB;\DDD)$,
the number of solutions of
$$
\sum_{i=1}^s (x_i^j-y_i^j) =d_j \quad (h\le j\le g), \tag{8.1}
$$
with $x_i,y_i \in \BB$ and $|d_j| \le D_j$ for each $j$.
Now set $\BB=\CPR$ and for non-negative integers $D$ define
$$
H(\a;D) = \frac{1}{D+1} \biggl| \sum_{|x| \le D} e(\a x) \biggr|^2 =
\sum_{|x| \le 2D} \pfrac{2D+1-|x|}{D+1} e(\a x).
$$
Define $f(\baa)$ as in section 5 and let
\def\Jt{{\widetilde{J}}}
$$
\Jt_{s,g,h}(\BB;\DDD) = \iu{g-h+1} |f(\baa)|^{2s} G(\baa)\, d\baa, \qquad
G(\baa)=H(\a_h;D_h) \cdots H(\a_g;D_g).
$$
Then $J^*_{s,g,h}(\BB;\DDD) \le \Jt_{s,g,h}(\BB;\DDD)$, because the latter
quantity counts the solutions of (8.1) each with weight
$$
w(\dd) = \prod_{j=h}^g \max \( 0, \frac{2D_j+1-|d_j|}{D_j+1} \). \tag{8.2}
$$
Since $G(\baa)$ is real and non-negative, we may follow the proof of 
Lemma 4.1 to bound $\Jt_{s,g,h}(\BB;\DDD)$.  We show the proof in
some detail, as this method may have other applications.

\proclaim{Lemma 8.1}  Suppose $h,\bg,g,r,s$ are positive integers with
$$
g\ge \bg \ge h \ge 9, \quad t=\bg-h+1, \quad h\le r\le \bg, \quad s\ge 2t.
$$
Further suppose that
$$
0\le D_j \le sP^j \quad (h\le j\le g), \qquad D_j=0 \quad (h\le j\le \bg)
$$
and
$$
R=P^\eta > g^2, \qquad |\CPR| \ge P^{1/2}, \qquad P>(8s 2^{g/s})^8.
$$
Then
$$
\multline
J^*_{s,g,h}(\CPR;\DDD) \le \max \biggl[ (8s)^{2s} (22t^2)^{2s/\eta} 2^g P^{s(1+
1/r)}, \\ 
4 g^{2t(1+1/(r\eta))} (P^{1/r}R)^{2s-2t+\frac12(r-h)(r-h+1)} 2^g P^t
 J^*_{s-t,g,h}(\CC(P^{1-1/r},R);\boldkey E) \biggr],
\endmultline
$$
where $E_j = \fl{\frac{2D_j}{P^{j/r}}}$ for $h\le j\le g$.
\endproclaim

\def\iuz{\iu{g-h+1}}
\demo{Sketch of proof}  First, $J^*_{s,g,h}(\BB;\DDD) \le \Jt_{s,g,h}
(\BB;\DDD)$,
and we follow the proof of Lemma 4.1 to bound  $S_0:=\Jt_{s,g,h}(\BB;\DDD)$.
Define $S_1,\ldots,S_4$ analogously,
and consider the same four cases.  When $S_1$ is the largest, we obtain
$$
S_0 \le (8s)^{2s} \iuz |f(\baa;P^{1/r})|^{2s} G(\baa)\, d\baa
\le (8sP^{1/r})^{2s}  \iuz G(\baa)\, d\baa.
$$
By (8.2), the last integral is $\le 2^{g-h+1} \le 2^g$, so
$S_0 \le  2^g (8sP^{1/r})^{2s}$.  However, the hypotheses imply
$S_0 \ge (P-1)^{s/2}$, giving a contradiction.
When $S_2$ is the largest,
$$
\split
S_0 &\le 4t^2 \iuz |f(\baa)^{2s-2} f(2\baa)| G(\baa)\, d\baa \\
&\le 4t^2 S_0^{1-1/s} \( \iuz |f(2\baa)|^{2s} G(\baa)\, d\baa \)^{\frac{1}{2s}}
  \( \iuz G(\baa)\, d\baa \)^{\frac{1}{2s}}.
\endsplit
$$
By considering the underlying Diophantine equations, the first integral on
the right is $\le S_0$, thus $S_0 \le 4t^2 S_0^{1-\frac{1}{2s}} 2^{g/(2s)}$,
whence $S_0 \le (4t^2)^{2s} 2^g$.  Again by the lower bound $S_0 \ge 
(P-1)^{s/2}$
and the assumed lower bound on $P$, this gives a contradiction.
Therefore, $S_0=4\max(S_3,S_4)$.

When $S_3$ is largest, we obtain
$$
\split
S_0 &\le (8s)^2 (8et^2)^{2/\eta} P^{1+1/r} \iuz |f(\baa)|^{2s-2} G(\baa)\,
  d\baa \\
&\le (8s)^2 (8et^2)^{2/\eta} P^{1+1/r} \( \iuz |f(\baa)|^{2s} G(\baa)\,
  d\baa \)^{1-\frac{1}{s}} \( \iuz G(\baa)\, d\baa \)^{\frac{1}{s}} \\
&\le  (8s)^2 (8et^2)^{2/\eta} P^{1+1/r} S_0^{1-1/s} 2^{g/s}.
\endsplit
$$
Therefore $S_0 \le (8s)^{2s} (22t^2)^{2s/\eta} 2^g P^{s(1+1/r)}$.

If $S_4$ is the largest, we add a factor $G(\baa)$ to each $X_i(\baa)$
and $Y_i(\baa)$ and obtain
$$
S_0 \le 4(P^{1/r} R)^{2s-2t} \max_{P^{\frac1{r}} < q \le P^{\frac{1}{r}} R}
W(q),
$$
where $W(q)$ counts solutions of
$$
\sum_{i=1}^t (x_i^j-y_i^j) + q^j \sum_{i=1}^{s-t} (u_i^j-v_i^j)=d_j
\qquad (h\le j\le g)
$$
each with weight $w(\dd)$.  Since $d_j=0$ for $h\le j\le \bg$, the argument
in the proof of Lemma 4.1 implies that there are at most
$g^{2t(1+1/(r\eta))} q^{(r-h)(r-h+1)/2} P^t$ possibilities for $\xx,\yy$
(note that here $t=\bg -h+1$). 
Let $\curly S$ be the set of possible $\xx,\yy$ and put
$$
F(\baa)=\sum_{(\xx,\yy)\in \curly S} e \( \sum_{j=h}^g \a_j (x_1^j-y_1^j +
\cdots + x_t^j-y_t^j) \).
$$
Putting $\widetilde{\baa}=(q^h \a_h, \ldots, q^g\a_g)$, we obtain
$$
\split
W(q) &\le \iuz |F(\baa)| |f(\widetilde{\baa};P/q)|^{2s-2t} G(\baa)\,
  d\baa \\
&\le g^{2t(1+1/(r\eta))} q^{(r-h)(r-h+1)/2} P^t
  \iuz |f(\widetilde{\baa};P/q)|^{2s-2t} G(\baa)\, d\baa.
\endsplit
$$
The integral on the right counts the solutions of
$$
q^j \sum_{i=1}^{s-t} (u_i^j-v_i^j)=d_j \qquad (h\le j\le g),
$$
each counted with weight $w(\dd)$.  Since $q>P^{1/r}$, this is at most
$2^g$ times the number of solutions of 
$$
 \sum_{i=1}^{s-t} (u_i^j-v_i^j)=e_j \qquad (h\le j\le g),
$$
with $u_i,v_i \in \CC(P^{1-1/r},R)$ and $|e_j| \le 2D_j/P^{j/r}$.
This proves the lemma in the last case.
\qed\enddemo

In Lemma 8.1  it is common that there are
more zeros among the numbers $E_i$ than among the numbers $D_i$.  Thus, as
Lemma 8.1 is iterated, $t$ steadily increases (if $t$ reaches $g-h+1$, then
one can apply the bounds from \S 4).  This is the primary source
of the improvement over Lemma 4.1, but the analysis of
the exponents of $P$ and the constants is much more complicated.
The analysis becomes even more
complex if repeat efficient differencing is used.  
By taking optimal parameters, using Lemma 8.1 in place of Lemma 4.1
has the potential to reduce $B$ by about $0.09$, or $\approx 2\%$.
\bigskip

Lastly, we indicate what is the limit of our method, i.e. the limit of
what could be accomplished with Lemma 5.1.  Assume now that
the lower bound (1.5) for $J_{s,k}(P)$ is close to the truth, i.e.
$J_{s,k}(P) \le C(k,s) P^{s}$ for $s\le \frac12 k(k+1)$.  Assume
also best possible upper bounds $J_{s,g,h}(\BB) \le C(s,g,h) P^{s}$
for $s\le \frac{t}{2}(g+h)$, valid for any $\BB \subset [1,P]$.
Adopt the notations from section 5.
With these assumptions, it turns out that the best choices for
$r,s,\mu_1,\mu_2$ are given by 
$$
r=\frac{k(k+1)}{2}, \quad s=\frac{t(g+h)}{2}, \quad \mu_1=\mu_2=\mu=\frac16.
$$
Also, one takes $\phi$ very close to (and larger than) $\frac{1}{1-\mu}$
and $\gamma$ very close to (and smaller than) $\frac{1}{1-\mu}$.
Plugging these values into (5.22) yields
$$
\lam^2 E = \frac{2}{27} - \eps,
$$
where $\eps \to 0^+$ as $\phi-\gamma \to 0$.  An application of
Lemma 7.3 (with $D=27/2+\eps'$) gives Theorem 1 with a constant
$B=\sqrt{2}+\eps''$ (valid for $\sigma\ge \frac{15}{16}$),
where $\eps',\eps''$ can be taken arbitrarily
small.

%
%
%
%

\bigskip

\head Appendix: Computer program listings \endhead

{
\input manmac
\eightpoint
\baselineskip 9pt
\beginlines
|/* PROGRAM 1. exponents and constants in Vinogradov's integral for small k.|
|   Used to prove the second part of Theorem 3; written 12/12/2000 K. Ford  */|
|#include <stdio.h>|
|#define max(x,y) (((x)>(y))?(x):(y))|
|#define min(x,y) (((x)>(y))?(y):(x))|
||
|double newdel(k,r,del)       /* returns delta_0(k,r,del) */|
|     double k,r,del;|
|{|
|  double y, sqrt(),p,tkr;|
|  long j,jj;|
||
|  if ((r<4.0) || (r>k)) return(2.0*del);  /* invalid r */|
|  tkr = 2.0*k*r; y=2.0*del-(k-r)*(k-r+1.0);|
|  if ((y<0.0)||(2.0*k/(tkr+y))<=1.0/(k+1.0)) return(del*2.0); /* invalid r */|
|  j = min((long) (0.5*(3.0+sqrt(4.0*y+1.0))), 9*r/10); |
|  p = 1.0/r;|
|  for (jj=j-1; jj>=1; jj--) |
|    p = 0.5/r+0.5*(1.0+(jj*jj-jj-y)/tkr)*p;|
|  return(del-k+0.5*p*(tkr-y));|
|}|
|main()|
|{|
|  long j,k,k0,k1,r,r0,r1,n,bestr,s;|
|  double kk,logk,del0,del1,sqrt(),log(),exp(),bestdel, goal, maxs, eta, om;|
|  double logH,logW,logC,k3,theta,thetamax;|
|  printf("enter k range : "); scanf("
|  maxs = 0.0;  thetamax=0.0;|
|  for (k=k0; k<= k1; k++) {|
|    kk=(double) k;|
|    logk=log(kk); k3 = kk*kk*kk*logk;|
|    om=0.5; for (j=1;j<=10;j++) om=1.5/(log(18.0*k3/om)-1.5);|
|    eta = 1.0+om;|
|    logW = (kk+1.0)*max(1.5+1.5/om, log(18.0/om*k3));|
|    del0 = 0.5*kk*kk*(1.0-1.0/kk);|
|    goal = 0.001*kk*kk;|
|    logH = 3.0*kk*logk+(kk*kk-4.0*kk)*log(eta); /* log(k^3k eta^(k^2-4k) */|
|    logC = kk*logk;                           /* upper bound for log(k!) */|
|    for (n=1;;n++) {|
|      r0 = (long) (sqrt(kk*kk+kk-2.0*del0)+0.5)-2; r1 = r0+4; /* r range */|
|      bestdel=kk*kk; bestr=-1;|
|      for (r=r0;r<=r1;r++) {|
|        del1=newdel(kk,(double) r,del0);|
|        if (del1<bestdel) { bestdel=del1; bestr=r; }|
|      }|
|      del1=bestdel; r=bestr;|
|      if ((del1 >= del0) || (r<r0)) exit(-1);|
|      logC += max(logH + 4.0*kk*n*log(eta),logW*(del0-del1));|
|      if (del1<=goal) {                             /* reached goal */|
|        s=(long) ((n+(del0-goal)/(del0-del1))*kk+1);|
|        theta = logC/k3;|
|        printf("
|            s/kk/kk,eta,theta);|
|        if ((s/kk/kk) > maxs) maxs=s/kk/kk;|
|        if (theta>thetamax) thetamax=theta;|
|        break;|
|      }|
|      del0=del1;|
|    }|
|  }|
|  printf("\n max s = 
|}|
\endlines

\beginlines
|/* PROGRAM 2.  Find optimal parameters for use in bounding S(N,t) for the|
|   Riemann zeta function  : intermediate lambda.  For Lemma 5.3,  |
|   lambda in [84,220].  By K. Ford 10/22/2001        */|
|#include <stdio.h>|
|#include <math.h>|
|long k,g,h,s,r,t, g0, h0,g1,h1,flag;|
|double mu1,mu2,xi,lam,lam1,lam2,D,sigma, Y, goal;|
|void calc(ex,c,pr)|
|     double *ex,*c; int pr;|
|{|
|  double kk,logk, k2, log(),exp(),pow(),floor(), ceil();|
|  double th,rr,ss,tt,gg,hh,rho,H,E1,E2,E3,m1,m2,Z0,Z1,reta,|
|     logC1,logC2,logC3,logC,dc;|
|  k=(long) (lam/(1.0-mu1-mu2)+0.000003);|
|  /* if (k<129) exit(-1); */|
|  kk=(double) k;|
|  logk=log(kk); k2=kk*kk;|
|  rho=3.21432; th=2.3291; | 
|  if (k<=199)  { rho=3.21734; th=2.3849; }    /* 150 to 199 */|
|  if (k<=149)  { rho=3.22313; th=2.4183; }    /* 129 to 149 */|
|  r = (long) (rho*k2+1.0);|
|  rr=(double) r;  ss=(double) s;|
|  gg=(double) g;  hh=(double) h; tt=(double) t;|
|    /* calculate minimum H = Z1 + lam*Z2 */|
|  m1 = floor(lam/(1.0-mu1));|
|  m2 = floor(lam/(1.0-mu2));|
|  Z0 = 0.5*((m1*m1+m1)*(1.0-mu1)+(m2*m2+m2)*(1.0-mu2)-hh*hh+hh-(1.0-mu1-mu2)*|
|	    (gg*gg+gg));|
|  Z1 = hh+gg-m1-m2-1.0;|
|  if (Z1<0.0) H = Z0 + lam2*Z1;|
|  else H=Z0 + lam1*Z1;   /* H is now the H' from Lemma 5.3 */|
|  reta = xi*pow(gg,1.5);   /* 1/eta */|
|  E1 = 0.001*k2;|
|  E2 = 0.5*tt*(tt-1.0)+hh*tt*exp(-ss/(hh*tt))+ss*ss/(2.0*tt*reta);|
|  E3 = log(Y*lam1*lam1)/(7.5*Y*lam1*lam1*lam1*lam1); |
|  *ex = (-E3 + (1.0/(2.0*rr*ss))*(H-mu1*E1-mu2*E2))*lam1*lam1;|
|  logC1=th*k2*kk*logk;|
|  logC2 = ss*ss/tt+10.5*xi*xi*tt*gg*gg*log(gg)*log(gg)/D;|
|  logC2 -= (ss*log(0.1*reta)*((reta+hh)*pow(1.0-1.0/hh,ss/tt)-h));|
|  logC3=1.04*reta*log(10.82*reta);|
|  logC = logC3/rr+(5.0*lam2*log(lam2)+logC1+logC2)/(2.0*rr*ss);|
|  *c = exp(logC)+1.0/kk;  /* constant for exponent ex */|
|  if (pr==1) {|
|    printf("
|    if (g>0) printf(" 
|		    s,g-g0,h1-h,t,1.0/(*ex)+0.00005,*c+0.00005);|
|    else printf("\n");|
|  }|
|}|
|main()|
|{|
|  double E,lam8,lam9,r[9],tmp,maxex,con,maxcon,bestth,bestcon,bp[5000]; |
|     /* bp[] are endpoints of intervals */|
|  long i,j,i0,w,n,m,maxm,bestg,besth, bests,s0,s1;|
|  mu1 = 0.1905; mu2 = 0.1603;|
|  goal=133.66;|
|  while (1) {  |
|    printf("enter Y  : "); scanf("
|    D = 0.1019*Y;|
|    printf("enter xi : "); scanf("
|    printf("enter sigma : "); scanf("
|    if (sigma<0.0) flag=1; else flag=0;  |
|    /* flag=1 means let the program find the best value of s */|
|    printf("enter lambda range: "); scanf("
|    if ((lam9<lam8) || (lam8<=80.0) || (lam9>=300.0)) continue;|
|    printf("    approx.\n");|
|    printf("  lambda range       k   s    a   b   t     exp      const\n");|
|    printf("----------------   ---- ---  --- --- ---  --------  --------\n");|
|    bp[1] = lam8; bp[2] = lam9; j=3;       /* make list of endpoints */|
|    i0 = (long) (lam9/(1.0-mu1-mu2))+10;|
|    for (i=1; i<=i0;i++) {|
|      w=(double) i;|
|      r[1]=w*(1.0-mu1);|
|      r[2]=w*(1.0-mu2);|
|      r[3]=(w-0.000003)*(1.0-mu1-mu2);|
|      for (m=1;m<=3;m++) if ((r[m]<lam9) && (r[m]>lam8)) bp[j++]=r[m];|
|    }|
|    n=j-1;      /* number of endpoints */|
|    for (i=1; i<=n-1; i++) for(j=i+1;j<=n;j++)   /* Bubble sort */|
|        if (bp[j]<bp[i]) { tmp=bp[i]; bp[i]=bp[j]; bp[j]=tmp; }|
|    maxex=0.0;  /* maximum exponent of N */|
|    maxcon = 0.0;  /* maximum constant */|
|    for (j=1; j<=n-1; j++) {|
|      lam = 0.5*(bp[j]+bp[j+1]);   /* midpoint of interval */|
|      lam1=bp[j]; lam2=bp[j+1];    /* endpoints */|
|      g0 = (long) (lam/(1.0-mu1)+1.0); g1=g0+1; /* g range */|
|      h1 = (long) (lam/(1.0-mu2)); h0=h1-1;     /* h range */|
|      bestg=-1; besth=-1; bestth=1.0e20; bestcon=1.0e40;|
|      for (g=g0;g<=g1;g++) for (h=h0;h<=h1;h++) {|
|       t=g-h+1;|
|       if ((g>=100) && ((double) g <= 1.254*lam1)) {  /* condition (5.16) */|
|         if (flag==0) {|
|           s0=(long) (sigma*h*t+1.0); s1=s0;|
|         }|
|         else {|
|           s0=h*(t-1)/4;|
|           s1=h*t/2;|
|         }|
|         for (s=s0; s<=s1; s++) {|
|           calc(&E,&con,0);|
|           if ((E>0.0) && (1.0/E < goal) && (con<bestcon)) {  |
|             /* look for best constant such that 1/exponent < goal */|
|             bestth=1.0/E; bestg=g; besth=h; bests=s;|
|           }|
|         }|
|       }|
|      }|
|      g=bestg; h=besth; t=g-h+1;|
|      s=bests;|
|      calc(&E,&con,1);|
|      if (1.0/E>maxex) maxex=1.0/E;|
|      if (con>maxcon) maxcon=con;|
|    }|
|    printf("  max. ex: 
|  }|
|}|
\endlines

%
%

\beginlines
|/* PROGRAM 3.  find optimal parameters for use in bounding S(N,t)|
|   for small lambda; Section 6.  Written by K. Ford 10/20/2001 */|
|#include <stdio.h>|
|#include <math.h>|
|#define max(x,y) (((x)>(y))?(x):(y))|
|#define min(x,y) (((x)<(y))?(x):(y))|
|long k,n0;|
|double kk, logk, logk1, pi, eta, logeta, L32, lam, lam4, lkf,k3,logA,B,C;|
|double Delta[10000], logC[10000]; /* Delta and log of constants */|
|double log(), exp(), pow(), sqrt();|
|/* #define DEBUG  */|
|double logV(double w)   /* log(V(w))  */|
|{|
| if (w==1.0) return(k3);|
|  if ((w<=0.5)&&(w>0.0)) return(max(1.5+1.5/w,k3+log(3.0/w)));|
|  exit(-1);|
|}|
|double F(double w)|
|{|
|  return((1.0+w)*exp(logA/B)-exp(logV(w)*C/B));|
|}|
|double bestomega(int n)  /* best omega value for Lemma 6.5 */|
|{|
|  double w0,w1,w2;|
|  B = kk*kk-Delta[n];            /* exponent of (1+w) */|
|  C = Delta[n];                  /* exponent of V     */|
|  if (F(1.0)<=0.0) return(1.0);  /* take w=1          */|
|  if (F(0.5)<=0.0) {             /* take w=1 or 1/2   */|
|    if (exp(logV(0.5)*C/B)<2.0*exp(logA/B)) return(0.5);|
|    else return(1.0);|
|  }                              /* solve F(w)=0 */|
|  w0=0.5; w1=0.2; while (F(w1)>=0.0) w1*=0.5;|
|  while (((w0-w1)/w1)>=0.0000001) {|
|    w2=0.5*(w0+w1);|
|    if (F(w2)>0.0) w0=w2; else w1=w2;|
|  }|
|  return(w1);|
|}|
|void calcparm()   /* calculate Delta_n and C_n */|
|{|
|  long n1,n,i;|
|  double f, s, logU, omega, logM1, logM2, AA, BB;|
|  kk=(double) k;|
|  logk=log(kk);|
|  logk1=log(kk-1.0);|
|  k3=3.0*logk+log(6.0*logk);  /* log(6k^3 log k) */|
|  lkf=0.0; for (i=2;i<=k;i++) lkf += log(((double) i));  /* log(k!) */|
|  logA = 3.0*logk+lkf+log(4.0);  /* log(4k^3 k!)   */|
|  logeta=log(eta);|
|  L32=log(32.0)-lkf;|
|  n1 = (long) (2.6*kk*logk+50);|
|  if (n1>=9999) n1=9998;  /* calculate constants up to n=n1    */|
|  for(i=1;i<=n0;i++) {    /* use trivial bound for 1<= n<= n0  */ |
|    Delta[i] = 0.5*kk*(kk-1.0); |
|    logC[i] = lkf;|
|  }|
|  f = 1.0-1.0/kk;|
|  for (n=n0+1; n<=n1+1; n++) Delta[n]=f*Delta[n-1];|
|  for (n=n0; n<=n1; n++) {|
|    s=kk*n;|
|    omega=bestomega(n);|
|    logM1 = max(logV(omega)*C,logA+B*log(1.0+omega));|
|                     /* M1=multiplier for constant in Lem. 6.5 */|
|    if (k>= 9) {     /* Lemma 6.7 only for k>=9  */|
|      AA =(kk*kk-Delta[n])*logeta+2.0*kk*log(s+kk)+L32;|
|      logU = (2.0*kk-2.0+(2.0*s+2.0)*logk1)/|
|        (2.0*s+2.0-0.5*kk*(kk+1.0)+Delta[n+1]);|
|      if (logU<logk) logU=logk;|
|      BB=Delta[n]*logU;|
|      logM2 = max(AA,BB);   /* M2=multiplier for constant in Lemma 6.7 */|
|    }|
|    else logM2=1.0e40; |
|    logC[n+1] = logC[n] + min(logM1,logM2);|
|#ifdef DEBUG|
|    printf("  logM1=
|#endif|
|  }|
|}|
|int exponent(n,c,pr)      /* from Lem. 6.3, 6.4 */|
|  int n,pr; double *c;    /* return constant in 'c' */|
|{|
|  double s,goal,logd,c1,e,mu,log(),pow(),exp();|
|  lam=kk-1.0; if (k==4) lam=lam4;  /* lower limit of lambda */|
|  mu=1.0-lam/(kk+1.0);             /* largest mu */|
|  s=kk*n;|
|  logd = log(4.0) + 0.5/s*(logC[n]+lkf+kk*log(2.0*kk*pi));|
|  logd = log(exp(logd)+2.0);  /* add 2 */|
|  goal = 133.66*lam*lam;      /* goal for denominator */|
|  e = (1.0-(1.0+Delta[n])*mu)/(2.0*s);|
|  if (e<1.0/goal) return(-1);   /* exponent not good enough */|
|  *c = exp(logd/e/goal); |
|  if (((*c) <= 10000.0) && (pr==1))|
|    printf("n=
|  return(0);|
|}|
|main()|
|{|
|  double log(),bestc,c,e,mu, CC[200];|
|  long bestn,bestn0,n, n2, i,k1,k2,j,nn[200], n00[200], n01, n02;|
|  pi=3.1416;  /* good enough upper bound */|
|  printf(" k range : "); scanf("
|  if (k1<4) exit(0);|
|  if (k1==4) {|
|    printf("enter lower bound on lam for k=4 : ");|
|    scanf("
|  }|
|  /* printf(" n0 range : "); scanf("
|  for (k=k1; k<=k2; k++) {|
|    if (k<=13) eta=1.308;|
|    else if (k<=32) eta=1.2609;|
|    else eta=1.12766;|
|    bestn0=0; bestn=0; bestc=1.0e40;|
|    for (n0=1; n0<=2*k; n0++) {|
|      calcparm();|
|      n2 = (long)  (kk*2.5*logk) + 50;|
|      for (n=k+1;n<=n2;n++) {|
|       if (exponent(n,&c,0)==0) {|
|         if (c<bestc) { bestc=c; bestn=n; bestn0=n0; }|
|       }|
|      }|
|      if (bestn<1) bestc=-99.99;|
|    }  /* for n0 */|
|    if (bestn0<1) CC[k]=-99.99;|
|    else {|
|#ifdef DEBUG|
|      for (n=bestn-25; n<=bestn+5 ; n++) exponent(n,&c,1);|
|#endif|
|      nn[k]=bestn; CC[k]=bestc+0.00005; n00[k]=bestn0;|
|      printf("k=
|          k,k-1,k,bestn0,bestn,bestc+0.000005); |
|    }|
|  }    /* for k */|
|  nn[k2+1]=999; CC[k2+1]=99.999;  /* print in TeX tabular format */|
|  i = (k1+k2)/2-k1+1;|
|  for (j=k1; j<=(k1+k2)/2; j++) {|
|    if (j==4) printf("&& 
|    else printf("&& 
|    printf("--
|    printf(" 
|      n00[j+i],nn[j+i],CC[j+i]);|
|  }|
|}|
\endlines

}
\enddocument
\bye